\documentclass{article}
\usepackage[a4paper]{geometry}
\usepackage[english]{babel}
\usepackage[utf8]{inputenc}
\usepackage{amsmath}
\usepackage[toc,page]{appendix}
\usepackage{amsthm}
\usepackage[small]{titlesec}
\usepackage{mathtools}
\usepackage{graphicx}
\usepackage{mathrsfs}
\usepackage{amsfonts} 
\usepackage{amssymb}
\usepackage{xcolor}
\usepackage{esint}
\usepackage{relsize}
\usepackage{comment}
\usepackage{bigints}
\usepackage{dsfont}
\usepackage{geometry}
\usepackage{float}
\usepackage{bbm}

\numberwithin{equation}{section}
\theoremstyle{plain}
\newtheorem{thm}{Theorem}[section]

\newtheorem{lem}[thm]{Lemma}

\newtheorem{mainresult}[thm]{Main Result}
\theoremstyle{definition}
\newtheorem{defn}{Definition}
\newtheorem{remark}{Remark}
\newtheorem{assumption}{Assumption}
\newenvironment{proofof}[1][Proof]{\noindent \textit{#1.} }{\ \qed}

\newcommand{\R}{\mathbb{R}^3}
\newcommand{\Rp}{\mathbb{R}^2}

\newcommand{\Ll}{\mathscr{L}}
\newcommand{\Rc}{\mathbb{R}^3}
\newcommand{\M}{\mathscr{M}_{+}\left(\Rc \right)}

\newcommand{\vp}{\varphi}
\newcommand{\Rl}{\mathbb{R}}
\newcommand{\Ii}{\int_0^{+\infty}}
\newcommand{\tx}{\tilde \xi}
\newcommand{\Su}{\sum_{n=0}^{+\infty}}

\theoremstyle{remark}

\title{Asymptotic Velocity Profiles for Homoenergetic Rayleigh–Boltzmann Flows under Dominant Shear}
\author{Nicola Miele, Alessia Nota, Juan J. L. Vel\'azquez}
\begin{document}

\maketitle

\begin{abstract}
In this paper, we study a particular class of solutions to the Rayleigh--Boltzmann equation, known in the nonlinear setting as \emph{homoenergetic solutions}. These solutions take the form $ g(x, v, t) = f(v - L(t)x, t),$ where the matrix $L(t)$ represents a shear flow deformation. We began our analysis in \cite{MNV}, where we rigorously proved the existence of a stationary non-equilibrium solution and established different behaviours of the solutions depending on the size of the shear parameter, for cut-off collision kernels with homogeneity parameter $0 \leq \gamma < 1$, thus including Maxwell molecules and hard potentials. In the present work, we focus on the regime in which the deformation term dominates the collision term for large times (hyperbolic-dominated regime). This scenario occurs for collision kernels with $\gamma < 0$; in particular, we focus on the range $\gamma \in (-1, 0)$. In this regime, it is challenging to obtain a clear and direct description of the long-time asymptotic behaviour of the solutions. Here we present a formal analysis of the velocity distribution's long-time asymptotics and derive for the first time the explicit form of the corresponding asymptotic profile. We also discuss the different asymptotic behaviour expected in the case of homogeneity $\gamma < -1$. In addition, we provide a probabilistic interpretation involving a stochastic process combining collisions with shear flow. The tagged particle velocity $\{v(t)\}_{t\geq 0}$ is a Markov process that arises from the combination of free flights in a shear flow along with random jumps caused by collisions.

\end{abstract}

\tableofcontents

\bigskip

\section{Introduction}
In this paper, we continue our study, initiated in  \cite{MNV}, of    homoenergetic solutions for the Rayleigh–Boltzmann Equation, which represent a particular class of non-equilibrium solutions. More specifically, we examine the Rayleigh–Boltzmann equation
\begin{align}
\partial_{t}g+v\cdot \partial_{x}g  &  =(\Ll g)\left(  v\right), 
	\;  \quad g=g\left(x,v,t\right): \R\times \R \times \Rl_+ \rightarrow  \mathbb{R}^+ \, ,
\label{A0_0}
    \\
	(\Ll g)\left(  v\right)   &  =\int_{\mathbb{R}^{3}}dv_{\ast}\int_{S^{2}%
	} \, B\left( n \cdot  \omega, |v-v_*| \right)
	\left(  M_*'g_{\ast}^{\prime}-M_*g\right) d\omega .   \label{A0_1} 
\end{align} 
which models the dynamics of a Rayleigh gas at the mesoscopic level. This model consists of a tracer particle interacting with a background gas at thermal equilibrium in the entire  $\mathbb{R}^3$ (e.g., see \cite{Spohn80}). The derivation of \eqref{A0_0}-\eqref{A0_1} from the Rayleigh gas model in the case of hard-sphere interactions and in a low-density limit, has been established in \cite{BLLS} (see also \cite{Spohn80} and \cite{LNW, MST, NVW}).

{
In \eqref{A0_1}, we use the conventional notation in kinetic theory 
$g=g\left(t, x, v\right)  , \ g_{\ast
}=g\left( t, x, v_{\ast}\right)  , \ g^{\prime}=g\left( t, x, v^{\prime}\right) , \ g_{\ast}^{\prime}=g\left( t, x,  v_{\ast}^{\prime}\right)  $ as well as $M_*=M(v_*), M'_*=M(v_*')$ where $M$ denotes the normalized Maxwellian density function $M(v_*)=\left(\frac{1}{2 \pi}\right)^{\frac{3}{2}}e^{-\frac{|v_*|^2}{2}}$. For simplicity, we normalize the temperature $T=\beta^{-1}>0$ to one, i.e., $T=1$. } 
Moreover,  $n=n\left(v,v_{\ast}\right)  =\frac{\left(  v-v_{\ast}\right)}{\left\vert v- v_{\ast}\right\vert }$ where the pair $(v,v_{\ast})$ identifies the velocity of the tagged particle and the velocity of one background particle before the collision respectively, while  $(v^{\prime},v_{\ast}^{\prime})$ consists of the corresponding pair of post-collisional velocities obtained by the collision rule 
\begin{equation}\label{eq:collisonRule}
    \begin{cases}
        v'=v-((v-v_*) \cdot \omega)\omega, & \\
        v_*'=v_*+((v-v_*) \cdot \omega)\omega.
    \end{cases}
\end{equation}
The vector  $\omega=\omega(v,V) \in S^2$ bisects the angle between the
incoming relative velocity $V=v-v_{\ast}$ and the outgoing relative velocity $V^{\prime}=v'-v_{\ast}^{\prime}$. The collision kernel $B\left(n \cdot \omega,|v|\right)  $ is a non-negative function which is proportional to the cross section for the scattering problem associated to the collision between two particles. 

We remark that the linear collision operator  $\Ll$ appearing in \eqref{A0_1} is such that $\Ll(g)=Q(M,g)$  where $Q$ denotes the standard nonlinear Boltzmann collision operator. Furthermore, we notice that, in the Rayleigh gas model, the energy of the tagged particle is not constant and changes due to the collisions with the background particles.

We will assume that the collision kernel is homogeneous  with respect to $|v|$ and denote its homogeneity by  $\gamma \in \mathbb{R}$, namely
\begin{equation}\label{eq:homB}
    B(n \cdot \omega,\lambda|v|)=\lambda^{\gamma}B(n \cdot \omega, |v|) \quad \text{for all $\lambda>0,\, v \in \R$ and for some $\gamma \in \Rl$.}
\end{equation} 
The homogeneity $\gamma$ is related to the properties of the interaction potential between particles. We recall that in the standard literature in kinetic theory, interaction potentials with the form $\phi(r)=\frac 1 {r^{s-1}}$ have  homogeneity  $\gamma=\frac{s-5}{s-1}$ in $d=3$ for the kernel $B$. Indeed, for this class of potentials, the collision kernels can be decomposed (e.g., cf.~\cite{V02}) as 
\begin{equation}\label{eq:decB}
    B(n \cdot \omega,|v|)=b(n\cdot \omega)|v|^{\gamma} 
 \quad b \in L^{\infty}([-1,1]),\, b   \geq 0, \, \gamma \in \mathbb{R}.
\end{equation}
since $n \cdot \omega = \cos \theta$ with $\theta 
\in [0,\pi]$ the scattering angle associated to the collision dynamics \eqref{eq:collisonRule}. 

The primary focus of this paper 
 is a specific class of solutions to the  Rayleigh–Boltzmann Equation known as \emph{homoenergetic solutions}, which provide insights into the properties of Boltzmann gases in open systems. Originally introduced for the nonlinear Boltzmann equation by Galkin \cite{Galkin1} and Truesdell \cite{T}, these solutions have since been considered both in the physical and in the mathematical literature. 
For instance, we refer to the books \cite{garzo} and \cite{TM}.  From the mathematical point of view the well-posedness of homoenergetic solutions to the nonlinear Boltzmann equation has been first addressed in \cite{CercArchive}, \cite{Cerc2000}
and more recently in \cite{JNV1}. The motivation to consider homoenergetic solutions to the Boltzmann equation is also justified by the existence of an invariant manifold of solutions to the classical molecular dynamics satisfying a number of symmetry properties, as constructed in  \cite{DJ, DJ1} and also in \cite{JQW}.

It is important to mention that the Boltzmann equation for homoenergetic flows can be thought of as a model for dilute gases having proper boundary conditions at infinity which yield mechanical deformations on the system. Hence, the solutions describing the long-time asymptotics of the problem cannot be expected, in general, to satisfy the detailed balance condition. In this regard, the model can be thought of as describing an open system and therefore we do not have thermalization to an equilibrium state. 
The long-time behaviour of homoenergetic solutions to the nonlinear Boltzmann equation is a novel topic which 
has been recently studied in a series of papers  \cite{BNV,DL1, DL2, JNV1,JNV2,JNV3,K1,K2,MT,NV}. We also refer to  \cite{CCDL}, which addresses the case of inelastic interactions.  The main goal of this paper is the study of the long-time behaviour of specific homoenergetic solutions for the Rayleigh-Boltzmann equation, that have not yet been  considered in the mathematical literature. As discussed further in this introduction, this particular class of solutions is worth studying in the linear setting as well, since they illustrate 
 interesting new features from both the mathematical and physical perspective. Due to the lack of detailed balance there are fluxes between the different states of the system, in the sense that there is no balance between the direct and the reverse collisions. The detailed description of these fluxes hence represents one of the main challenges addressed in this paper.  
 
To be precise, homoenergetic solutions are solutions to \eqref{A0_0}-\eqref{A0_1} of the form
\begin{equation}
	g\left(x,v,t\right)  =f\left(w,t \right)  \text{   with $w=v-\xi\left(
	t,x\right)$}.  \label{B1_0}
\end{equation}
As pointed out in \cite{MNV}, the transformation given by \eqref{B1_0} acts on both the tagged particle and the background particles. This implies that also  the Maxwellian distribution for the velocities of the background particles  is modified by the field $\xi(x,t)$. In fact, as a consequence of  \eqref{B1_0},   the background velocities are transformed as  $w_* = v_*-\xi(x,t) $, from which $M(w_\ast)$ is then transformed into
\begin{equation}\label{B1_1}
\widetilde{M}(w_\ast)=(2\pi)^{- \frac{3}{2}}e^{-\frac{|w_*|^2}{2}} .
\end{equation}
This requirement is what makes the homoenergetic ansatz  physically consistent in the linear setting as well; we refer to the appendix of \cite{MNV} for further details. 
Under suitable integrability conditions, every solution to 
the nonlinear Boltzmann equation of the form \eqref{B1_0} yields only time-dependent density $ \rho
	\left(  x,t\right)  =\rho\left(  t\right)$ and internal energy $\varepsilon\left(  x,t \right)  =\varepsilon\left(  t\right)$ with
\begin{equation}
\rho\left(t\right)=\int_{\mathbb{R}^3} f(w,t) dw , \quad  \varepsilon\left(t\right)= \frac 1 {2\rho\left(t\right)}
   \int_{\mathbb{R}^{3}}f\left(  w,t\right)  \vert w\vert^{2} dw  \,.\label{eq:Hom1}
\end{equation}

Nonetheless, for solutions of \eqref{A0_0} the only relevant macroscopic observable is the density, and thus for homoenergetic solutions of \eqref{A0_0} the density satisfies $\rho(x,t)=\rho(t)$. From direct computations we obtain that in order to have solutions of \eqref{A0_0} of the form \eqref{B1_0} for a sufficiently large class of initial data we must have
\begin{equation}
	\frac{\partial\xi_j}{\partial x_k} \; \; \text{independent of $x$ and \; 
	$\partial_{t}\xi+\xi\cdot\nabla\xi=0$} .\label{B2_0}
\end{equation}
The first condition yields that $\xi$ is an affine function of $x$, namely $\xi(x,t)=L(t)x+b(t)$ where $L\left(  t\right)  \in \mathrm{M}_3\left(  \mathbb{R}\right)  $ is a
$3\times3$ real matrix and $b(t)\in \mathbb{R}^3$. However, for simplicity, we will restrict our attention to the case in which $\xi$ is a linear function of $x$ (see \cite{JNV3}, Remark 2.1 for more details), whence  
\begin{equation}\label{eq:xitv}
  \xi (x,t) = L (t) x \, .  
\end{equation} 
On the other hand, the second equation implies  
\begin{equation}
	\frac{d L\left(  t\right)  }{d t}+\left(  L\left(  t\right)  \right)
	^{2}=0, \quad L (0) = A \label{B3_0}%
\end{equation} 
where we have added an initial condition $A\in \mathrm{M}_3\left(  \mathbb{R}\right)$. 
Classical results of ODE theory provide a unique solution of \eqref{B3_0}, which reads
\begin{equation}
	L\left(  t\right)  =\left(  I+tA\right)  ^{-1}A=A\left(  I+tA\right)
	^{-1} \label{B7_0}%
\end{equation}
and that is defined on the maximal interval of existence $[0,a)$ on which $\det \left( I+tA \right)>0 $. Then, combining \eqref{B1_0} and \eqref{eq:xitv}, equation \eqref{A0_0} becomes 
\begin{equation}\label{eq:homBol}
    \partial_t f-L(t)w \cdot \partial_w f = \Ll(f), \quad f=f(w,t), 
\end{equation}
where $\Ll$ is the collision operator  defined in \eqref{A0_1}. In the following we will denote the term  $L(t)w \cdot \partial_w f$ as the \emph{hyperbolic term}.   
We are interested in the long-time asymptotics of the 
function $\xi\left(  x,t \right)  =L\left(  t\right)  x = A\left(  I+tA\right)^{-1} x$.  
In order to study its asymptotic behaviour, the key idea is to analyse the form of the matrix $L(t)$ as $t \to  + \infty$ in a particular orthonormal basis using the Jordan canonical form for real $3\times 3$ matrices. Such a classification has been obtained in \cite{JNV1} Theorem $3.1$, where all the possible long-time asymptotics of the matrix $L(t)$ has been obtained under the assumption $\det(I+tA)>0$ for all $t \geq 0$. We also remark  that the different asymptotics of the matrix $L(t)$ describe the mechanical deformation acting on the gas. In this work we will focus on the case of \emph{simple shear deformation}, i.e. when the matrix $L(t)$ has the form
\begin{equation}\label{eq:ShearMat}
L = \left( \begin{array}{ccc} 0 &  -K & 0 \\
	0 & 0 & 0  \\
    0 & 0 & 0 
\end{array} \right), \quad K>0\,. 
 \end{equation}
With this choice, the linear Boltzmann equation for homoenergetic flows, as given by \eqref{eq:homBol}, becomes 
\begin{equation}\label{eq:Bshear}
	\partial_tf {+}K w_2 \partial_{w_1}f=\Ll(f).
\end{equation}
Notice that we adopt an unconventional convention by assigning the shear parameter a negative sign in order to ensure that the support of the function lies in the half-plane $\xi_1 > 0$, thereby simplifying the computations. 
Such a deformation is paradigmatic, since it preserves the mass conservation property of \eqref{A0_0}, while at the same time its lack of symmetry provides new and interesting behaviours for the solution of \eqref{eq:homBol}.

Our main interest is to study the long-time behaviour of the solutions to \eqref{eq:Bshear}. Even though the model that we are considering is linear its analysis  has both mathematical and physical interest. In fact, the interplay between the external  deformation due to the shear and the collisional dynamics, given by the interaction of the tagged particle with the background thermal bath, produces new features. 
In order to study the asymptotic behaviour of the solutions to \eqref{eq:Bshear}, we start by looking at the physical dimensions of the three terms appearing in the equation, namely
\begin{equation}\label{eq:eqMagitude}
    \frac{[f]}{[t]},[f],[w]^{\gamma}[f]
\end{equation} 
where $[w]$ is the order of magnitude of $w$, $\gamma$ is the homogeneity of the collision kernel $B$,  as in \eqref{eq:homB}, and $[f]$ is the order of magnitude of $f$. The long-time behaviour of the solutions of \eqref{eq:Bshear} is strongly  dependent on $\gamma$. As $t \rightarrow + \infty$, one expects three different possibilities: either the collision term
is much larger than the transport term, 
or the transport term is much larger than the other two, or the transport term and the collision term are of the same order of magnitude. According to the terminology introduced in \cite{JNV1}-\cite{JNV3} for the nonlinear Boltzmann equation one refers to the first case as \emph{collision-dominated regime}, taking place for $\gamma>0$, and to the second case as \emph{hyperbolic-dominated regime}, taking place for $\gamma<0$.  The case in which the three terms in \eqref{eq:Bshear} have the same order of magnitude is obtained for collision kernels describing Maxwell molecules interactions between the gas particles; this corresponds to homogeneity $\gamma=0$.  
Notice that the values of $\gamma$ characterizing the different regimes 
are related to the  specific form of the deformation matrix $L(t)$ acting on the system.  In \cite{MNV} we presented a rigorous analysis of the solutions to \eqref{eq:Bshear} in the case of cut-off collision kernels with homogeneity parameter $0\leq \gamma<1$, i.e. in the collision-dominated  and in the balance regime. 
We established a well-posedness result for this class of solutions in the space of non-negative Radon measures as well as the existence of a stationary non-equilibrium solution. Moreover, for the case of Maxwell molecules, i.e.,  when $\gamma=0$, we obtained information on the behaviour of the
solution also for large values of the shear parameter, hence establishing  different behaviours of the solutions for small and large values of the shear parameter. 

In this paper, we will focus on the long-time behaviour of solutions to \eqref{eq:Bshear} in the hyperbolic-dominated regime, i.e., when $\gamma<0$. We recall that collision kernels with homogeneity parameter $\gamma<0$ correspond to the power law interaction potentials $\phi(r)=\frac 1 {r^{s-1}}$ with $2<s<5$ in $d=3$, usually referred to as soft potentials. More specifically, we concentrate on the case of $-1<\gamma<0$. The distinction between the cases  $\gamma<-1$ and $-1 <\gamma<0$ arises from the expectation of significantly different  long-time behaviours in these two regimes. Indeed, when $\gamma<-1$, the collision rate is expected to be exponentially decreasing
in time and hence negligible as $t\to  + \infty$, as will be explained in Section \ref{sec:frozenCollision}.  In this case, the effect of the
collisions becomes frozen for large times, which motivates the name \emph{frozen collision} to refer to this regime. More precisely, the tagged particle, whose dynamics is described by \eqref{eq:Bshear}, experiences a finite number of collisions with probability $1$. On the other hand, in the regime considered here, i.e.,   $-1<\gamma<0$,  
the collision rate is
small but  not negligible  and still plays a significant role in determining the asymptotics of the solutions to \eqref{eq:Bshear}, the main effect being that collisions take place with increasingly larger mean free paths. The tagged particle experiences an infinite number of collisions with probability $1$. The description of the long
time asymptotics for the particle distribution in this case is a challenging problem, as already emphasized in \cite{JNV1,JNV3}. Indeed, in \cite{JNV3}, the authors considered a simplified model that incorporates the effects of shear and a basic mechanism for collisions, and obtained estimates suggesting  
a complex behaviour for the particles distribution, which differs greatly from a
self-similar distribution or a Maxwellian distribution.
In fact, in this regime, we will show that there is a not a single time-dependent velocity scale $\vert w\vert \approx t^{\alpha}$
, that characterizes the scale of velocities in which most of the particles of the system as well as the energy is contained at a given time, as usually happens in the case of self-similar solutions.

In this paper, we provide a contribution in this direction. More precisely, using formal asymptotics, we obtain a characterization of the asymptotic profile of the solution to \eqref{eq:Bshear} for long times. To the best of our knowledge, this result is new and provides a first understanding of the detailed particle distribution for large times in the case in which the hyperbolic terms are much larger than the collision 
terms.   

Setting for simplicity $a=|\gamma|$, we can state the following formal theorem on the long-time behaviour of the solutions to \eqref{eq:Bshear}, which is the main result of this paper.

\begin{mainresult}[Asymptotics]\label{thm:ConvergenceF} 
    Let $f \in C([0,+\infty),\M)$ be 
    a solution of \eqref{eq:Bshear} 
    and let $a \in (0,1)$.  Set
$$f(w,t)= \frac{1}{t^{\frac{3}{a}-2}} F\left(\xi, \tau \right),$$ where  $\xi=(\xi_1, \xi_2, \xi_3)\in  \mathbb{R}^3,$ with $\xi_1=\frac{w_1}{t^{\frac{1}{a}}}, \ \ \xi_i=\frac{w_i}{t^{\frac{1}{a}-1}},$ for $i=2,3,$ and  $ \tau= \log t$. 
Then we can construct, using matched asymptotic expansions, a solution with the following asymptotic behavior
\begin{equation}\notag
F(\xi,\tau) \sim F\left(0,\frac{(\xi_2)^{\frac{1}{a}}}{(\xi_2-\xi_1)^{\frac{1}{a}-1}},\frac{(\xi_2)^{\frac{1}{a}-1}\xi_3}{(\xi_2-\xi_1)^{\frac{1}{a}-1}},\tau\right)\exp \left(-\frac{4 \pi\xi_1^{1-a}}{(1-a)\xi_2}\right)\left(1-\frac{\xi_1}{\xi_2}\right)^{-\left(\frac{3}{a}-2\right)}  \ \ \text{as} \ \ \tau\to + \infty 
\end{equation}
with 
\begin{align}
   F(0,\xi_2,\xi_3,\tau) \sim \frac{8aM_0}{(1-a)\log\left(\frac{1}{1-a}\right)\xi_2(|2\tx|)^{2+a}e^{-a\tau}\tau}  \int_1^{+\infty}R\left( \frac{1}{z}\right)\exp \left(-\frac{4 \pi z^{-a}}{1-a} \left(\frac{2|\tx|}{e^{\tau}}\right)^{-a}\right)\frac{dz}{z^2} \notag
\end{align}
where $\tx=(\xi_2,\xi_3)$, $\xi_1,\xi_2 \geq 0$, $M_0>0$ is the initial mass, and   the function $R:[0,1] \rightarrow \Rl$ is (see \eqref{eq:R} below) 
\begin{equation*}
    R(s)=\frac{1}{\sqrt{2}}\frac{s^a}{\sqrt{1-s^2}}\left(\sqrt{1+\sqrt{1-s^2}}+\sqrt{1-\sqrt{1-s^2}}\right).
\end{equation*}
\end{mainresult}
\smallskip

\begin{remark}
    Notice that the result we are proving in this paper is an asymptotic result in which the only information that we use about the collision kernel is the asymptotic behaviour in the velocity variable. Specifically, the fact that the kinetic part of the kernel behaves as a power law. More precisely, as it will be discussed in Section \ref{ssec:AsymModel} below, we will assume by simplicity that $B(n \cdot \omega, |w|)=|w|^{\gamma}$ with $\gamma \in (-1, 0)\,$. Observe that the singularities of the collision kernel could  affect the well-posedness of the solutions of \eqref{eq:homBol}.  Actually, the asymptotics we obtain in this paper can be expected to hold for kernels of the form $B(n \cdot \omega, |w|)=(1+|w|)^{\gamma}$ which do not have singularities in the region of small velocities. 
\end{remark}

The plan of the paper is the following. In Section \ref{ssec:AsymModel} we introduce a formal approximation for the adjoint equation of \eqref{eq:Bshear} that holds true for large velocities. As discussed in Section \ref{sec:heuristics},  from this asymptotic model it is possible to deduce a probabilistic interpretation of the microscopic process underlying the dynamics.  This is achieved by computing the typical \emph{flight probabilities} and \emph{jump probabilities} of a stochastic process associated with the approximate model. The details of this computations are presented in Section \ref{sec:probabilities}. 
   As emphasized in Section \ref{sec:heuristics},  the main interest in the probabilistic description is that it provides
 the correct scaling under which it is feasible to identify the asymptotic behaviour of the solutions to \eqref{eq:Bshear}. In particular, from this analysis, we obtain that the expected time 
 between collision increases with time, and the 
 probability distribution of the tagged particle turns out to be strongly anisotropic, very
elongated along the direction of $v_1$, as a consequence of the shear deformation \eqref{eq:ShearMat}. 

\begin{figure}[H]
    \centering
    \includegraphics[width=0.5\linewidth]{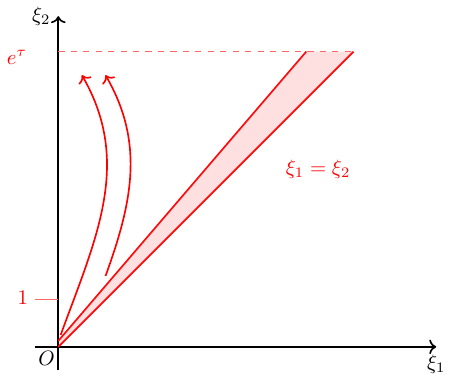}
    \label{fig}
    \caption{The mass $M$ of $G$ lives in $\xi_2 \in [1,e^{\tau}], \xi_1 \approx \xi_2$. The arrows represent the effect of the collisions.}
\end{figure}

Section \ref{sec:longtime} constitutes the main body of this paper, in which we obtain  the long-time asymptotics for the solution of \eqref{eq:Bshear} at the formal level. To this aim, we consider a suitable rescaling of the distribution function. Recalling that for simplicity we set $a=|\gamma|$, we look for solution $f$ with the following form:
   $$f(w,t)= \frac{1}{t^{\frac{3}{a}-2}} F\left(\xi_1, \xi_2, \xi_3, \tau \right), \quad   \xi_1=\frac{w_1}{t^{\frac{1}{a}}}, \, \xi_2=\frac{w_2}{t^{\frac{1}{a}-1}}, \, \xi_3=\frac{w_3}{t^{\frac{1}{a}-1}},\quad \tau= \log t\, .$$  A crucial point upon which we will rely is that the equation satisfied by the rescaled distribution $F$ can be conveniently reformulated as a suitable boundary value problem.    Essentially, this involves a simplified equation that retains only the loss term of the collision operator, while the gain term is adjusted to represent the boundary condition. 
   More precisely, we will rely on the analysis of a reduced distribution function 
\begin{equation}\label{def:G_2d}\notag
    G(\xi_1,\xi_2, \tau) =\int_{\Rl} F(\xi_1,\xi_2,\xi_3,\tau)\, d \xi_3\, 
\end{equation}
which satisfies a closed equation. Specifically, we can obtain that the function $G(0,\xi_2, \tau)$ satisfies approximately a delay equation for which we can obtain an explicit asymptotic solution. A visual  
description  illustrating   where the mass of $G$ is concentrated is presented in Figure \ref{fig}. 

This analysis allows to determine the asymptotic profile of $F$. The most relevant feature is that the mass of $F$ is concentrated in the region $\xi_2,\xi_3 \in [1,e^{\tau}],$ $ \xi_1 \approx \xi_2$. 
In other words, the mass of the solution is concentrated in an increasing number of dyadic intervals as $\tau\to \infty$. 
As explained in Section \ref{sec:noSelf}, we could have expected $F$ to behave in a self-similar manner for large time, thus having its mass concentrated in the region $\xi\approx 1$. Actually, this is not the case, as the shear deformation elongates the profile of $F$ in the direction of  $\xi_2$.
Consequently, this spreading of mass over various scales of $\xi_2$ disrupts the self-similar structure, leading to its inconsistency. 

We conclude the paper with a discussion on the expected behaviour of the solutions to \eqref{eq:Bshear} in the frozen collision regime, i.e., when $\gamma <-1$, which is presented in Section \ref{sec:frozenCollision}. In this regime, we expect the collision term in \eqref{eq:eqMagitude} to decay very fast in time. Therefore,  the tagged particle experiences only a finite number of collisions and the dynamics of its probability distribution is driven only by the shear term. Using a proper rescaling for the distribution function, which differs from the one considered in the complementary regime $-1<\gamma<0$ discussed above, and which reads
\begin{equation*}
    f(w,t)=\frac{1}{t}F(\xi,t),\quad \xi_1= \frac{w_1}{t}, \, \xi_2=w_2,\, \xi_3=w_3, \quad \tau= \log t\, ,
\end{equation*}
we conjecture that $F$ converges to the stationary distribution $\delta(\xi_1-K\xi_2)F_{\infty}(\xi)$, with $F_{\infty}$ only depending on the initial datum, as $\tau \rightarrow + \infty$, (see \textbf{Main Result} 
\ref{thm:ConvergenceFrozen} below). 

It is worth noticing that, so far, we have excluded the threshold case $\gamma=-1$ from both the frozen collision regime and the case studied in this paper. More precisely, when $\gamma=-1$, we expect the velocities $w$ to increase linearly as $t \rightarrow+\infty$, thus exhibiting an asymptotic behaviour analogous to the one derived in this paper. Nonetheless, for $\gamma=-1$ we expect to have additional logarithmic corrections in the derivation of the long-time behaviour, even at the formal level; therefore this particular case requires a separate and more careful treatment that will be the object of further studies by the authors.  In conclusion, we emphasize that the results obtained in this paper,  albeit formal, reveal new   
interesting behaviours for the solution of the Rayleigh-Boltzmann equation under an external force field, whose asymptotics greatly differs from classical  results established in the literature.

\section{The approximate model}\label{ssec:AsymModel}

In this section, we introduce an approximate model for the adjoint formulation of the Rayleigh–Boltzmann equation \eqref{eq:Bshear}. Specifically, we consider an approximation  $\mathcal{L}$ of the adjoint of the collision operator $\Ll$ for large velocities, i.e., 
 $|w|\gg1$. 
The resulting equation will be the object of study of this paper. 

As mentioned earlier in the introduction, in the kinetic theory literature it is customary to consider collision kernels of the form \eqref{eq:decB}, namely $
    B(n \cdot \omega, |w|)=b(n \cdot \omega)|w|^{\gamma},$ with $\gamma \in \Rl.$
These kernels are homogeneous with respect to $|w|$. 
In this paper we will restrict 
our attention to the case of \emph{soft} interaction \emph{potentials} corresponding to homogeneity $\gamma<0$. More precisely, here we outline the  
precise assumptions that will be made on the collision kernel. 
\begin{assumption}\label{ass:B}
The collision kernel $B$ has  the form \eqref{eq:decB} and satisfies  
\begin{equation}\label{eq:assB}
   B(n \cdot \omega, |w|)=|w|^{\gamma} \quad \text{with}\quad \gamma \in (-1, 0)\, . 
\end{equation} 
\end{assumption}
Note that, for simplicity,  we have assumed the angular part  of the collision kernel to be constant and, without loss of generality, we set this constant equal to $1$, i.e.,  $b(n\cdot \omega)=1$. Additionally, we remark that we are only considering a specific range of negative values for $\gamma$. As mentioned in the introduction, to simplify notation, we will denote $a:=|\gamma|$ in the rest of this paper. Consequently, Assumption \ref{ass:B} becomes $B(n \cdot \omega, |w|)=|w|^{-a}$.  It seems possible to adapt the arguments of this paper to the more general case of non-constant $b( n \cdot \omega)$.  
However, for clarity, we will restrict ourselves to the case where $b(n \cdot \omega)=1$.

We first introduce the adjoint equation of \eqref{eq:Bshear}
\begin{align}
    \label{eq:homAdj}
    &\partial_t \vp  - Kw_2 \partial_{w_1}\vp = \Ll^*(\vp), \quad \vp=\vp(w,t) \\& 
    \quad \Ll^*(\vp) = \int_{\R}dw_* \int_{S^2} B(n \cdot \omega, |w-w_*|)(\vp'-\vp) \,M_*\, d \omega \, , \nonumber 
\end{align}
with $M_*=M(w_*)$. Using  Assumption \ref{ass:B}  with $a=|\gamma|$ , we have that $\Ll^*$  reads
\begin{equation}
\label{eq:adjL}
    \Ll^*(\vp) =  \int_{\R}dw_* \int_{S^2}   |w-w_*|^{-a} (\vp'-\vp)\, M_* d \omega \, .
\end{equation}
As it will be discussed in detail in Section  \ref{sec:longtime}, equation \eqref{eq:homAdj} is motivated by the weak formulation \eqref{eq:weakf} for \eqref{eq:Bshear} (see Definition \ref{def:weakSol} in Section \ref{sec:longtime}).  
At this point, we approximate the reference distribution   in \eqref{eq:adjL} as  $M(w_*)\sim \delta(w_*)$, where $\delta(w-w_0)$ denotes the Dirac delta measure concentrated at $w_0\in \mathbb{R}^3$. Thus, the collision operator \eqref{eq:adjL} becomes
\begin{equation}\label{eq:apprL}
    \Ll^*(\varphi)(w)=|w|^{-a} \int_{S^2} (\varphi(w-(w \cdot \omega)\omega)-\varphi(w)) d \omega \, \,  \end{equation} 
The relation between pre- and post-collisional velocities  given by \eqref{eq:collisonRule} simplifies to the following asymptotic collision rule  
\begin{equation}\label{eq:approxColl}
    w'= w-(w \cdot \omega) \omega = \mathrm{P}^{\perp}_{\omega}w,
\end{equation}
where $\mathrm{P}^{\perp}_{\omega}w$ denotes the orthogonal projection of $w$ onto $\omega\in S^2$. Roughly speaking, this approximation can be justified as follows. Suppose we take $R>0$ to be sufficiently large and of order one. 
Then, we have $\int_{|w_*| \leq R}M_* dw_* \approx 1$. If we consider velocities $w\in \mathbb{R}^3$ such that $|w|\gg1$,  it follows that $|w-w_*| \approx |w|$,  which implies that $|w_*|\approx 0$, given that the velocity $w$ is only slightly modified during collisions. Consequently, we have 
$M(w_*) \sim \delta(w_*)$. Thus,  by using \eqref{eq:assB}, 
the adjoint collision operator becomes  \eqref{eq:apprL}. From \eqref{eq:homAdj}, we arrive at
\begin{equation}\label{eq:asymModel}
    \partial_t \vp-Kw_2 \partial_{w_1}\vp= |w|^{-a} \int_{S^2} (\varphi(w-(w \cdot \omega)\omega)-\varphi(w)) d \omega	.
\end{equation} 
that corresponds to  \eqref{forwardadj}  in Section \ref{sec:preliminary}. 
Moreover, without loss of generality,  we can rescale time and velocities in such a way that we obtain a shear parameter $K=1$. Hence, from now on, we will always assume $K=1$.

\section{Derivation of the characteristic scales for the free flights and the velocities}  
\label{sec:heuristics}

In this section, we provide some heuristic considerations that will be useful to study the long time behaviour of solutions to \eqref{eq:Bshear}. In order to do this, we will rely on the adjoint formulation of the problem, that is  \eqref{eq:homAdj}. More precisely, we will assume $K=1$ for simplicity and study the approximate model for \eqref{eq:homAdj} proposed in Section \ref{ssec:AsymModel} (cf.~\eqref{eq:asymModel}),  namely 
\begin{equation*}
    \partial_t \vp - w_2 \partial_{w_1}\vp= 
    |w|^{-a}\int_{S^2} (\varphi(w-(w \cdot \omega)\omega)-\varphi(w)) d \omega,	
\end{equation*} 
where we recall that we are denoting  $a=|\gamma|$  to simplify the notation. Notice that \eqref{eq:asymModel} is the evolution 
equation for the semigroup associated with 
the stochastic process that combines free flights with speed $Kw_2=w_2$ and jumps with rate $|w|^{-a}$.

We refer to Section \ref{sec:longtime} for the precise definitions. 
Here, we consider the stochastic process underlying \eqref{eq:asymModel}. This will enable us to determine the characteristic scales of the system, hence suggesting the correct time and velocity scales.

The stochastic process describing the dynamics of the tagged particle of the Rayleigh gas under shear is the result of the composition of two effects: first the particle experiences a free flight which takes place with different velocities depending on the value of $w_2$ due to the shear flow, then at random times it changes its direction according
to random kicks, i.e. it collides with another particle belonging to the background and thus the particle jumps to a new velocity. In order to study the stochastic process we first introduce  the collisionless flow $ T_{s} : \R \rightarrow \R $ defined by
\begin{equation}\label{def:T_t} 
       T_{s}(v_1,v_2,v_3)=(v_1+ s v_2,v_2,v_3), \quad v=(v_1,v_2,v_3)\in \mathbb{R}^3, \quad s \geq 0,
    \end{equation} 
which describes the action of the shear flow deformation on the velocity of the tagged particle. Then we define two sequences of random variables $(\xi^m)_m$, $(v^m)_m$ as
	\begin{align}
		\xi^m &=T_{t^m}v^m; \\
		v^{m+1} &= \mathrm{P}_{\omega}^{\perp}\xi^m
	\end{align}
with $t^0,v^0 \neq 0$ and where $(t^{m})_m$ are the times in between collisions for $m \in \mathbb{N}$. The reason for introducing two different sequences of velocity (random) variables is due to the coexisting effect of the shear and collision dynamics that we need to analyse separately.  More precisely, the random variable $\xi^m$ denotes the tagged particle velocity due to 
the collisionless shear flow  up to a time $t^{m}$ for any $m\in \mathbb{N}$, where the action of the shear is described by the map $T_{s}$. On the other hand, the random variable $v^m$ denotes tagged particle velocity resulting from an instantaneous collision with a background particle, as prescribed by the asymptotic collision rule \eqref{eq:approxColl}. Here $\omega$ is a random impact parameter which is assumed to be uniformly distributed on the sphere $S^2$. The times between two collisions $(t^{m})_m$ or, equivalently, the free flight times during which the particle is only subject to the action of the shear, are distributed according to an inhomogeneous Poisson process with rate $|T_{t^{m}}v^m|^{-a}$ whose probability density is given by
\begin{equation}\label{eq:time}
	    q(t^{m}|\, v^m) := |T_{t^{m}}v^m|^{-a} \exp\left(-\int_0^{t^{m}}|T_sv^m|^{-a} d s \right).
\end{equation}
Therefore the process as a whole can be schematized as
\begin{equation}
    v^0 \xrightarrow{t^0} \xi^0 \xrightarrow{} v^1 \xrightarrow{t^1} \xi_2 \ldots v^{m-1} \xrightarrow{t^m} \xi^m. 
\end{equation}
Moreover,  we observe that if we start from $v^0 \neq 0$ then $\xi^{m+1}$ is obtained as 
\begin{align}
    \xi^{m+1}=(T_{t^m}\mathrm{P}^{\perp}_{\omega^m})(\xi^m) = (T_{t^m}\mathrm{P}^{\perp}_{\omega^m})(T_{t^{m-1}}\mathrm{P}^{\perp}_{\omega^{m-1}})\ldots (T_{t^1}\mathrm{P}^{\perp}_{\omega^1})v_0 
\end{align}
where we have assumed $\omega^0=0$. Expanding the above equation using \eqref{def:T_t}, the fact that $T_{t+s}v=T_t(T_s(v))$ for $t,s \geq 0$ and the formula for $\mathrm{P}^{\perp}_{\omega}$ we obtain the following equation for $\xi^{m+1}$
\begin{equation}
    \xi^{m+1}=T_{\sum_{k=1}^m t^k}(v^0)+\omega^m\left(\prod_{k=1}^m(\omega^k \cdot T_{t^k}(\omega^{k-1}))(\omega^1 \cdot T_{t^0}(v^0))\right).
\end{equation}

To simplify the notation, in what follow we do not write explicitly the indexes of the random variables. It is possible to compute explicitly the \emph{flight probability}, i.e. the probability $p(\xi| \, v)$ that the tagged particle flies with an initial  velocity $v$ 
up to a time $t$
when the velocity is updated to a velocity $\xi$. More precisely, we have 
\begin{lem}[Flight probabilities] \label{lem:pFlight}
Let $p(\xi|\, v)$ denote the probability of flight of reaching a point $\xi=(\xi_1,\xi_2,\xi_3)\in \mathbb{R}^3$ in the velocity space from a point $v=(v_1,v_2,v_3)\in \mathbb{R}^3$. Then
\begin{align}\label{flightProb}
		p(\xi|\, v) & = 
        G(\xi_1| \, v)\delta(\xi_2-v_2)\delta(\xi_3-v_3)
	\end{align}
	where $G(\xi_1|v)$ is given by  
    \begin{align}
		G(\xi_1|\, v) := \frac{1}{|v_2|}\left(\xi_1^2+v_2^2+v_3^2 \right)^{-\frac{a}{2}}	&\Bigg[\exp\left(- \frac{1}{|v_2|}\int_{v_1}^{\xi_1}\ \left(\eta^2+v_2^2+v_3^2 \right)^{-\frac{a}{2}}d \eta \right)\chi_{\{v_2\geq 0\}}  \notag \\
		& + 
        \exp\left(-\frac{1}{|v_2|}\int_{\xi_1}^{v_1}(\eta^2+v_2^2+v_3^2 )^{-\frac{a}{2}}d \eta  \right)\chi_{\{v_2<0\}}\Bigg].\label{eq:defGprob}
	\end{align}
\end{lem}
When the velocity of the particle $v$ 
is updated to $\xi$, 
the particle immediately experiences a collision with a background obstacle that further updates its velocity to $v=\mathrm{P}_{\omega}^{\perp}\xi$. The probability of updating the velocity $\xi$ to $v$, i.e. the transition probability $p(\xi|v)$ from $v$ to $\xi$, can be computed explicitly. This is the content of the following lemma.   
\begin{lem}[Jump probabilities] \label{lem:pJump}
  Let $ p(v| \, \xi)$ denote the probability of jumping from a velocity $\xi\in \mathbb{R}^3 $ to $v\in \mathbb{R}^3$ due to an instantaneous collision with a background particle. Then, \begin{equation}\label{eq:defpJump}
            p(v| \, \xi) = H(v,\xi)\,\delta\left(\Big |v-\frac{\xi}{2}\Big |^2-\frac{|\xi|^2}{4} \right)
\end{equation}
where
\begin{equation}\label{eq:defpH}
            H(v,\xi) := \frac{1}{2\pi|\xi|} \frac{1}{|v-\xi|}.
\end{equation} 
\end{lem}
\begin{remark}
    In \eqref{eq:defpJump} we have obtained that the collisions are supported in the sphere $S_{\xi/2}(\xi/2)$. This is due to the approximation performed to obtain the asymptotic model in Section \ref{ssec:AsymModel} and, specifically, to the fact that in the approximate collision operator \eqref{eq:asymModel} the background particles have no dispersion of the velocities. For the original Rayleigh–Boltzmann equation \eqref{eq:homBol}, we have instead that the collisions are supported in a region that includes the sphere $S_{\xi/2}(\xi/2)$, but it has a thickness of order one and the support has some small elongation at infinity due the presence of the Maxwellian distribution $M_*$ in \eqref{A0_1}.
\end{remark}

It is worth noticing that, if we consider only the Markov process given by the sequence of random variables $(\xi^m)_m$, we can compute the transition probabilities $p(\xi^{m+1}| \, \xi^m)$ relying on the classical Chapman-Kolmogorov equations (e.g., \cite{Kall}). For simplicity, we set $\xi^{m+1}=\xi,\xi^m=\eta, v^m=v$. We have the following.  
\begin{lem}[Transition probabilities]\label{lem:transitionp}
    Let $p(\xi| \, \eta)$ denote the transition probability from $\eta\in \mathbb{R}^3$ to $\xi \in \mathbb{R}^3$. 
    Then,
\begin{align}
		p(\xi|\eta) = \frac{1}{\sqrt{\frac{|\eta|^2}{4}-\left(\xi_2-\frac{\eta_2}{2} \right)^2-\left(\xi_3- \frac{\eta_3}{2}\right)^2}} 
		\, \mathlarger{\sum}_{\ell =1, 2}\,\frac{G(\xi_1,v_1^{\ell},\xi_2,\xi_3)}{4 \pi \sqrt{(v_1^{\ell}-\eta_1)^2+(\xi_2-\eta_2)^2+(\xi_3-\eta_3)^2}} \label{transitionFinal}
	\end{align}
    where $G$ is defined as in \eqref{eq:defGprob} and with
    \begin{equation}
		v_1^{\ell}=\frac{\eta_1}{2}\pm\sqrt{\frac{|\eta|^2}{4}-\left(\xi_2-\frac{\eta_2}{2} \right)^2-\left(\xi_3- \frac{\eta_3}{2}\right)^2}, \quad \ell \in \{1,2\}.
	\end{equation}
    
    \end{lem}
    The proofs of the above Lemmas \ref{lem:pFlight}, \ref{lem:pJump} and \ref{lem:transitionp} are postponed to Section \ref{sec:probabilities}. 

\subsection{Characteristic scales of the system}\label{ssec:scalinA>1}
In order to study the asymptotic behaviour of the solution to \eqref{eq:asymModel} for $0 <a< 1$, we look for a convenient set of variables that encodes the scaling properties of the solution. From \eqref{def:T_t} we have that, after a flight in which the velocity of the tagged particle is affected by the sole effect of the shear deformation, the first component of its velocity $v^m_1 \approx v^m_2 t^m$ and $|v^m| \approx |v^m_1| $. Furthermore, recalling that the probability $q(t^m|v^{m-1})$ introduced in \eqref{eq:time} at time $t^m$ follows an inhomogeneous Poisson process we obtain 
\begin{equation}
    \frac{d\, q(t^m|v^{m-1})}{dt^m}=-|v^m|^{-a} q(t^m|v^{m-1}) \sim -|v^m_1|^{-a}q(t^m|v^{m-1}) \sim |v^m_2|^{-a}(t^m)^{-a}q(t^m|v^{m-1}).
\end{equation}
Integrating we get 
\begin{equation}
    q(t^m|v^{m-1}) \sim \exp \left(\int_0^{t^m} \frac 1 {|v_2^m|^{a} \,(s^m)^{a}} ds\right) .
\end{equation}
The right-hand side must be approximately $1$ when a jump takes place, since $q(t^m|v^{m-1})$ is a probability density. Therefore we obtain that the integral inside the exponential has to be approximately $1$ as well, that is
\begin{gather}\notag
   \int_0^{t^m} \frac 1 {|v_2^m|^{a} \,(s^m)^{a}}ds \approx 1
\end{gather}
which then gives
\begin{equation}\label{eq:approxtn}
    |v_2^m|^{-a}(t^m)^{1-a} \approx 1 \quad \text{and $t^m \sim |v_2^m|^{\frac{a}{1-a}}$}\, .
\end{equation}
 Note that \eqref{eq:approxtn} implies that the time between jumps increases together with $|v_2^m|$. From the definition of $T_{s}$ we have that $|v^{m+1}_2|\sim |v^m_1|$ and that $t^m \sim \frac{|v^m_1|}{|v^m_2|}$. Moreover
\begin{equation}
    |v^{m+1}_2|\sim |v^m_1|= |v^m_2|t^m=|v^m_2||v^m_2|^{\frac{a}{1-a}}=|v_2^m|^{\frac{1}{1-a}}.
\end{equation}
Defining then $\xi^m=\log(|v_2^m|)$ gives the recurrence law
\begin{equation}
    \xi^{m+1}=\frac{1}{1-a}\xi^m
\end{equation}
from which it follows that $\xi^m=C_0\left(\frac{1}{1-a}\right)^m$. This gives the following
\begin{align}
    |v_2^m| &= \exp\left( C_0 \left(\frac{1}{1-a}\right)^m\right); \\
    |v_1^m| &= \exp\left( C_0 \left(\frac{1}{1-a}\right)^{m+1}\right).
\end{align}
Hence, from the fact that $t^m=\displaystyle \frac{|v_1^m|}{|v^m_2|}$ we get
\begin{align}
    t^m & = \exp\left( C_0 \left(\frac{1}{1-a}\right)^m\left[\frac{1}{1-a}-1\right]\right) = 
    \exp\left( C_0 \left(\frac{1}{1-a}\right)^m\frac{a}{1-a}\right) \notag \\
    & = \exp\left( C_0 \left(\frac{1}{1-a}\right)^m\right)^{\frac{a}{1-a}}=|v_2^m|^{\frac{a}{1-a}} \label{eq:exptn}
\end{align}
and that
\begin{align}
     |v_1^m| = \exp\left( \frac{C_0}{1-a} \left(\frac{1}{1-a}\right)^m\right)=|v_2^m|^{\frac{1}{1-a}}.
\end{align}
We notice that \eqref{eq:exptn} implies that 
\begin{align}
    \frac{t^{m+1}}{t^m}= \exp \left(C_0 \left(\frac{1}{(1-a)^{m+1}}-\frac{1}{(1-a)^m}\right) \right) = \exp\left(\frac{C_0 a}{(1-a)^{m+1}}\right) \label{eq:increasingJumps}
\end{align}
which tends to $+\infty$ as $m \to +\infty$ since $\frac{1}{1-a}>1$. Moreover, we obtain
\begin{align}
    |v_2^m|&=(t^m)^{\frac{1-a}{a}}; \\
    |v_1^m|&=|v_2^m|^{\frac{1}{1-a}}=(t^m)^{\frac{1}{a}}
\end{align}
and at the end we get 
\begin{equation}
    |v_1^m|\sim (t^m)^{\frac{1}{a}}; \; |v_2^m|\sim \frac{|v_1^m|}{t^m}\sim (t^m)^{\frac{1-a}{a}}.
\end{equation}
It is worth noting that  such a scaling can be heuristically derived by directly examining  \eqref{eq:Bshear} and comparing the orders of magnitude, denoted by $[\cdot]$, of the three terms appearing in the equation, as previously discussed in the introduction. Assuming that the three terms in \eqref{eq:Bshear} are comparable we obtain
\begin{equation}
    \frac{[f]}{[t]}
    \approx [w_2]\frac{[f]}{[w_1]}\approx [w]^{-a}[f].
\end{equation}
Since $a>0$, balancing the second and third terms in \eqref{eq:homBol} we must have that $[w_2] \ll [w_1] $  as $\vert w\vert\to +\infty$, which implies $[w] \approx [w_1].$ To balance the first and third terms, we must have $[w_1] \approx t^{\frac{1}{a}}$, and balancing the first with the second term indicates $[w_2]=\frac{[w_1]}{t}$. We can expect $[w_3]\approx [w_2]$ since the effect of the collisions produces a similar change for both components  $w_2, \ w_3$. Therefore, to ensure mass conservation for $f$, we must fulfill the condition 
$[f][w_1][w_2][w_3]\approx 1$, which leads to $[f]\approx \frac{1}{t^{\frac{3}{a}-2}}$. 
In summary, we are led to consider the following scaling: 
\begin{equation}\label{def:scalingA<1}
   f(w,t)= \frac{1}{t^{\frac{3}{a}-2}} F\left(\xi_1, \xi_2, \xi_3, \tau \right), \quad   \xi_1=\frac{w_1}{t^{\frac{1}{a}}}, \, \xi_2=\frac{w_2}{t^{\frac{1}{a}-1}}, \, \xi_3=\frac{w_3}{t^{\frac{1}{a}-1}}, \quad \tau= \log t \ .  
\end{equation}

Relying on the change of variables  \eqref{def:scalingA<1} in the next section we will study the long-time behaviour of solutions to \eqref{eq:Bshear}. It is important to emphasize that, although this scaling would suggest a self-similar behaviour for $F$, the actual asymptotic behavior of $F$ is very different. This discrepancy arises from the fact 
that the mass is not concentrated in the region where $\xi_1 \approx1, \xi_2 \approx 1, \xi_3 \approx 2$. As  will be shown in Section \ref{ssec:conclusion}, the mass 
is distributed across many different scales for $\xi_1 \approx \xi_2$ and for $\xi_2,\xi_3 \in [1,e^{\tau}]$. Nevertheless, the scaling \eqref{def:scalingA<1} 
proves to be very useful when performing computations and when deriving the long-time asymptotics of $f$.


\section{Preliminary tools and Strategy} \label{sec:preliminary} 
We recall that we are interested in investigating the long-time asymptotics of the solutions to \eqref{eq:Bshear}, assuming $K=1$ without loss of generality, and with $B$ satisfying Assumption \ref{ass:B}, namely 
\begin{equation*} 
	\partial_tf {+}  w_2 \partial_{w_1}f=\int_{\mathbb{R}^{3}}dw_{\ast}\int_{S^{2}%
	} \,  |w-w_*|^{-a}
	\left(  M_*'f_{\ast}^{\prime}-M_*f\right) d\omega , 
\end{equation*}
with $0<a< 1$, in the asymptotic regime  $|w| \gg 1$. First, we introduce the functional framework in which we study \eqref{eq:Bshear}  and briefly discuss how to achieve a well-posedness result for the Cauchy problem associated with  \eqref{eq:Bshear} (cf. Section \ref{ssec:WP}). We then outline the strategy employed to obtain information on the long-time asymptotics of the solutions to \eqref{eq:Bshear}. 

\subsection{Functional Framework} 
\label{ssec:WP} 
We start by introducing the precise notion of weak solutions of \eqref{eq:Bshear} that we will use. 
\begin{defn}\label{def:weakSol}
 	Let $T>0$. A function $f \in C\left( [0,+\infty),\M\right)$ is a \emph{weak solution} of \eqref{eq:Bshear} if for every function $\vp \in C^1\left( [0,+\infty),C^{\infty}_c(\R) \right) $ one has
 \begin{align}
     \label{eq:weakf}
		&\int_{\R}\vp(w,T)f(d w,T)- \int_{\R}\vp(w,0)f(dw,0) = \notag \\ &
		\int_0^T dt \int_{\R} f(dw,t)\left[\partial_t\vp {+} w_2\partial_{w_1}\vp+\Ll^*(\vp) \right](w,t)
	\end{align}
where $\Ll^*$ is the adjoint operator of $\Ll$, i.e. 
\begin{equation}
\left(\Ll^* \vp\right)(w)=\int_{\R}d w_* M_* \int_{S^2}B(n\cdot \omega,|w-w_*|)(\vp'-\vp)d \omega.
\end{equation}
\end{defn}

We observe that it is possible to prove a well-posedness result for the Cauchy problem associated with \eqref{eq:Bshear}. More precisely, it is possible to prove that for any initial data $f_0$ in the space of non-negative Radon measures $\M$ there exists a unique weak solution $f \in C\left( [0,+\infty),\M \right)$ of \eqref{eq:Bshear}. 
The following theorem holds.
\begin{thm}\label{thm:wp1}
Let $f_0 \in \M$. There exists a unique weak solution $f \in C\left( [0,+\infty),\M \right)$ in the sense of Definition \ref{def:weakSol} with initial datum $f_0 \in \M$ to the Cauchy problem associated with \eqref{eq:Bshear}.  
\end{thm}
 To prove this result one works with the adjoint equation of \eqref{eq:Bshear}, relying  on the duality provided by the classical Riesz-Markov-Kakutani Theorem. We will not present here the details of the proof, since this is not the aim of this paper. We only mention that such a well-posedness result can be proved by a suitable adaptation of the arguments in \cite{MNV} (cf.~Remark 4) and we refer to \cite{Miele} for the complete proof. Motivated by the weak formulation in Definition \ref{def:weakSol}, we define the \emph{backward in time} adjoint problem which reads as  
\begin{equation}\label{adjEq}
	-\partial_t  {\vp} {- w_2 \partial_{w_1} {\vp} }= \Ll^*({\vp}),
 \end{equation}
 with final condition $\vp(T,v)=\vp_T(v),$ for $ \ T > 0 $ with $\vp_T \in C^{\infty}_c(\R)$. We aim at studying the \emph{forward in time} equation corresponding to \eqref{adjEq}. To this end, we change the time
variable $t\to s=T-t$, with $s \in [0,T]$ for $T>0$, and obtain $\partial_s\vp+w_2\partial_{w_1}\vp-\Ll^*(\vp)=0$. Renaming $s$ as $t$ with a slight abuse of notation we then obtain the \emph{forward in time} adjoint equation 
   \begin{equation}\label{forwardadj}
		{\partial_t\vp- w_2 \partial_{w_1} \vp = \Ll^*(\vp)\, .}
\end{equation}
\paragraph{Approximated model} As discussed in Section \ref{ssec:AsymModel} and Section \ref{sec:heuristics} the model we are actually interested in studying is the approximate model for \eqref{forwardadj} described by \eqref{eq:asymModel}, namely 
\begin{equation}\label{eq:asymDirac}
    \partial_t \vp - w_2 \partial_{w_1}\vp= |w|^{-a} \int_{S^2} (\varphi(w-(w \cdot \omega)\omega)-\varphi(w)) d \omega	.
\end{equation} 
 which arises, following the same steps as above, from the  weak formulation 
 \begin{align}
     \label{eq:weakfAsym} &\int_{\R}\vp(w,T)f(d w,T)- \int_{\R}\vp(w,0)f(dw, 0) = \notag \\ &
		\int_0^T dt \int_{\R} f(dw,t)\left[ {\partial_t\vp(w,t) + w_2\partial_{w_1}\vp(w,t)} +  
        |w|^{-a}\int_{S^2}(\vp(w-(w\cdot\omega)\omega,t)-\vp(w,t))d\omega \right] dt 
\end{align}
for any $\vp \in C^1\left( [0,+\infty),C^{\infty}_c(\R) \right) $. This provides the counterpart of \eqref{eq:weakf} once we approximate $M(w_*)\sim \delta(w_*)$, where $\delta(w_*)$ denotes the Dirac delta measure  $0$.  Notice that the term in $[\dots]=0$ in \eqref{eq:weakfAsym} is equivalent to \eqref{eq:asymDirac} using the time-reversal $t\to -t$.
\bigskip

\subsection{Strategy} 

We are interested in the long-time behaviour of solutions of \eqref{eq:asymModel} and will focus on this from now on. As described in Section in Section \ref{ssec:scalinA>1}, we  look for solutions of the form 
\begin{equation}\label{def:scaling}
    f(w,t)= \frac{1}{t^{\frac{3}{a}-2}} F\left(\xi_1, \xi_2, \xi_3, \tau \right), \;\;\xi_1=\frac{w_1}{t^{\frac{1}{a}}}, \ \ \xi_i=\frac{w_i}{t^{\frac{1}{a}-1}} , \ i=2,3, \quad \tau= \log t \ .  
\end{equation}
Notice that, due to the anisotropy of the
change of variables, also the collision rule and the collision kernel are affected by the scaling \eqref{def:scaling}. 
In order to determine the approximate collision rule for the variable $\xi$ we introduce the matrix
\begin{equation*}
\label{eq:defmatrixS}
    S(t)=\frac{1}{t^{\frac{1}{a}}}\left( \begin{array}{ccc}
        1 & 0 & 0  \\
        0 & t & 0  \\
        0 & 0 & t
    \end{array}\right)
\end{equation*}
with which we can write $\xi=S(t)w$. The approximate collision rule $w'=w-(w \cdot \omega)\omega$ then reads
\begin{equation}\label{eq:CollRule_NewVar}
    w'= T(t)\xi-(T(t)\xi\cdot \omega)\omega
\end{equation}
where we have denoted by $T(t)$ the inverse matrix of $S(t)$ i.e. $T(t)=(S(t))^{-1}$. Therefore applying $S(t)$ to both sides of the above equation yields
\begin{equation}\label{eq:approxCollXi}
    \xi'=\xi-(T(t)\xi \cdot \omega)S(t)\omega.
\end{equation}

Obtaining the strong formulation of the equation satisfied by the rescaled measure $F$ introduced in \eqref{def:scaling}, is rather involved due to the complicated structure of the collision operator. In fact, we look for an equation having the structure
\begin{align}\label{eq:selfSimStrongHyp}
	\partial_{\tau}F & -\frac{1}{a}\partial_{\xi_1}(\xi_1F)-\left(\frac{1}{a}-1\right)\partial_{\xi_2}(\xi_2F)- \left(\frac{1}{a} -1\right)\partial_{\xi_3}(\xi_3F) + \xi_2 \partial_{\xi_1}F  
    =(\mathcal{C} F)(\xi) \, ,
\end{align}
where the linear collision operator $\mathcal{C} (F) $ has to be determined. To this end, it is  convenient to consider the weak formulation of \eqref{eq:selfSimStrongHyp},  relying  on \eqref{eq:approxCollXi}. We have
\begin{defn}\label{weakSolAsymF}
 	Let $\tau^*>0$ and let $F_0 \in \M$. A function $F\in C\left( [0,+\infty),\M\right)$ is a \emph{weak solution} to  \eqref{eq:selfSimStrongHyp} if for every function $\psi \in C^1\left( [0,+\infty),C^{\infty}_c(\R) \right) $ of the form $\psi(\xi,\tau)=\vp(T(t)w,t)$ with $\vp \in C^1([0,+\infty),C_c^{\infty}(\R))$ one has
 \begin{align}
     \label{eq:weakF}
		&\int_{\R}\psi(\xi,\tau^*)F(d\xi,\tau^*)- \int_{\R}\psi(\xi,0)F_0(d\xi) = \notag \\ &
		\int_0^{\tau^*} d\tau \int_{\R} F(d\xi,\tau)\left[\partial_\tau \psi(\xi,\tau)-\frac{1}{a}\xi_1\partial_{\xi_2}\psi(\xi,\tau)-\left(\frac{1}{a}-1\right)\xi_2\partial_{\xi_2}\psi(\xi,\tau)-\left(\frac{1}{a}-1\right)\xi_3\partial_{\xi_3}\psi(\xi,\tau) \right. \notag \\
       &\left. +\xi_2 \partial_{\xi_1}\psi(\xi,\tau)+|T(\tau)\xi|^{-a} \int_{S^2}(\psi(\xi-(T(\tau)\xi\cdot\omega)S(\tau)\omega,\tau)-\psi(\xi,\tau))d\omega \right]. 
	\end{align}
\end{defn}


The derivation of the strong formulation of the equation \eqref{eq:weakF} and, specifically, of the collision operator $\mathcal{C} (F)$ will be performed in Section \ref{ssec:appeq}. In the end, we arrive at the  following strong formulation  
\begin{equation}
 \label{eq:selfSimStrong2_A} 
    \partial_{\tau}F-\frac{1}{a}\partial_{\xi_1}(\xi_1F)-\left(\frac{1}{a}-1\right)\partial_{\xi_2}(\xi_2F)-\left(\frac{1}{a}-1\right)\partial_{\xi_3}(\xi_3F)
   + \xi_2 \partial_{\xi_1} F =
     \left(\mathcal{C}^+F\right)(\xi)-\left(\mathcal{C}^-F\right)(\xi)
\end{equation}
where
\begin{align}\label{eq:selfSimStrong2_operators}
 & \left(\mathcal{C}^+F\right)(\xi):= \frac{8 \delta(\xi_1)}{(2|\tilde \xi|)^{1+a}t^{1-a}}\int_0^{\pi}d\theta \frac{\sin(\theta)}{|\sin(2\theta)|^{1-a}}\int_{\Rp} F\left(\frac{2|\tilde \xi|}{t \sin(2\theta)},\eta, \tau\right)d\eta \\&
  \left(\mathcal{C}^-F\right)(\xi):=-4\pi |\xi_1|^{-a} F \label{eq:selfSimStrong2_operatorsLoss}
\end{align}
with $\xi=(\xi_1,\tilde \xi)\in \mathbb{R}^3$, $\tilde \xi=(\xi_2,\xi_3)\in \mathbb{R}^2$.

\subsubsection{Derivation of the strong formulation \eqref{eq:selfSimStrong2_A} for the approximate equation}  \label{ssec:appeq}

In this section we justify \eqref{eq:selfSimStrong2_A}, which gives the strong formulation of equation \eqref{eq:weakF} and, in particular, of the collision operator $\mathcal{C} (F)$. 
We begin by noticing that, since in the approximation for large velocities in \eqref{eq:asymModel} we expect to have $[w_1]\gg [w_2],[w_3]$, we can assume $[w] \approx [w_1]$  and thus, in the variables $\xi$ \eqref{def:scalingA<1}, we can assume $\xi \sim \xi_1 e_1$. We now observe that the matrix $T(t)$ introduced in \eqref{eq:CollRule_NewVar} behaves asymptotically as 
\begin{equation}\nonumber
    T(t) \sim t^{\frac{1}{a}}\left( \begin{array}{ccc}
       1  & 0 & 0  \\
        0 & 0 & 0 \\
        0 & 0 & 0
    \end{array}\right) \quad \text{as $\, t \rightarrow + \infty$.}
\end{equation}
Then, $w=\left(t^{\frac{1}{a}} \xi_1, t^{\frac{1}{a}-1}\xi_2, t^{\frac{1}{a}-1}\xi_3\right)$ and, as $t \rightarrow +\infty$, $w= t^{\frac{1}{a}}\left( \xi_1,0,0\right)$. Therefore, it follows that
\begin{equation}\label{eq:wprime}
w'=(I-\omega \otimes \omega)w \sim t^{\frac{1}{a}}\xi_1(I-\omega \otimes \omega)e_1=t^{\frac{1}{a}} \xi_1(P_\omega^1(e_1),P_\omega^2(e_1),P_\omega^3 (e_1))
\end{equation}
where $P_i$, $i=1,2,3$ denote the components of the orthogonal projection, i.e. $P_{\omega}(w)=(I-\omega \otimes \omega)w=(P_1,P_2,P_3).$
In particular, if $w=e_1$ we have 
\begin{equation*}\label{eq:Pcomp}
  P_{\omega}(e_1)=(P_\omega^1(e_1),P_\omega^2(e_1),P_\omega^3 (e_1))=   
(1-\omega_1^2, -\omega_2\omega_1,-\omega_3\omega_1).
\end{equation*}
Therefore, using \eqref{def:scaling} in \eqref{eq:wprime}, we obtain
\begin{equation}\label{eq:xiprime}
\xi'=  \left(\xi_1P_\omega^1(e_1),t\xi_1P_\omega^2(e_1),t\xi_1P_\omega^3(e_1)\right) 
\end{equation}
 and thus
\begin{equation}\label{eq:Psiw'}
\psi(\xi',\tau)\sim \psi(\xi_1P_\omega^1(e_1),t\xi_1P_\omega^2(e_1),t\xi_1P_\omega^3(e_1),\tau).
\end{equation} 
We now observe that if the test function $\psi(\xi,t)$ is compactly supported or, alternatively, gives a relevant contribution in the region where $\vert \xi\vert \approx1$, using \eqref{eq:xiprime}, the test function $\psi(\xi',t)$ gives a relevant contribution for  $\vert t\xi_1\vert \leq C$ or, equivalently, $\vert \xi_1\vert \leq \frac C t$ as $t \rightarrow + \infty$ and for most of the values of $\omega$. 
Hence we can take $\xi_1  P_\omega^1(e_1)\sim 0$ in \eqref{eq:Psiw'} for large $t$ and this gives $\psi(\xi',\tau)\sim \psi(0,t\xi_1P_\omega^2(e_1),t\xi_1P_\omega^3(e_1))$. This approximation can be justified \emph{a posteriori} by the fact that, for the solution that we will compute later, the collisions take place for $|\xi_1| \approx 1$. This allows us to use the approximation $\omega=\omega_1 e_1$ even if $|\xi_1| \approx 1$.

Then, taking a compactly supported test function in the weak formulation \eqref{eq:weakF}, we obtain
\begin{align}\label{eq:weak2}
    \qquad \qquad \qquad &\int_{\R}F(\xi,\tau_*)\psi(\xi,\tau_*)d\xi-\int_{\R}F_0(\xi)\psi(\xi,0)d\xi\nonumber \\ &= \int_0^{\tau_*}d\tau\int_{\R}F(\xi,\tau)\left(\partial_{\tau}\psi-\frac{1}{a}\xi_1\partial_{\xi_1}\psi-\left(\frac{1}{a}-1\right)\xi_2\partial_{\xi_2}\psi -\left(\frac{1}{a}-1\right)\xi_3\partial_{\xi_3}\psi+\xi_2 \partial_{\xi_1}\psi\right)(\xi,\tau)
     \nonumber\\&
    \quad +\int_0^{\tau_*}d\tau\int_{\R}F(\xi,\tau)
    |\xi_1|^{-a}\int_{S^2}(\psi( 0
    ,\xi_1tP_\omega^2(e_1),\xi_1tP_\omega^3(e_1),\tau)-\psi(\xi,\tau))d\omega d\xi
\end{align}
or, equivalently,
\begin{align}\label{eq:weakTimeDep}
    \qquad \qquad \qquad &\frac{d}{d\tau}\int_{\R}F(\xi,\tau)\psi(\xi,\tau)d\xi \notag\\ &= \int_{\R}F(\xi,\tau)\left(\partial_{\tau}\psi-\frac{1}{a}\xi_1\partial_{\xi_1}\psi-\left(\frac{1}{a}-1\right)\xi_2\partial_{\xi_2}\psi -\left(\frac{1}{a}-1\right)\xi_3\partial_{\xi_3}\psi+\xi_2 \partial_{\xi_1}\psi\right)(\xi,\tau)+
     \nonumber\\
     &+\int_{\R}F(\xi,\tau)
    |\xi_1|^{-a}\int_{S^2}(\psi( 0
    ,\xi_1tP_\omega^2(e_1),\xi_1tP_\omega^3(e_1),\tau)-\psi(\xi,\tau))d\omega d\xi.
\end{align}
From now on we will assume that $F$ is smooth in the region where $0<|\xi_1|<|\xi_2|$. Our goal is to write the strong formulation of the problem satisfied by $F$, starting from \eqref{eq:weak2}. Set $\tilde \xi=(\xi_2,\xi_3)$. Due to the form of the collision operator in \eqref{eq:weak2}, we have that the collision dynamics is such that the outgoing velocities are always sent to a neighbourhood of the plane $\{\xi_1=0\}$. We will assume in the following that the jumps send the particle directly to the plane $\{\xi_1=0\}$. Consider first in \eqref{eq:weak2} the term containing the gain operator, namely 
\begin{equation*}
    |\xi_1|^{-a} \int_{S^2}\psi(0,\xi_1tP_\omega^2(e_1),\xi_1tP_\omega^3(e_1))d \omega .
\end{equation*}
Then, we can write 
\begin{align*}
    &\int_{\R}d\xi \  F(\xi,\tau) \Ll^{\ast}_{+}(\psi)(\xi)
     \notag \\
    &=\int_{\R}d\xi \  |\xi_1|^{-a}F(\xi,\tau)\int_{S^2}d\omega\int_{\R}\delta(y_1)\delta(y_2-\xi_1tP_\omega^2(e_1))\delta(y_3-\xi_1tP_\omega^3(e_1))\psi(y_1,y_2,y_3,\tau)dy\notag \\
&=\int_{\R} dy \ \psi(y_1,y_2,y_3,\tau)\int_{\R}d\xi\  |\xi_1|^{-a}F(\xi,\tau)\int_{S^2}\delta(y_1)\delta(y_2-\xi_1tP_\omega^2(e_1))\delta(y_3-\xi_1tP_\omega^3(e_1))d\omega.
\end{align*}
We consider the equation
\begin{equation}\label{eq:operator1}
\delta(y_1)\int_{\R}d\xi \ |\xi_1|^{-a}F(\xi,\tau) \int_{S^2}\delta(y_2-\xi_1tP_\omega^2(e_1))\delta(y_3-\xi_1tP_\omega^3(e_1))d\omega.
\end{equation}
We set $\tilde{y}=(y_2,y_3) \in \mathbb{R}^2$ and observe that it is possible to write the vector $(0,\xi_1tP_\omega^1(e_1),\xi_1tP_\omega^2(e_1))$ as $Q(\mathrm{P}_{\omega}^{\perp}(\xi_1te_1))$ where $Q$ is the orthogonal projection on $e_1$, i.e. $Q(v)=v-(v \cdot e_1)e_1$ with $v \in \R$. This implies
\begin{equation*}
    Q(\mathrm{P}_{\omega}^{\perp}(e_1))=Q(e_1-(\omega \otimes \omega)e_1)=-Q(\omega)(\omega \cdot e_1).
\end{equation*}
Using spherical coordinates on the unit sphere, with the North pole at $e_1$, the vector $\omega$ then reads
\begin{equation*}
    \omega=\left( \begin{array}{c}
         \cos \theta \\
         \sin \theta \cos \vp \\
         \sin \theta \sin \vp
    \end{array} \right), \quad \theta \in [0,\pi], \vp \in [0,2\pi].
\end{equation*}
 Thus, 
\begin{equation}\label{eq:Qomega}
 -Q(\omega)(\omega \cdot e_1)=-\cos(\theta)\left( \begin{array}{c}
      0 \\
      \sin \theta \cos \vp \\
      \sin \theta \sin \vp
 \end{array}\right).
\end{equation}
Combining \eqref{eq:Qomega} with \eqref{eq:operator1}, we then have
 \begin{align}
 &\
\delta(y_1)\int_{S^2}d\omega \ \delta(y_2-\xi_1tP_\omega^2(e_1))\delta(y_3-\xi_1tP_\omega^3(e_1))\nonumber\\
&= \delta(y_1)\int_0^{\pi}d\theta \int_0^{2\pi} d\vp \sin(\theta)\delta\left(\tilde y-\xi_1t\cos(\theta)\left(
    \sin(\theta)\cos(\vp), \sin(\theta)\sin(\vp)\right)\right)  \, .\label{eq:operator2}
 \end{align} 
Setting  $\tilde y= |\tilde y|(\cos(\psi_0),\sin(\psi_0))$ we can rewrite 
\eqref{eq:operator2} as 
\begin{align} 
 \delta(y_1)& \int_0^{\pi}d\theta \int_0^{2\pi} d\vp\, \sin(\theta) \delta(|\tilde y|\cos(\psi_0)-\xi_1t\cos(\theta)\sin(\theta)\cos(\vp)) \; \delta(|\tilde y|\sin(\psi_0)-\xi_1t\cos(\theta)\sin(\theta)\sin(\vp))\nonumber\\&
    = \ \int_0^{\pi}d\theta \int_0^{2\pi} d\vp\, \sin(\theta) \delta(X_1(\theta,\vp))\delta(X_2(\theta,\vp))\notag
\end{align}
with
\begin{align}
    X_1(\theta,\vp) &:=|\tilde y|\cos(\psi_0)-\xi_1 t\cos(\theta)\sin(\theta)\cos(\vp); \label{eq:F1} \\
    X_2(\theta,\vp) &:= |\tilde y|\sin(\psi_0)-\xi_1t\cos(\theta)\sin(\theta)\sin(\vp) \label{eq:F2}.
\end{align}
We now set $X_1(\theta,\vp)=\zeta_1,\, X_2(\theta,\vp)=\zeta_2$ so that 
\begin{equation}\label{eq:operator4} 
    \int_0^{\pi}d\theta \int_0^{2\pi} d\vp\, \sin(\theta) \delta(X_1(\theta,\vp))\delta(X_2(\theta,\vp))= \int_{\Rl}  d\zeta_1 \int_{\Rl} d\zeta_2 \sin(\theta)\,\delta(\zeta_1)\delta( \zeta_2) \left(\Big| \frac{\partial(X_1,X_2}{\partial(\theta,\vp)} \Big|\right)^{-1}.
\end{equation}
We now look for solutions to $X_1(\theta,\vp)=0$, and $X_2(\theta,\vp)=0$. Dividing \eqref{eq:F2} by \eqref{eq:F1}  we obtain the following trigonometric relation
\begin{equation*}
\tan(\vp)=\tan(\psi_0)
\end{equation*}
which is satisfied for $\vp=\psi_0$ 
or for $\vp=\psi_0+n\pi$ with $n\in \mathbb{Z}$. Therefore, from \eqref{eq:F1}, we have
\begin{equation*}
    \begin{cases}
    \xi_1=\frac{2|\tilde y|}{t\sin(2 \theta)} & \text{for $\vp=\psi_0$} \quad \text{with} \;\ 0<\theta<\frac \pi 2;\\
        \xi_1=-\frac{2|\tilde y|}{t\sin(2 \theta)} & \text{for $\vp=\psi_0+\pi$} \quad \text{with} \;\ \frac \pi 2<\theta<\pi.
    \end{cases}
\end{equation*}
Notice that since $(X_1,X_2)$ depends linearly on $\xi_1$ we have $\frac{\partial(X_1,X_2)}{\partial(\theta,\vp)}=\frac{\partial(X_1,X_2)}{\partial(\theta,\xi_1)}$. Hence,
\begin{equation*}
\frac{\partial(X_1,X_2)}{\partial(\theta,\xi_1)}=\Bigg|\det \left(\begin{array}{cc}
       -\frac{1}{2}\xi_1 t\sin(2\theta)\sin(\vp)  & \frac{1}{2}\xi_1 t\sin(2\theta)\cos(\vp)  \\
        \frac{1}{2}t\sin(2\theta)\cos(\vp) & \frac{1}{2}t\sin(2\theta)\sin(\vp)
\end{array}\right) \Bigg|=\frac{1}{4}\vert \xi_1\vert t^2 \sin^2(2\theta).
\end{equation*}
In conclusion, we then obtain
\begin{equation}\label{eq:operator5}
\left(\delta(\vp-\psi_0)\delta\left(\xi_1-\frac{2|\tilde y|}{t\sin(2\theta)}\right)+\delta(\vp-\psi_0+\pi)\delta\left(\xi_1+\frac{2|\tilde y|}{t\sin(2\theta)}\right)\right)\frac{4}{\vert \xi_1\vert t^2\sin^2(2\theta)}.
\end{equation}
Recalling that $\tilde{\xi}=(\xi_2,\xi_3)$ and $\tilde y= |\tilde y|(\cos(\psi_0),\sin(\psi_0))$, plugging \eqref{eq:operator5}   into \eqref{eq:operator1} and using \eqref{eq:operator4} we then obtain
\begin{align}\label{eq:operator6}
&\delta(y_1)\left\{\int_{\mathbb{R}}d\xi_1 |\xi_1|^{-a}\int_{\mathbb{R}^2}d\tilde{\xi} \, F(\xi_1,\tilde{\xi},\tau) \int_0^{\pi}d\theta \int_0^{2\pi}d\vp\frac{4 \,\sin(\theta) }{\vert\xi_1\vert t^2\sin^2(2\theta)} \left[\delta(\vp-\psi_0)\delta\left(\xi_1-\frac{2|\tilde y|}{t\sin(2\theta)}\right)  \right]\right. \nonumber \\&
    \quad + \left.\int_{\mathbb{R}}d\xi_1 |\xi_1|^{-a}\int_{\mathbb{R}^2}d\tilde{\xi}\, F(\xi_1,\tilde{\xi},\tau) \int_0^{\pi}d\theta\int_0^{2\pi}d\vp\frac{4\, \sin(\theta)}{\vert\xi_1\vert t^2\sin^2(2\theta)}\left[\delta(\vp-\psi_0+\pi)\delta\left(\xi_1+\frac{2|\tilde y|}{t\sin(2\theta)}\right)  \right] \right\}\nonumber \\&=\delta(y_1)\left\{\int_0^{\pi} d\theta \left(\frac{2|\tilde y|}{t \sin(2\theta)}\right)^{-a}\int_{\Rp}d\tilde{\xi} \, F\left(\frac{2|\tilde y|}{t \sin(2\theta)},\tilde \xi, \tau\right) \frac{4\sin(\theta)}{t^2\sin^2(2\theta)\frac{2|\tilde y|}{t \sin(2\theta)}}\right. \nonumber \\&
    \quad + \left. \int_0^{\pi}d\theta \left(\frac{2|\tilde y|}{t \sin(2\theta)}\right)^{-a}\int_{\Rp}d\tilde{\xi} \, F\left(-\frac{2|\tilde y|}{t \sin(2\theta)},\tilde \xi, \tau\right) \frac{4\sin(\theta)}{t^2\sin^2(2\theta)\frac{2|\tilde y|}{t \sin(2\theta)}}\right\} =:  (I_1+I_2)\delta(y_1)
\end{align} 
Setting $\theta=\pi-\tilde \theta$ in the second integral, since $\sin(\theta)=\sin(\pi-\tilde\theta)=\sin(\tilde\theta)$ and $\sin(2\theta)=-\sin(2\tilde\theta)$, then we get $I_1=I_2$ and thus from \eqref{eq:operator6} we obtain
\begin{align*}
 (I_1+I_2)\delta(y_1)&=2I_1\delta(y_1) \nonumber \\&=2\delta(y_1)\int_0^\pi d\theta \left(\frac{2|\tilde y|}{t \sin(2\theta)}\right)^{-a} \int_{\Rp} d\tilde {\xi} \, F\left(\frac{2|\tilde y|}{t \sin(2\theta)},\tilde \xi,\tau\right) \frac{2 \sin(\theta)}{|\tilde y|t \sin(2\theta)} \nonumber \\&
    =\frac{4 \delta(y_1)}{2^{a}(|\tilde y|)^{1+a}t^{1-a}}\int_0^{\pi}d\theta \frac{\sin(\theta)}{|\sin(2\theta)|^{1-a}}\int_{\Rp} F\left(\frac{2|\tilde y|}{t \sin(2\theta)},\tilde{\xi}, \tau\right)d\tilde{\xi} \ .
\end{align*} 
Renaming the variable $y$ to $\xi$, we finally obtain \eqref{eq:selfSimStrong2_A}, namely:
\begin{equation*}
    \partial_{\tau}F-\frac{1}{a}\partial_{\xi_1}(\xi_1F)-\left(\frac{1}{a}-1\right)\partial_{\xi_2}(\xi_2F)-\left(\frac{1}{a}-1\right)\partial_{\xi_3}(\xi_3F)
   + \xi_2 \partial_{\xi_1}F =
     \left(\mathcal{C}^+F\right)(\xi)-\left(\mathcal{C}^-F\right)(\xi)
\end{equation*}
where $\left(\mathcal{C}^+F\right)$ and $\left(\mathcal{C}^- F\right)$ are as in \eqref{eq:selfSimStrong2_operators}-\eqref{eq:selfSimStrong2_operatorsLoss}
with $\xi=(\xi_1,\tilde \xi)$, $\tilde \xi=(\xi_2,\xi_3)$.  We observe that for \eqref{eq:selfSimStrong2_A} the mass conservation property holds true, as it can be checked by means of  direct computations, see Section \ref{sec:auxiliary} for the  detailed proof.

\bigskip

\section{Long-time asymptotics (Proof of Main Result \ref{thm:ConvergenceF})} 
\label{sec:longtime}

The goal of this section is to investigate the long-time asymptotics of the profile $F(\xi,\tau)$ defined as in  \eqref{def:scaling}. The evolution equation for $F(\xi,\tau)$, with $\xi\in\mathbb{R}^3$, $\tau\in \mathbb{R}_{+}$  is given by \eqref{eq:selfSimStrong2_A} which has been introduced in Section \ref{sec:preliminary}. It reads 
\begin{align} \label{eq:selfSimStrong2}
    \partial_{\tau}F &- \frac{1}{a}\partial_{\xi_1}(\xi_1F)-\left(\frac{1}{a}-1\right)\partial_{\xi_2}(\xi_2F)-\left(\frac{1}{a}-1\right)\partial_{\xi_3}(\xi_3F)
   + \xi_2 \partial_{\xi_1}F  \nonumber \\
   & = \frac{8 \delta(\xi_1)}{(2|\tilde \xi|)^{1+a}t^{1-a}}\int_0^{\pi}d\theta \frac{\sin(\theta)}{|\sin(2\theta)|^{1-a}}\int_{\Rp} F\left(\frac{2|\tilde \xi|}{t \sin(2\theta)},\eta, \tau\right)d\eta-4\pi |\xi_1|^{-a} F
\end{align}
where $\xi=(\xi_1,\tilde \xi)$ with $\tilde{\xi}=(\xi_2,\xi_3)$. We recall that we are interested in the case $0<a< 1$ where $a=  \vert \gamma \vert $.

\subsection{Long-time behaviour of the solution of a reduced two-dimensional problem}\label{sec:strategy}

To study the long-time behaviour of the solution $F=F(\xi_1,\xi_2,\xi_3,\tau)$ to \eqref{eq:selfSimStrong2}, we first consider a reduced two-dimensional problem obtained by integrating in the variable $\xi_3$. We then recover the long-time behaviour of 
$F(\xi_1,\xi_2,\xi_3,\tau)$ (cf. Section \ref{ssec:longTimeF}) from the asymptotics of the reduced solution $G=G(\xi_1,\xi_2,\tau)$ obtained in this section. More precisely, we set 
\begin{equation}\notag
    G(\xi_1,\xi_2, \tau) =\int_{\Rl} F(\xi_1,\xi_2,\xi_3,\tau)\, d \xi_3\, .
\end{equation}
Assuming that $F$ vanishes sufficiently fast at infinity, the reduced function $G$ satisfies
\begin{align} 
    \partial_{\tau}G &-\frac{1}{a}\partial_{\xi_1}(\xi_1G)-\left(\frac{1}{a}-1\right)\partial_{\xi_2}(\xi_2G) + \xi_2 \partial_{\xi_1}G \nonumber \\
    & = \frac{16 \, \delta(\xi_1)}{t^{1-a}}\int_0^{\pi}d\theta \frac{\sin(\theta)}{|\sin(2\theta)|^{1-a}}\int_{0}^{+\infty} \frac{d\xi_3}{(2|\tilde \xi|)^{1+a}}\int_{\Rl} G\left(\frac{2|\tilde \xi|}{t \sin(2\theta)},\eta, \tau\right)d\eta-4\pi |\xi_1|^{-a} G \label{eq:strongG}. 
\end{align}
This equation is closed in $G$, so the long-time behaviour of $G$ can be derived independently of $F$'s asymptotics. A crucial point is that 
equation \eqref{eq:strongG} can be reformulated as a suitable boundary-value problem — specifically, a delay equation.  More precisely, 
\begin{align} 
    \partial_{\tau}G&-\frac{1}{a}\partial_{\xi_1}(\xi_1G)-\left(\frac{1}{a}-1\right)\partial_{\xi_2}(\xi_2G) + \xi_2 \partial_{\xi_1}G= -4\pi |\xi_1|^{-a} G\, , \label{eq:EvolG}
    \\&
    G(0,\xi_2,\tau)  = \frac{16 }{\xi_2t^{1-a}}\int_0^{\pi}d\theta \frac{\sin(\theta)}{|\sin(2\theta)|^{1-a}}\Ii \frac{d\xi_3}{(2|\tilde \xi|)^{1+a}}\int_{\Rl} G\left(\frac{2|\tilde \xi|}{t \sin(2\theta)},\eta, \tau\right)d\eta . \label{eq:Gboundary}
\end{align}

Our goal is to look for a proper reformulation of the expression on the right-hand side of \eqref{eq:Gboundary}. By means of direct computations it is possible to prove the following. 
\begin{lem}\label{lem:bdvalue}
Let $G(0,\xi_2,\tau)$ be defined as in \eqref{eq:Gboundary}. We define the function $H: [0,1] \rightarrow \Rl$ as 
\begin{equation}\label{eq:defHs}
    H(s) = 2^{-\frac{3}{2}} \int_{-\pi/2}^{\pi/2} \frac{d \theta}{\sqrt{1-\sqrt{1-s^2}\sin \theta}} .
\end{equation}
Then, we can rewrite  the boundary value $G(0,\xi_2,\tau)$ as 
\begin{equation}
G(0,\xi_2,\tau)= \frac{16}{\xi_2 t} \int_{0}^{+\infty} d \eta \int_{\frac{2|\xi_2|}{t}}^{+\infty} dy\, G(y,\eta,\tau) \frac{1}{y^{1+a}} \, H\left(\frac{2\xi_2}{ty}\right) \, .
\end{equation}
\end{lem}

Moreover, it is possible to show that the function $H(s)$ introduced above diverges logarithmically close to the origin while it stays bounded as $s\to 1$. In fact, we have the following.
\begin{lem}\label{lem:asymHs}
The function $H(s)$ defined as in \eqref{eq:defHs} has the following asymptotic behaviour
\begin{equation} \notag
    H(s)\sim -\log(s) ,  \quad \text{as }\ \ s \rightarrow 0^+\
\end{equation}
and we have that $H(1)= 2^{-\frac{3}{2}}\pi.$
\end{lem}
The proofs of Lemmas \ref{lem:bdvalue} and \ref{lem:asymHs} are postponed to Section \ref{sec:auxiliary}.
\medskip

\noindent\textbf{Characteristic curves} We now consider  \eqref{eq:Gboundary} and study the characteristics associated with it. They read
\begin{align}\label{eq:characteristics}
    \frac{d \xi_1}{ds} & = -\frac{1}{a} \xi_1+ \xi_2; \notag \\
    \frac{d \xi_2}{ds} & = -\left(\frac{1}{a}-1\right) \xi_2; \notag \\
    \frac{d G}{ds} & =  \left(\frac{2}{a}-1\right)G-4 \pi |\xi_1(s)|^{-a}G
\end{align}
with initial data, at time  $\tau_0$,  $\xi_{1}(\tau_0)=0, \,\xi_{2}(\tau_0)=\xi_{2,0}$. 
Solving   the system \eqref{eq:characteristics} for $0 \leq \tau_0 < \tau$ we then obtain
\begin{align}
    \xi_1(\tau) &=\xi_{2,0}  e^{-\left(\frac{1}{a}-1\right)(\tau-\tau_0)}- \xi_{2,0} e^{-\frac{1}{a}(\tau-\tau_0)}= \xi_{2,0} e^{-\frac{1}{a}(\tau-\tau_0)} \left(e^{(\tau-\tau_0)}-1 \right)\label{eq:xi1}; \\
    \xi_2(\tau) & = \xi_{2,0} e^{-\left(\frac{1}{a}-1\right)(\tau-\tau_0)}\label{eq:xi2};  \\
    G(\xi_1,\xi_2,\tau) & = G(0,\xi_{2,0},\tau_0) \exp\left(\left(\frac{2}{a}-1\right)(\tau-\tau_0)\right)\exp\left(- 4 \pi \int_{\tau_0}^{\tau}|\xi_1(s)|^{-a}ds\right)\label{eq:charG}.
\end{align}
Notice that from \eqref{eq:xi1} and \eqref{eq:xi2} we have 
\begin{equation}\label{eq:characteristics2}
e^{-(\tau-\tau_0)} =1- \frac{\xi_1(\tau)}{\xi_2(\tau)}.
\end{equation}
Moreover, the integral inside the exponential in \eqref{eq:charG} can be computed explicitly and gives
\begin{align}
    \int_{\tau_0}^{\tau}\xi_1(s)^{-a}ds & = (\xi_{2,0})^{-a} \int_{\tau_0}^{\tau} \left( e^{-\frac{1}{a}(s-\tau_0)} \left(e^{(s-\tau_0)}-1 \right)\right)^{-a}ds = (\xi_{2,0})^{-a} \int_0^{\tau-\tau_0} \left( e^{-\frac{1}{a}s} \left(e^{s}-1 \right)\right)^{-a}ds  \notag \\
    & = (\xi_{2,0})^{-a} \int_0^{\tau-\tau_0} e^s(e^s-1)^{-a}ds = \frac{(\xi_{2,0})^{-a} }{1-a} \left(e^{\tau-\tau_0}-1\right)^{1-a} \notag \\
    &= \frac{\xi_2^{-a}}{1-a}e^{-(1-a)(\tau-\tau_0)}(e^{\tau-\tau_0}-1)^{1-a}\notag
\end{align}
where in the last equality we have used \eqref{eq:xi2}. From \eqref{eq:characteristics2} we also have that
\begin{equation*}
    e^{\tau-\tau_0}-1=\frac{1}{1-\xi_1/\xi_2}-1 = \frac{\xi_2}{\xi_2-\xi_1}-1=\frac{\xi_1}{\xi_2-\xi_1}.
\end{equation*}
Therefore using \eqref{eq:xi2}  we obtain
\begin{align}
    \int_{\tau_0}^{\tau}\xi_1(s)^{-a}ds= \frac{\xi_2^{-a}}{1-a}\left(\frac{\xi_2-\xi_1}{\xi_2}\right)^{1-a} \left(\frac{\xi_1}{\xi_2-\xi_1}\right)^{1-a}= \frac{1}{1-a}\frac{\xi_1^{1-a}}{\xi_2}.\label{eq:lossCharacteristics}
\end{align}
Combining  \eqref{eq:charG} and  \eqref{eq:lossCharacteristics} we arrive at
\begin{equation}\label{eq:evGchar}
    G(\xi_1,\xi_2,\tau)=G\left(0,\frac{(\xi_2)^{\frac{1}{a}}}{(\xi_2-\xi_1)^{\frac{1}{a}-1}},\tau_0\right) \left(1-\frac{\xi_1}{\xi_2}\right)^{-\left(\frac{2}{a}-1\right)} \exp\left(-\frac{4 \pi}{1-a}\frac{\xi_1^{1-a}}{\xi_2}\right)
\end{equation}
where we have assumed that $0\leq \xi_1 \leq \xi_2$ or that $-\xi_2\leq-\xi_1\leq0$. Without loss of generality, we will then take $0\leq \xi_1\leq \xi_2$. \medskip

\noindent\textbf{Delay equation for $G(0,\xi_2,\cdot)$.}  
Combining Lemma \ref{lem:bdvalue} and \eqref{eq:evGchar}, equation \eqref{eq:Gboundary} becomes
\begin{align}
    &G(0,\xi_2,\tau) =   \notag \\
    &\frac{16}{\xi_2t} \int_{\frac{2|\xi_2|}{t}}^{+\infty} d\eta \int_{\frac{2|\xi_2|}{t}}^{\eta} \frac{dy}{y^{a+1}} G \left(0, \frac{\eta^{\frac{1}{a}}}{(\eta-y)^{\frac{1}{a}-1}},\tau_0 \right) \left(1-\frac{y}{\eta}\right)^{-\left(\frac{2}{a}-1\right)}\exp \left(-\frac{4 \pi}{1-a} \frac{y^{1-a}}{\eta}\right)   H \left(\frac{2\xi_2}{ty}\right). \label{eq:Gboundary1}
\end{align}

Notice that in \eqref{eq:Gboundary1} we restricted the variables to the domain $0 \leq y \leq \eta,$ and $\frac{2|\xi_2|}{t} \leq \eta$ as a consequence of the fact that we are only considering $0 \leq \xi_1 \leq \xi_2$. 

It is worth noticing that \eqref{eq:Gboundary1} is a delay equation for the function $G$ since the value at the boundary $\xi_1=0$ and at time $\tau$ is computed from $G$ itself evaluated at previous times $0 \leq\tau_0 \leq \tau$.

We now perform the change of variables $y \to \omega= \frac{\eta^{\frac{1}{a}}}{(\eta-y)^{\frac{1}{a}-1}}$. Hence, we can write $y=\eta \left(1-\left(\frac{\eta}{\omega}\right)^{\frac{a}{1-a}}\right)$ and then 
$1-\frac{y}{\eta}=\left(\frac{\eta}{\omega}\right)^{\frac{a}{1-a}}$, $\left(1-\frac{y}{\eta}\right)^{-\left(\frac{2}{a}-1\right)}=\left(\frac{\eta}{\omega}\right)^{-\frac{2-a}{1-a}}$,  with differential 
$dy=\frac{a}{1-a} \left(\frac{\eta}{\omega}\right)^{\frac{1}{1-a}}d\omega$. Thus \eqref{eq:Gboundary1} becomes
\begin{align}
    G(0,\xi_2,\tau) =& \frac{16}{\xi_2t} \frac{a}{1-a} \int_{\frac{2|\xi_2|}{t}}^{+\infty} d\eta \int_{\omega_0}^{+\infty} G(0,\omega,\tau_0)\frac{d\omega}{\eta^{a+1}\left(1-\left(\frac{\eta}{\omega}\right)^{\frac{a}{1-a}}\right)^{a+1}} \left(\frac{\eta}{\omega}\right)^{\frac{1}{1-a}}\left(\frac{\eta}{\omega}\right)^{-\frac{2-a}{1-a}} \notag \\
    & \times\exp\left(-4 \pi\frac{\eta^{1-a}\left(1-\left(\frac{\eta}{\omega}\right)^{\frac{a}{1-a}}\right)^{1-a}}{(1-a)\eta}\right) H \left(\frac{2\xi_2}{t \eta \left(1-\left(\frac{\eta}{\omega}\right)^{\frac{a}{1-a}}\right)}\right) \notag \\
    &= \frac{16}{\xi_2t} \frac{a}{1-a} \int_{\frac{2|\xi_2|}{t}}^{+\infty} d\eta \int_{\omega_0}^{+\infty} G(0,\omega,\tau_0)\frac{d\omega}{\eta^{a+1}\left(1-\left(\frac{\eta}{\omega}\right)^{\frac{a}{1-a}}\right)^{a+1}}\frac{\omega}{\eta}\notag \\
    & \times \exp\left(-\frac{4 \pi\,\eta^{-a}}{1-a}\left(1-\left(\frac{\eta}{\omega}\right)^{\frac{a}{1-a}}\right)^{1-a}\right) H \left(\frac{2\xi_2}{t \eta \left(1-\left(\frac{\eta}{\omega}\right)^{\frac{a}{1-a}}\right)}\right)    \label{eq:GboundaryOmega}
\end{align}
with $\omega_0=\frac{\eta^{\frac{1}{a}}}{\left(\eta-\frac{2\xi_2}{t}\right)^{\frac{1}{a}-1}}$.
\bigskip
So far in this subsection we have not introduced any new approximation besides those already used to derive \eqref{eq:selfSimStrong2} while here we introduce the two fundamental approximations 
which we will rely on in the following.  

 \medskip
It is possible to approximate $\tau_0$ as a function of $\tau$ when $\tau$ is large enough in \eqref{eq:GboundaryOmega} because the contribution of some logarithmic terms is negligible. 
This approximation is suggested by the fact that $\frac{t}{t_0}=e^{\tau-\tau_0}$ and by the fact that the main contribution to the collision is given for $\tau-\tau_0 \gg 1$ as suggested by \eqref{eq:increasingJumps}.
This is due to the effect of the shear flow on the system which elongates the trajectories of the particle increasingly with time, or equivalently increasingly with every jump. It is important to emphasize that the system that we are considering in this paper exhibits a behaviour far from that of a diffusion process in which the main effects are due to small deflections of the velocities in short times.  Using this assumption and the fact that $e^{-(\tau-\tau_0)}= \left(\frac{\eta}{\omega}\right)^{\frac{a}{1-a}}$ we can expect the main contribution to the integrals to be given for $\eta/\omega \ll 1$. Going back to \eqref{eq:GboundaryOmega}  and using that $t=e^{\tau}$ we obtain
\begin{align}\label{eq:GBoundary1}
   G(0,\xi_2,\tau) & = \frac{16}{\xi_2}e^{-\tau} \frac{a}{1-a} \int_{\frac{2|\xi_2|}{e^{\tau}}}^{+\infty} \frac{\exp\left(-\frac{4 \pi \eta^{-a}}{1-a}\right)}{\eta^{a+2}}   H \left(\frac{2\xi_2}{e^{\tau} \eta}\right) d\eta \int_{\omega_0}^{+\infty} \omega G(0,\omega,\tau_0) \;d\omega 
\end{align}
with $\omega_0=\frac{\eta^{\frac{1}{a}}}{\left(\eta-\frac{2\xi_2}{e^{\tau}}\right)^{\frac{1}{a}-1}}$. 
We now define the function
\begin{equation}\label{def:U}
    e^{-2\tau}U(\zeta,\tau):=G(0,\xi_2,\tau)\quad \text{ with }\ \zeta=\xi_2e^{-\tau}
\end{equation}
and we notice that $G(0,\omega,\tau_0)=e^{-2\tau_0}U(\omega e^{-\tau_0},\tau_0)$. Recalling that $\tau_0=\tau-\frac{a}{1-a} \log \left(\frac{\omega}{\eta}\right)$, we obtain that
\begin{align}\label{eq:Uomega}
    U(\zeta,\tau) &= \frac{16a}{(1-a)\zeta} e^{-2\tau}\int_{2\zeta}^{+\infty} \frac{\exp\left(-\frac{4 \pi \eta^{-a}}{1-a}\right)}{\eta^{2+a}}H\left(\frac{2\zeta}{\eta}\right) d \eta \notag \\
    & \times \int_{\omega_0}^{+\infty} \omega \left(\frac{\omega}{\eta}\right)^{\frac{2a}{1-a}} U \left( \left(\frac{\omega}{\eta^a}\right)^{\frac{1}{1-a}}e^{-\tau}, \tau-\frac{a}{1-a} \log \left(\frac{\omega}{\eta}\right)\right) d \omega.
\end{align}
Setting $\omega=e^{(1-a)\tau}\sigma$ gives
\begin{align}\label{eq:Usigma}
     U(\zeta,\tau) &= \frac{16a}{(1-a)\zeta} \int_{2\zeta}^{+\infty} \frac{\exp\left(-\frac{4 \pi \eta^{-a}}{1-a}\right)}{\eta^{2+a}}H\left(\frac{2\zeta}{\eta}\right) d \eta \notag \\
    & \times \int_{\omega_0e^{-(1-a)\tau}}^{+\infty} \sigma \left(\frac{\sigma}{\eta}\right)^{\frac{2a}{1-a}} U \left( \left(\frac{\sigma}{\eta^a}\right)^{\frac{1}{1-a}}, (1-a)\tau-\frac{a}{1-a} \log \left(\frac{\sigma}{\eta}\right)\right) d \omega.
\end{align}

We can now state the second approximation that we use in the following. \\ 

\medskip
We assume that 
    \begin{equation}
        \int_{\omega_0e^{-(1-a)\tau}}^{+\infty} d\omega \, [\dots]\sim \int_0^{+\infty}d\omega \, [\dots], \quad \text{where} \quad \omega_0=\frac{\eta^{\frac{1}{a}}}{(\eta-2\zeta)^{\frac{1}{a}-1}}.
    \end{equation}  
This assumption can be justified by the fact that the main contribution to the integrals in $\eta$ comes from the sets where $\eta$ are of order one. Moreover, $\zeta$ is of order $1$ since we are assuming $\xi_2 \sim e^{\tau}$. Hence, we have that 
\begin{equation}
    \omega_0e^{-(1-a)\tau}= \frac{\eta^{\frac{1}{a}}}{(\eta-2\zeta)^{\frac{1}{a}-1}}e^{-(1-a)\tau} \rightarrow 0^+ \quad \text{as $\tau \to +\infty$.}
\end{equation}
Notice that the region $\eta \approx2 \zeta$ gives a negligible contribution since in this set $U(\zeta,\tau) \approx 0$. \\
Finally, setting $y=\left(\frac{\sigma}{\eta^a}\right)^{\frac{1}{1-a}}$ in \eqref{eq:Usigma} gives the final formula for $U(\zeta,\tau)$, namely
\begin{equation}\label{eq:Ufinal}
    U(\zeta,\tau) \sim \frac{16a}{\zeta}\int_{2\zeta}^{+\infty} \frac{\exp\left(-\frac{4 \pi \eta^{-a}}{1-a}\right)}{\eta^{2+a}}H\left(\frac{2\zeta}{\eta}\right) d \eta \Ii yU(y,(1-a)\tau)dy.
\end{equation}

The term containing $\log(y/\eta)$ has been neglected since the integrals give the main contribution for $\eta,y$ of order one. Hence, also $\log \left(y/\eta\right)$ is of order one and thus $\log(y/\eta) \ll \tau$ as $\tau \rightarrow + \infty$. It is convenient to define the functions
\begin{align} 
\lambda(\tau) & :=\Ii yU(y,\tau)dy \, ; \label{def:lambda}\\
\Phi(\zeta) &:= \frac{16 a}{\zeta}\int_{2\zeta}^{+\infty} d\eta  \frac{1}{\eta^{a+2}}  \exp\left(-\frac{4 \pi\eta^{-a}}{1-a}\right) H \left(\frac{2\zeta}{ \eta}\right) \,\label{eq:defPhi} .
\end{align}
Notice that $\Phi$ is bounded as $\zeta \to 0^{+}$ and by Lemma \ref{lem:asymHs} it follows that
\begin{equation}\label{eq:asymPhi}
\Phi(\zeta)\sim \frac {C}{\zeta^{2+a}} \log \zeta \quad \text{as}
\quad \zeta\to +\infty, \quad C>0,
\end{equation}
which in turns implies that $\lambda(\tau)$ is well-defined.
Therefore, \eqref{eq:Ufinal} can then be written in term of $\lambda,\Phi$ as
\begin{equation}\label{eq:U2}
    U(\zeta,\tau)=\Phi(\zeta)\lambda((1-a)\tau).
\end{equation}
From \eqref{eq:U2} it is possible to determine $\lambda(\tau)$. In fact, multiplying both sides of \eqref{eq:U2} by $\zeta$ and integrating over $[0,+\infty)$, we obtain
\begin{equation}\label{eq:conditioOnBeta}
    \lambda(\tau)=\left(\Ii \zeta 
    \Phi(\zeta) \, d \zeta \right) \lambda((1-a)\tau).
\end{equation}
This is a delay equation which is well-known to have solutions of the form $\lambda(\tau)=A\tau^{\beta}$ for some $A,\beta \in \Rl$ to be determined. In Section \ref{ssec:conclusion}, see \eqref{eq:valueA}, we will show that, in order to obtain conservation of mass, $A= \frac{M_0}{\log \left(\frac{1}{1-a}\right)}$, where $M_0$ is the initial mass of the solution $G$ to \eqref{eq:strongG}. Substituting $\lambda(\tau)=A\tau^{\beta}$ into \eqref{eq:conditioOnBeta} gives
\begin{equation}
    A\tau^{\beta}=\left(\Ii \zeta 
    \Phi(\zeta) \, d \zeta \right)(1-a)^{\beta}A\tau^{\beta}
\end{equation}
from which it follows that $\beta$ must be such that $(1-a)^{\beta}\Ii \zeta \Phi(\zeta) d\zeta = 1$ or more explicitly  
\begin{equation}\label{eq:EQbeta}
1=16a(1-a)^{\beta} \Ii d \zeta  \int_{2\zeta}^{+\infty}  \frac{\exp\left(-\frac{4 \pi \eta^{-a}}{1-a}\right)}{\eta^{a+2}}   H \left(\frac{2\zeta}{ \eta}\right) d\eta \, . 
\end{equation}
In order to determine $\beta$, we further introduce the function
\begin{equation}\label{def:Ka}
  K(a) =  \Ii d \zeta  \int_{2\zeta}^{+\infty}  \frac{1}{\eta^{a+2}} \,\exp\left(-\frac{4 \pi\eta^{-a}}{1-a}\right)   H \left(\frac{2\zeta}{ \eta}\right) d\eta .
\end{equation}
The function $K(a)$ can be computed explicitly; this is the content of the following lemma, whose proof is postponed to Section \ref{sec:auxiliary}.
\begin{lem}\label{lem:Ka}
    Let $K(a)$ be defined as in \eqref{def:Ka}. We have
    \begin{equation}\label{eq:Ka}
        K(a)=\frac{(1-a)}{16a}.
    \end{equation}
\end{lem}

From \eqref{eq:EQbeta}, using \eqref{eq:Ka} we obtain
\begin{align}
    1=16a(1-a)^{\beta} \frac{(1-a)}{16 a}
\end{align}
which implies 
\begin{equation}\label{eq:beta-1}
    \beta=-1.
\end{equation}

Finally, we have the following formula for $G(0,\xi_2,\tau)$
\begin{align}\label{eq:finalG_0}
    G(0,\xi_2,\tau) & \sim e^{-2 \tau}U(\xi_2e^{-\tau},\tau) \sim e^{-2\tau} \Phi(\xi_2 e^{-\tau})\lambda((1-a)\tau) \notag \\
    & = \frac{M_0 e^{-2\tau}}{(1-a) \log \left(\frac{1}{1-a}\right)\tau}\Phi(\xi_2e^{-\tau})\quad\text{as $\tau \rightarrow + \infty$}.
\end{align}

We now return to \eqref{eq:evGchar}, that is, the formula for $G$ obtained integrating by characteristics. Using
\eqref{eq:finalG_0} in \eqref{eq:evGchar}, as well as  the formula $\tau_0=\tau+\log\left(1-\frac{\xi_1}{\xi_2}\right)$, we arrive at 
\begin{align} \label{eq:GFinal}
    G(\xi_1,\xi_2,\tau)&=G\left(0,\frac{(\xi_2)^{\frac{1}{a}}}{(\xi_2-\xi_1)^{\frac{1}{a}-1}},\tau_0\right) \left(1-\frac{\xi_1}{\xi_2}\right)^{-\left(\frac{2}{a}-1\right)} \exp\left(-\frac{4 \pi\xi_1^{1-a}}{(1-a)\xi_2}\right)\notag \\
    &=  e^{-2\tau_0}\lambda((1-a)\tau_0)  \, \Phi\left( \frac{(\xi_2)^{\frac{1}{a}}}{(\xi_2-\xi_1)^{\frac{1}{a}-1}} e^{-\tau_0}
    \right) \exp\left(-\frac{4 \pi\xi_1^{1-a}}{(1-a)\xi_2}\right)\left(1-\frac{\xi_1}{\xi_2}\right)^{-\left(\frac{2}{a}-1\right)} \notag \\
    & = \frac{A e^{-2\tau_0}}{(1-a)\tau_0}  \, \Phi\left( \frac{(\xi_2)^{\frac{1}{a}}}{(\xi_2-\xi_1)^{\frac{1}{a}-1}} e^{-\tau_0}
    \right) \exp\left(-\frac{4 \pi\xi_1^{1-a}}{(1-a)\xi_2}\right)\left(1-\frac{\xi_1}{\xi_2}\right)^{-\left(\frac{2}{a}-1\right)} \notag\\
   & = \frac{M_0 e^{-2\tau_0}}{(1-a) \log \left(\frac{1}{1-a}\right)\tau_0}  \, \Phi\left( \frac{(\xi_2)^{\frac{1}{a}}}{(\xi_2-\xi_1)^{\frac{1}{a}-1}} e^{-\tau_0}
    \right) \exp\left(-\frac{4 \pi\xi_1^{1-a}}{(1-a)\xi_2}\right)\left(1-\frac{\xi_1}{\xi_2}\right)^{-\left(\frac{2}{a}-1\right)} \notag\\
    & = \frac{M_0 e^{-2\tau}}{(1-a) \log \left(\frac{1}{1-a}\right)(\tau+\log(1-\xi_1/\xi_2))} \; \Phi\left( \frac{(\xi_2)^{\frac{1}{a}+1}}{(\xi_2-\xi_1)^{\frac{1}{a}}} e^{-\tau}
    \right) \exp\left(-\frac{4 \pi\xi_1^{1-a}}{(1-a)\xi_2}\right)\left(1-\frac{\xi_1}{\xi_2}\right)^{-\left(\frac{2}{a}+1\right)} .
\end{align}

To conclude, we present an equivalent formula for $G$ which will be useful to better highlight where the mass is distributed. To this end we notice that
\begin{align}
    \Phi\left( \frac{(\xi_2)^{\frac{1}{a}+1}}{(\xi_2-\xi_1)^{\frac{1}{a}}} e^{-\tau}\right) & = \Phi \left(\frac{1}{\left(\frac{(\xi_2-\xi_1)e^{a\tau}}{\xi_2^{1+a}}\right)^{\frac{1}{a}}}\right); \notag \\
    e^{-2 \tau}\left(1-\frac{\xi_1}{\xi_2}\right)^{-\left(\frac{2}{a}+1\right)}& = \frac{1}{\left(\frac{(\xi_2-\xi_1)e^{a\tau}}{\xi_2^{1+a}}\right)^{\frac{2}{a}+1}}\frac{1}{\xi_2^{2+a}e^{-a\tau}}
\end{align}
and we set
\begin{equation}\label{eq:Z}
Z\left(\frac{(\xi_2-\xi_1)e^{a\tau}}{\xi_2^{1+a}}\right)= \Phi \left(\frac{1}{\left(\frac{(\xi_2-\xi_1)e^{a\tau}}{\xi_2^{1+a}}\right)^{\frac{1}{a}}}\right)\frac{1}{\left(\frac{(\xi_2-\xi_1)e^{a\tau}}{\xi_2^{1+a}}\right)^{\frac{2}{a}+1}}.
\end{equation}
Since $Z(s)=\Phi\left(s^{-{\frac{1}{a}}}\right)s^{-\left(\frac{2}{a}+1\right)}$, using \eqref{eq:defPhi},\eqref{eq:asymPhi} we obtain the following asymptotics for $Z$
\begin{equation}\label{eq:asymZ}
    \begin{cases}
        Z(s) \sim \frac{C_1}{\left(s^{-\frac{1}{a}}\right)^{2+a}} \log \left(s^{-\frac{1}{a}}\right)s^{-\left(\frac{2}{a}+1\right)} \sim -\frac{C_1}{a} \log s & \text{as $s \rightarrow 0^+$,} \\
        Z(s) \sim \frac{C_2}{s^{-\frac{1}{a}}}s^{-\left(\frac{2}{a}+1\right)} \sim \frac{C_2}{s^{1+\frac{1}{a}}} & \text{as $s \rightarrow + \infty$}
    \end{cases}
\end{equation}
with $C_1,C_2>0$. Hence we can write \eqref{eq:GFinal} in the equivalent form
\begin{equation}\label{eq:GfinalZ}
      G(\xi_1,\xi_2,\tau) \sim  \frac{M_0}{(1-a) \log \left(\frac{1}{1-a}\right) (\tau+\log(1-\xi_1/\xi_2))} \exp\left(-\frac{4\pi}{1-a}\frac{\xi_1^{1-a}}{\xi_2}\right) \frac{e^{a \tau}}{\xi_2^{2+a}}Z\left(\frac{(\xi_2-\xi_1)e^{a\tau}}{\xi_2^{1+a}}\right).
\end{equation}

\subsection{Conservation of mass for the function $G$}\label{ssec:conclusion}

In  this subsection,  we compute the mass of the function $G$ obtained in \eqref{eq:GFinal}. We recall that $\lambda(\tau)=\frac{A}{\tau}$. Imposing conservation of mass will determine the value of the constant $A$. Furthermore we will obtain that most of the mass of $G$ is in the region 
\begin{equation}\notag
    1-\frac{\xi_1}{\xi_2} \leq (\xi_2e^{-\tau})^a \quad \text{for $0 \lesssim \xi_2 \lesssim e^{\tau}$.}
\end{equation}
This means that the mass is spread across many \emph{size scales} of the velocity and it does not exhibit the typical behaviour of a standard self-similar solution as we explain in Section \ref{sec:noSelf}. Namely, if that were the case the mass would concentrate in a single scale, typically for $|\xi| \approx 1$. This behaviour is novel and has analogies with the toy model introduced to study simple shear deformations in \cite{JNV3}. To determine where the mass lies we start by considering the formula for $G$ obtained in \eqref{eq:GFinal}. To compute the mass in the region $0 \leq \xi_1 \leq \xi_2 \leq e^{\tau}$, using \eqref{eq:GFinal}, we have
\begin{align}
    &\int_0^{+\infty}  \int_0^{\xi_2}G(\xi_1,\xi_2,\tau)d \xi_2 d \xi_1 \notag \\ 
    & \sim    \frac{A e^{-2\tau}}{(1-a)} \int_0^{+\infty}\int_0^{\xi_2} \; \frac{1}{\tau+\log(1-\xi_1/\xi_2)} \Phi\left( \frac{(\xi_2)^{\frac{1}{a}+1}}{(\xi_2-\xi_1)^{\frac{1}{a}}} e^{-\tau}
    \right) \exp\left(-\frac{4 \pi\xi_1^{1-a}}{(1-a)\xi_2}\right)\left(1-\frac{\xi_1}{\xi_2}\right)^{-\left(\frac{2}{a}+1\right)} d \xi_1 d\xi_2.
\end{align}
Consider now the change of variables $\xi_1 \to \omega$  with $\omega= \frac{(\xi_2)^{\frac{1}{a}+1}}{(\xi_2-\xi_1)^{\frac{1}{a}}}e^{-\tau}$, where we have $d\xi_1=  ae^{\tau} \left(1-\frac{\xi_1}{\xi_2}\right)^{\frac{1}{a}+1}d\omega$. Then we obtain
\begin{align}
   \int_0^{+\infty} d\xi_2 \int_{0}^{\xi_2} &  G(\xi_1,\xi_2,\tau) d \xi_1\notag \\
   &= A\frac{a}{(1-a)} \int_0^{+\infty} \frac{ d \xi_2}{\xi_2} \int_{\xi_2e^{-\tau}}^{{+\infty}} \frac{\omega}{(1-a)\tau+a \log \left(\frac{\xi_2}{\omega}\right)}   \exp\left(-\frac{4 \pi}{1-a}\xi_2^{-a}\left(1-\left(\frac{\xi_2e^{-\tau}}{\omega}\right)^{a}\right)^{1-a}\right)  \Phi(\omega)d\omega .\label{eq:massG2} 
\end{align}
To estimate the last integral, fix $R>0$ large enough and consider the two cases $\xi_2 \in (R,+\infty)$ and $\xi_2 \in [0,R]$. In the first case we have that 
\begin{align}
    \Bigg|\exp\left(-\frac{4 \pi}{1-a}\xi_2^{-a}\left(1-\left(\frac{\xi_2e^{-\tau}}{\omega}\right)^{a}\right)^{1-a}\right) -1 \Bigg| \leq CR^{-a}, \quad C>0.
\end{align}
Therefore, we can write
\begin{align}\label{eq:I123}
    &A\frac{a}{(1-a)} \int_0^{+\infty} \frac{ d \xi_2}{\xi_2} \int_{\xi_2e^{-\tau}}^{{+\infty}} \frac{\omega}{(1-a)\tau+a \log \left(\frac{\xi_2}{\omega}\right)}   \exp\left(-\frac{4 \pi}{1-a}\xi_2^{-a}\left(1-\left(\frac{\xi_2e^{-\tau}}{\omega}\right)^{a}\right)^{1-a}\right)  \Phi(\omega)d\omega \notag \\
    &=A\frac{a}{(1-a)} \int_0^{R} \frac{ d \xi_2}{\xi_2} \int_{\xi_2e^{-\tau}}^{{+\infty}} \frac{\omega}{(1-a)\tau+a \log \left(\frac{\xi_2}{\omega}\right)}   \exp\left(-\frac{4 \pi}{1-a}\xi_2^{-a}\left(1-\left(\frac{\xi_2e^{-\tau}}{\omega}\right)^{a}\right)^{1-a}\right)  \Phi(\omega)d\omega \notag \\
    & + A\frac{a}{(1-a)} \int_R^{+\infty} \frac{ d \xi_2}{\xi_2} \int_{\xi_2e^{-\tau}}^{{+\infty}} \frac{\omega}{(1-a)\tau+a \log \left(\frac{\xi_2}{\omega}\right)}   \left(\exp\left(-\frac{4 \pi}{1-a}\xi_2^{-a}\left(1-\left(\frac{\xi_2e^{-\tau}}{\omega}\right)^{a}\right)^{1-a}\right)-1\right)  \Phi(\omega)d\omega \notag \\
    & + A\frac{a}{(1-a)} \int_R^{+\infty} \frac{ d \xi_2}{\xi_2} \int_{\xi_2e^{-\tau}}^{{+\infty}} \frac{\omega \Phi(\omega)}{(1-a)\tau+a \log \left(\frac{\xi_2}{\omega}\right)}  d\omega \notag \\
    & = I_1+I_2+I_3. 
\end{align}

Before estimating the integrals, a preliminary discussion is needed. Indeed, if $\xi_2 -\xi_1 \rightarrow 0^+$ the quantity $\tau_0=\tau+\log \left(1-\frac{\xi_1}{\xi_2}\right)$ becomes negative. This would result in a  non-integrable singularity at the region where $\tau_0$ vanishes. Notice however that the characteristics for which $\tau_0<0$ do not depart from the plane $\{\xi_1=0\}$ but from the initial time $\tau=0$. Therefore, if the initial datum $G_0$ is bounded, the contribution  from this region is exponentially small in $\tau$. In fact, if $\tau_0 \leq \delta \tau$, with $\delta>0$ small, we would have 
$$ \log \left(1-\frac{\xi_1}{\xi_2}\right) \leq -(1-\delta) \tau \quad \text{from which $1-\frac{\xi_1}{\xi_2} \leq e^{-(1-\delta)\tau}$}.$$
Moreover, the characteristic for $\xi_2$ in this region is $\xi_2(\tau)=\xi_{2,0}e^{- \left(\frac{1}{a}-1\right)(\tau-\tau_0)}$. Then if $\xi_2(\tau) \in [1,e^{\tau}]$ we would have $\xi_{2,0} \gg 1$ whereas in this case $G$ would increase exponentially since the characteristic curve comes from the region where $\xi_0$ is very large. Therefore, assuming that  $G_0(\xi_{0})$ decreases as a power law for $\xi_0\in[1, e^{\tau}]$ we would obtain that $G(\xi,\tau)$ tends to zero in the region where $\tau_0 \leq 0$.  
Finally, if $|\xi| \to 0$ we would obtain an additional decay due to the term 
$\exp \left(-\int_{\tau_0}^{\tau}|\xi_1(s)|^{-a}ds\right)$ that tends to zero as $\tau \rightarrow + \infty$. Therefore, with a slight abuse of the notation we will keep the notation in \eqref{eq:I123} understanding that the contribution due to the region where $\tau_0 \leq 0$ is negligible.  
We can then estimate $I_1$ as follows
\bigskip
\begin{align}
    I_1 & = A\frac{a}{(1-a)} \int_0^{R} \frac{ d \xi_2}{\xi_2} \int_{\xi_2e^{-\tau}}^{{+\infty}} \frac{\omega}{(1-a)\tau+a \log \left(\frac{\xi_2}{\omega}\right)}   \exp\left(-\frac{4 \pi}{1-a}\xi_2^{-a}\left(1-\left(\frac{\xi_2e^{-\tau}}{\omega}\right)^{a}\right)^{1-a}\right)  \Phi(\omega)d\omega \notag \\
    & \leq A\frac{C_a}{\tau} \int_0^R e^{-C \xi_2^{-a}} \frac{d \xi_2}{\xi_2} \int_{\xi_2 e^{-\tau}}^{+ \infty} \omega \Phi(\omega)d \omega + o(1) \notag \\
    &\leq A\frac{C_{a,R}}{\tau}+o(1)
\end{align}
where the remainder $o(1)$ contains the contribution coming from the regions with $|\xi_2-\xi_1| \leq e^{-(1-\delta)\tau}|\xi_2|$ and $|\omega-\xi_2 e^{-\tau}| \leq \frac{1}{2}|\xi_2e^{-\tau}|$. On the other hand, $I_2$ can be estimated as
\begin{align}
    I_2 &= A\frac{a}{(1-a)} \int_R^{+\infty} \frac{ d \xi_2}{\xi_2} \int_{\xi_2e^{-\tau}}^{{+\infty}} \frac{\omega}{(1-a)\tau+a \log \left(\frac{\xi_2}{\omega}\right)}   \left(\exp\left(-\frac{4 \pi}{1-a}\xi_2^{-a}\left(1-\left(\frac{\xi_2e^{-\tau}}{\omega}\right)^{a}\right)^{1-a}\right)-1\right)  \Phi(\omega)d\omega \notag \\
    & \leq A \frac{C_a}{\tau} \int_r^{+\infty} \frac{ d\xi_2}{\xi_2^{1+a}} \int_{\xi_2 e^{-\tau}}^{+\infty} \omega \Phi(\omega)d \omega \notag \\
    &\leq A\frac{C_{a,R}}{\tau}.
\end{align}

Finally, to estimate $I_3$ we must take into account the fact that mass of $G$ is distributed across several size scales of $\xi_2$. In fact, take $\delta >0$ small, we have that 

\begin{align}\label{eq:sumI_3}
    I_3 & =  A\frac{a}{(1-a)} \int_R^{+\infty} \frac{ d \xi_2}{\xi_2} \int_{\xi_2e^{-\tau}}^{{+\infty}} \frac{\omega \Phi(\omega)}{(1-a)\tau+a \log \left(\frac{\xi_2}{\omega}\right)}  d\omega \notag \\
    &=  A \frac{a}{(1-a)} \sum_{n= \lfloor R \rfloor}^{\lfloor\delta e^{\tau}\rfloor} \int_n^{n+1} \frac{d \xi_2}{\xi_2} \int_{\xi_2e^{-\tau}}^{{+\infty}} \frac{\omega \Phi(\omega)}{(1-a)\tau+a \log \left(\frac{\xi_2}{\omega}\right)}  d\omega \notag \\
    &+ A \frac{a}{(1-a)}\int_{\delta e^{\tau}}^{+\infty}\frac{d \xi_2}{\xi_2} \int_{\xi_2e^{-\tau}}^{{+\infty}} \frac{\omega \Phi(\omega)}{(1-a)\tau+a \log \left(\frac{\xi_2}{\omega}\right)}  d\omega. 
\end{align}
Using the asymptotic behaviour  \eqref{eq:asymPhi} for $\Phi$, the last term in the sum can be estimated as
\begin{align}
    A \frac{a}{(1-a)}\int_{\delta e^{\tau}}^{+\infty}\frac{d \xi_2}{\xi_2} \int_{\xi_2e^{-\tau}}^{{+\infty}} \frac{\omega \Phi(\omega)}{(1-a)\tau+a \log \left(\frac{\xi_2}{\omega}\right)}  d\omega & \leq A\frac{C_a}{\tau} \int_{\delta e^{\tau}}^{+\infty} \frac{d \xi_2}{\xi_2(\xi_2 e^{-\tau})^{\alpha}} \notag\\
    & = A\frac{C_a}{\tau} \int_{\delta}^{+\infty} \frac{dy}{y^{1+\alpha}} \notag \\
    &  \leq A \frac{C_{a,\alpha,\delta}}{\tau}
\end{align}
where $\alpha>0$. Moreover, recalling that the integral in $\omega$ has its main contribution for $\omega$ of order $1$, the contribution of the sum in \eqref{eq:sumI_3} is
\begin{align}\label{eq:sumApprox}
     &A \frac{a}{(1-a)} \sum_{n= \lfloor R \rfloor}^{\lfloor\delta e^{\tau}\rfloor} \int_n^{n+1} \frac{d \xi_2}{\xi_2} \int_{\xi_2e^{-\tau}}^{{+\infty}} \frac{\omega \Phi(\omega)}{(1-a)\tau+a \log \left(\frac{\xi_2}{\omega}\right)}  d\omega \notag \\
     & \sim A \frac{a}{(1-a)} \sum_{n= \lfloor R \rfloor}^{\lfloor\delta e^{\tau}\rfloor} \int_n^{n+1} \frac{d \xi_2}{\xi_2} \int_{\xi_2e^{-\tau}}^{{+\infty}} \frac{\omega \Phi(\omega)}{(1-a)\tau+a \log \xi_2}  d\omega \notag \\
     & \sim A \frac{a}{(1-a)}  \sum_{n= \lfloor R \rfloor}^{\lfloor\delta e^{\tau}\rfloor}  \frac{1}{n} \frac{1}{(1-a)\tau+a \log n} \int_{ne^{-\tau}}^{{+\infty}} \omega \Phi(\omega)  d\omega \notag \\
     & \approx A a \sum_{n= \lfloor R \rfloor}^{\lfloor\delta e^{\tau}\rfloor}  \frac{1}{n} \frac{1}{(1-a)\tau+a \log n}
\end{align}
where in the last line we have approximated $\int_{n e^{-\tau}}^{+\infty} \omega \Phi(\omega) d \omega \approx \Ii \omega \Phi(\omega) d \omega$ and we have used the fact that $\frac{1}{1-a}\Ii \omega \Phi(\omega) d \omega=1$ from \eqref{eq:EQbeta}. Approximating the sum as an integral we obtain
\begin{align}
    I_3 & \sim A a \sum_{n= \lfloor R \rfloor}^{\lfloor\delta e^{\tau}\rfloor}  \frac{1}{n} \frac{1}{(1-a)\tau+a \log n} \notag \\
    & \approx A a \, \frac{1}{a} \left[\log ((1-a)\tau+a \log (\delta e^{\tau}))-\log ((1-a)\tau)+a \log(R)\right] \notag \\
    & = A \log \left(\frac{\tau+a\log \delta}{(1-a)\tau+\log R}\right) \notag \\
    & \sim A \log \left(\frac{1}{1-a}\right) \quad \text{as $\tau \rightarrow + \infty$.}
\end{align}

In the end, letting $M_0$ be the initial mass of the solution $G$, we arrive at  
\begin{equation}\label{eq:valueA}
    M_0= \int_0^{+\infty}d\xi_2 \int_0^{\xi_1} G(\xi_1,\xi_2,\tau)d \xi_1 \sim A \log \left(\frac{1}{1-a}\right)
\end{equation}
from which, choosing $A=\frac{M_0}{\log \left(\frac{1}{1-a}\right)}$, we obtain the conservation of mass.

To conclude this subsection, we heuristically show that the asymptotic behaviour given by \eqref{eq:GfinalZ} better highlights the distribution of the mass. In fact, using \eqref{eq:GfinalZ} the mass behaves in a manner which is analogous to a mollifier of a Dirac delta $\delta(\xi_2-\xi_1)$ for $0 \lesssim\xi_2 \lesssim e^{\tau}$. To be precise, we recall that
\begin{equation}
      G(\xi_1,\xi_2,\tau) \sim  \frac{M_0}{(1-a) \log \left(\frac{1}{1-a}\right) (\tau+\log(1-\xi_1/\xi_2))} \exp\left(-\frac{4\pi}{1-a}\frac{\xi_1^{1-a}}{\xi_2}\right) \frac{e^{a \tau}}{\xi_2^{2+a}}Z\left(\frac{(\xi_2-\xi_1)e^{a\tau}}{\xi_2^{1+a}}\right).
\end{equation}
Formally, for $0 \lesssim \xi_2 \lesssim e^{\tau}$, we have that $\frac{e^{a \tau}}{\xi_2^{1+a}}Z\left(\frac{(\xi_2-\xi_1)e^{a \tau}}{\xi_2^{1+a}}\right)$ behaves like a Dirac delta in $\xi_1$, i.e.
\begin{align}\label{eq:DiracApprox}
  \frac{e^{a \tau}}{\xi_2^{1+a}}Z\left(\frac{(\xi_2-\xi_1)e^{a \tau}}{\xi_2^{1+a}}\right)\approx Z_0\delta(\xi_2-\xi_1) \quad \text{as $\tau \rightarrow + \infty$,}
\end{align}
with
\begin{align}
    Z_0= \Ii Z(s)ds = \Ii\Phi\left(s^{-\frac{1}{a}}\right)s^{-\left(\frac{2}{a}+1\right)}ds=a \Ii\Phi(y)ydy=a(1-a).
\end{align}
This can be seen by multiplying $ \frac{e^{a \tau}}{\xi_2^{1+a}}Z\left(\frac{(\xi_2-\xi_1)e^{a \tau}}{\xi_2^{1+a}}\right)$ by a test function $\vp(\xi_1)$ and integrating with respect to $\xi_1$, namely
\begin{equation*}
    \int_{0}^{+\infty} \frac{e^{a \tau}}{\xi_2^{1+a}}Z\left(\frac{(\xi_2-\xi_1)e^{a \tau}}{\xi_2^{1+a}}\right) \vp(\xi_1)d\xi_1 = -\int_0^{+\infty} Z(s) \vp\left(\frac{\xi_2^{1+a}}{e^{\tau}}s\right)ds.
\end{equation*}
Hence, taking the limit as $\tau \rightarrow + \infty$ and using the fact that $Z(s)$ is integrable gives
\begin{equation}\label{eq:DeltaZ}
     \int_{0}^{+\infty} \frac{e^{a \tau}}{\xi_2^{1+a}}Z\left(\frac{(\xi_2-\xi_1)e^{a \tau}}{\xi_2^{1+a}}\right) \vp(\xi_1)d\xi_1 \sim \left(\Ii Z(s)ds \right) \vp(0)= Z_0 \vp(0).
\end{equation}
\smallskip

Therefore, \eqref{eq:DeltaZ} implies that the mass lies in the region $0 \lesssim \xi_2 \lesssim e^{\tau}$ and from \eqref{eq:DiracApprox} we also obtain that for large $\tau$ we have that
\begin{equation}
    \frac{\xi_2-\xi_1}{\xi_2} \leq (\xi_2e^{-\tau})^a \ll 1
\end{equation}
from which it follows that $\xi_1 \approx \xi_2$ as $\tau \rightarrow + \infty$.

\subsection{Long-time behaviour for the solution $F$ of \eqref{eq:selfSimStrong2}}\label{ssec:longTimeF}
In this section we derive the asymptotic behaviour of the function $F$ by means of \eqref{eq:GFinal}. We recall that $F$ satisfies the equation
\begin{equation}\notag
    \partial_{\tau}F-\frac{1}{a}\partial_{\xi_1}(\xi_1F)-\left(\frac{1}{a}-1\right)\partial_{\xi_2}(\xi_2F)-\left(\frac{1}{a}-1\right)\partial_{\xi_3}(\xi_3F)
   + \xi_2 \partial_{\xi_1} F =
     \left(\mathcal{C}^+F\right)(\xi)-\left(\mathcal{C}^-F\right)(\xi)
\end{equation}
with
\begin{align}
 & \left(\mathcal{C}^+F\right)(\xi)= \frac{8 \delta(\xi_1)}{(2|\tilde \xi|)^{1+a}t^{1-a}}\int_0^{\pi}d\theta \frac{\sin(\theta)}{|\sin(2\theta)|^{1-a}}\int_{\Rp} F\left(\frac{2|\tilde \xi|}{t \sin(2\theta)},\eta, \tau\right)d\eta; \notag \\&
  \left(\mathcal{C}^-F\right)(\xi)=-4\pi |\xi_1|^{-a} F \notag . 
\end{align}
In analogy with the boundary problem \eqref{eq:Gboundary}, it is  possible to write a  boundary value problem also for the original solution $F$. More precisely, we consider the following
\begin{align}
     \partial_{\tau}F  -\frac{1}{a}\partial_{\xi_1}(\xi_1F) & -\left(\frac{1}{a}-1\right)\partial_{\xi_2}(\xi_2F)-\left(\frac{1}{a}-1\right)\partial_{\xi_3}(\xi_3F)
   + \xi_2 \partial_{\xi_1} F = -4 \pi |\xi_1|^{-a}; \notag \\
   F(0,\xi_2,\xi_3,\tau) &= \frac{8}{\xi_2(2|\tilde \xi|)^{1+a}t^{1-a}}\int_0^{\pi}d\theta \frac{\sin(\theta)}{|\sin(2\theta)|^{1-a}}\int_{\Rp} F\left(\frac{2|\tilde \xi|}{t \sin(2\theta)},\tilde \eta, \tau\right)d\tilde \eta \label{eq:BVPF}
\end{align}
and we notice that
\begin{equation}\notag
    \int_{\Rp} F\left(\frac{2|\tilde \xi|}{t \sin(2\theta)},\tilde \eta, \tau\right)d\tilde \eta = \int_{\Rl} G \left(\frac{2|\tilde \xi|}{t \sin(2\theta)},\eta, \tau\right)d\eta.
\end{equation}
In the same spirit of  Lemma \ref{lem:bdvalue}, we establish the following lemma for  $F(0,\xi_2,\xi_3,\tau)$, whose proof is postponed to Section \ref{sec:auxiliary}.
\begin{lem}\label{lem:R}
    Assume that $F(0,\xi_2,\xi_3,\tau)$ is the boundary value given in \eqref{eq:BVPF}. Then 
\begin{equation}\label{eq:BoundaryF_R}
    F(0,\xi_2,\xi_3,\tau)= \frac{8}{\xi_2(2|\tx|)^{1+a}t^{1-a}} \int_{0}^{+\infty} d\eta \int_{\frac{2|\tx|}{t}}^{+\infty} G(y,\eta,\tau) R\left(\frac{2|\tx|}{ty}\right) \frac{dy}{y}
\end{equation}
where $R:[0,1] \rightarrow \Rl$ is defined as
\begin{equation}\label{eq:R}
    R(s)=\frac{1}{\sqrt{2}}\frac{s^a}{\sqrt{1-s^2}}\left(\sqrt{1+\sqrt{1-s^2}}+\sqrt{1-\sqrt{1-s^2}}\right).
\end{equation}
\end{lem}
\smallskip
Combining Lemma \ref{lem:R} with the asymptotic behaviour for $G$ obtained in \eqref{eq:GFinal}
we obtain
\begin{align}
    F(0,\xi_2,\xi_3,\tau)  = \frac{8}{\xi_2(|2\tx|)^{1+a}t^{1-a}} \int_{\frac{2|\tx|}{t}}^{+\infty} d \eta \int_{\frac{2|\tx|}{t}}^{\eta} & \frac{dy}{y}  R \left(\frac{2|\tx|}{ty}\right) \frac{M_0 e^{-2\tau}}{(1-a) \log \left(\frac{1}{1-a}\right)(\tau+\log(1-y/\eta))}\notag \\
    & \times \Phi\left( \frac{(\eta)^{\frac{1}{a}+1}}{(\eta-y)^{\frac{1}{a}}} e^{-\tau}
    \right)  \exp\left(-\frac{4 \pi y^{1-a}}{(1-a)\eta}\right)\left(1-\frac{y}{\eta}\right)^{-\left(\frac{2}{a}+1\right)} .\notag
\end{align}
where we assumed $F=0$ if $y > \eta \geq 0$. Performing the changes of variables $y=\frac{2|\tx|}{t}s$ and $\eta=\frac{2|\tx|}{t}z$ we obtain
\begin{align}
    F(0,\xi_2,\xi_3,\tau) & = \frac{8M_0 e^{-2 \tau}}{(1-a) \log \left(\frac{1}{1-a}\right)\xi_2(|2\tx|)^{1+a}t^{1-a}} \int_1^{+\infty} \frac{2|\tx|}{t} \exp \left(-\frac{4 \pi s^{1-a}}{(1-a)z} \left(\frac{2|\tx|}{t}\right)^{-a}\right)dz \notag \\
    & \times\int_1^{z} \frac{1}{\tau + \log (1-s/z)} R\left(\frac{1}{s}\right)\Phi\left(\frac{z}{(1-s/z)^{\frac{1}{a}}}\frac{2|\tx|}{t^2}\right)\left(1-\frac{s}{z}\right)^{-\left(\frac{2}{a}+1\right)} \frac{ds}{s}. \notag
\end{align}
To obtain a simpler expression for $F(0,\xi_2,\xi_3,\tau)$, we consider the change of variables $ s\to \zeta$ given by
$$ \zeta= \frac{z}{(1-s/z)^{\frac{1}{a}}}\frac{2|\tx|}{t^2} \quad \text{with $d\zeta = \frac{1}{a}\frac{2|\tx|}{t^2} \frac{ds}{(1-s/z)^{\frac{1}{a}+1}} $} $$
from which it follows that
\begin{align}
     F(0,\xi_2,\xi_3,\tau) &= \frac{8aM_0}{(1-a) \log \left(\frac{1}{1-a}\right)\xi_2(|2\tx|)^{1+a}t^{1-a}} \frac{t^2}{2|\tx|}\int_1^{+\infty} \frac{2|\tx|}{t}\exp \left(-\frac{4 \pi( s(\zeta))^{1-a}}{(1-a)z} \left(\frac{2|\tx|}{t}\right)^{-a}\right)dz \notag \\
     &\times \int_{\zeta_0}^{z} \frac{1}{\tau+\log (s(\zeta)/z)}R\left( \frac{1}{s(\zeta)}\right) \Phi(\zeta) \left(1-\frac{s(\zeta)}{z}\right)^{-\left(\frac{2}{a}+1\right)}\frac{t^2}{2|\tx|}\left(1-\frac{s(\zeta)}{z}\right)^{\frac{1}{a}+1} \frac{d\zeta}{s(\zeta)}
\end{align}
with $\zeta_0= \frac{z}{(1-1/z)^{\frac{1}{a}}}\frac{2|\tx|}{t^2}$. Moreover, the change of variables $s \to \zeta$ also gives that 
\begin{align}
    &(1-s/z) = \left(\frac{z}{\zeta}\right)^a \left(\frac{2|\tx|}{t^2}\right)^a;\notag \\
    s  &= z \left(1-\left(\frac{z}{\zeta}\right)^a \left(\frac{2|\tx|}{t^2}\right)^a\right); \notag \\
    &(1-s/z)^{-\frac{1}{a}}  = \frac{\zeta}{z} \frac{t^2}{2|\tx|}. \notag 
\end{align}
Hence 
\begin{align}
     F(0,\xi_2,\xi_3,\tau) &= \frac{8aM_0 e^{-2 \tau}}{(1-a) \log \left(\frac{1}{1-a}\right)\xi_2(|2\tx|)^{1+a}t^{1-a}\tau} \int_1^{+\infty} \frac{2|\tx|}{t}\exp \left(-\frac{4 \pi z^{1-a}\left(1-\left(z/\zeta\right)^a\left(\frac{2|\tx|}{t^2}\right)^a\right)^{1-a}}{(1-a)z} \left(\frac{2|\tx|}{t}\right)^{-a}\right)dz \notag \\
     &\times \int_{\zeta_0}^{z}R\left( \frac{1}{s(\zeta)}\right) \Phi(\zeta) \frac{t^4}{4|\tx|^2} \frac{\zeta}{z}  \frac{1}{(1-a)\tau + a\log \left( \frac{z}{\zeta}\frac{2|\tx|}{t} \right) }\frac{d \zeta}{s(\zeta)} \notag \\
     & = \frac{8aM_0}{(1-a) \log \left(\frac{1}{1-a}\right)\xi_2(|2\tx|)^{1+a}t^{1-a}\tau} \frac{t}{2 |\tx|} \int_1^{+\infty}\exp \left(-\frac{4 \pi z^{-a}\left(1-\left(z/\zeta\right)^a\left(\frac{2|\tx|}{t^2}\right)^a\right)^{1-a}}{1-a} \left(\frac{2|\tx|}{t}\right)^{-a}\right)\frac{dz}{z} \notag \\
     &\times \int_{\zeta_0}^{z}R\left( \frac{1}{s(\zeta)}\right) \Phi(\zeta)\zeta\frac{1}{(1-a)\tau + a\log \left( \frac{z}{\zeta}\frac{2|\tx|}{t} \right) } \frac{d \zeta}{s(\zeta)}.
\end{align}

Since the integrals with respect to $z,\zeta$ give the main contribution for $z,\zeta$ of order one, in the set where $\frac{2|\tx|}{t}$ is of order one we have that $\zeta_0 \approx 0$, $ s \approx z$, $\frac{z}{\zeta}\frac{2|\tx |}{t^2} \approx 0$ and $\log \left(\frac{2|\tx|}{t^2}\right) =\log \left(\frac{2|\tx|}{t}\right)- \log t \sim -\tau$ as $t \rightarrow + \infty$. This gives, as $\tau \to \infty$, 
\begin{align}
    F(0,\xi_2,\xi_3,\tau)  & \sim\frac{8aM_0}{ \log \left(\frac{1}{1-a}\right)\xi_2(|2\tx|)^{2+a}t^{-a}(1-a)^2\tau}  \int_1^{+\infty}R\left( \frac{1}{z}\right)\exp \left(-\frac{4 \pi z^{-a}}{1-a} \left(\frac{2|\tx|}{t}\right)^{-a}\right)\frac{dz}{z^2} \int_{0}^{+\infty} \zeta\Phi(\zeta) d \zeta \notag \\
    & \approx  \frac{8aM_0}{(1-a)\log\left(\frac{1}{1-a}\right)\xi_2(|2\tx|)^{2+a}t^{-a}\tau}  \int_1^{+\infty}R\left( \frac{1}{z}\right)\exp \left(-\frac{4 \pi z^{-a}}{1-a} \left(\frac{2|\tx|}{t}\right)^{-a}\right)\frac{dz}{z^2} \label{eq:finalF0}
\end{align}
where we have used the fact that $ \frac{1}{1-a}\Ii \zeta \Phi(\zeta)d\zeta=1$ from \eqref{eq:EQbeta}. \\

We can now describe the asymptotic profile of $F$. From \eqref{eq:selfSimStrong2} the characteristics read
\begin{align}\label{eq:characteristicsF}
    \frac{d \xi_1}{ds} &=-\frac{1}{a}\xi_1+\xi_2; \\
    \frac{d\xi_2}{ds} & = -\left(\frac{1}{a}-1\right)\xi_2; \notag \\
    \frac{d\xi_3}{ds} & = -\left(\frac{1}{a}-1\right)\xi_3; \notag \\
    \frac{d F}{ds} &=\left(\frac{3}{a}-2\right)F-4 \pi|\xi_1(s)|^{-a}F
\end{align}
which, taking $\xi_{1,0}=0,\xi_2(0)=\xi_{2,0},\xi_3(0)=\xi_{3,0}$ as initial conditions, are solved by
\begin{align}
    \xi_1(\tau) &=\xi_{2,0}e^{-\frac{1}{a}(\tau-\tau_0)}(e^{\tau-\tau_0}-1); \label{eq:xi1F} \\
    \xi_2(\tau) & =\xi_{2,0}e^{-\left(\frac{1}{a}-1\right)(\tau-\tau_0)}; \label{eq:xi2F} \\
    \xi_3(\tau) & =\xi_{3,0}e^{-\left(\frac{1}{a}-1\right)(\tau-\tau_0)}; \label{eq:xi3F} \\
    F(\xi_1,\xi_2,\xi_3,\tau) & =F(0,\xi_{2,0},\xi_{3,0},\tau_0)\exp\left(\left(\frac{3}{a}-2\right)(\tau-\tau_0)\right)\exp \left(-4 \pi\int_{\tau_0}^{\tau}|\xi_1(s)|^{-a}\right) \label{eq:charF}.
\end{align}
Moreover from \eqref{eq:lossCharacteristics} we obtain 
\begin{equation}
F(\xi_1,\xi_2,\xi_3,\tau)  =F\left(0,\xi_{2,0},\xi_{3,0},\tau_0\right)\exp \left(-\frac{4 \pi \,\xi_1^{1-a}}{(1-a)\xi_2}\right)\left(1-\frac{\xi_1}{\xi_2}\right)^{-\left(\frac{3}{a}-2\right)}, 
\end{equation}
thus, using \eqref{eq:characteristicsF}, we also get 
\begin{equation}\label{eq:characteristicsFIntermediate}
F(\xi_1,\xi_2,\xi_3,\tau)  =F\left(0,\frac{(\xi_2)^{\frac{1}{a}}}{(\xi_2-\xi_1)^{\frac{1}{a}-1}},\frac{(\xi_2)^{\frac{1}{a}-1}\xi_3}{(\xi_2-\xi_1)^{\frac{1}{a}-1}},\tau_0\right)\exp \left(- \frac{4 \pi \, \xi_1^{1-a}}{(1-a)\xi_2}\right)\left(1-\frac{\xi_1}{\xi_2}\right)^{-\left(\frac{3}{a}-2\right)}.
\end{equation}

Therefore, from \eqref{eq:finalF0} we obtain
\begin{align}
    F(\xi_1,\xi_2,\xi_3,\tau) & = \frac{8aM_0}{(1-a) \log \left(\frac{1}{1-a}\right)\frac{(\xi_2)^{\frac{1}{a}}}{(\xi_2-\xi_1)^{\frac{1}{a}-1}}\frac{\xi_2^{\frac{2}{a}-1-a}}{(\xi_2-\xi_1)^{\frac{2}{a}-1-a}}(2|\tx|)^{2+a}e^{-a\tau_0}\tau_0}\exp \left(- \frac{4 \pi \, \xi_1^{1-a}}{(1-a)\xi_2}\right)\left(1-\frac{\xi_1}{\xi_2}\right)^{-\left(\frac{3}{a}-2\right)}  \notag \\
    &\times \int_1^{+\infty} R\left(\frac{1}{z}\right) \exp \left(-\frac{4 \pi z^{-a}}{1-a} \left(\frac{\xi_2}{(\xi_2-\xi_1)}\right)^{-(1-a)}\left(\frac{2|\tx|}{e^{\tau_0}}\right)^{-a}\right) \frac{dz}{z^2} \notag \\
    & = \frac{8aM_0}{(1-a) \log \left(\frac{1}{1-a}\right)\left(\frac{\xi_2}{\xi_2-\xi_1}\right)^{-a}\xi_2(2|\tx|)^{2+a}e^{-a\tau_0}\tau_0}\exp \left(- \frac{4 \pi \, \xi_1^{1-a}}{(1-a)\xi_2}\right)  \notag \\
    &\times \int_1^{+\infty} R\left(\frac{1}{z}\right) \exp \left(-\frac{4 \pi z^{-a}}{1-a} \left(\frac{\xi_2}{(\xi_2-\xi_1)}\right)^{-(1-a)}\left(\frac{2|\tx|}{e^{\tau_0}}\right)^{-a}\right) \frac{dz}{z^2}.\notag
\end{align}
Moreover using the formula $\tau_0=\tau+\log(1-\xi_1/\xi_2)$ in the equation above yields
\begin{align}
    F(\xi_1,\xi_2,\xi_3,\tau) & =  \frac{8aM_0}{(1-a) \log \left(\frac{1}{1-a}\right)\left(\frac{\xi_2}{\xi_2-\xi_1}\right)^{-a}\xi_2(2|\tx|)^{2+a}e^{-a\tau}(1-\xi_1/\xi_2)^{-a}(\tau+\log(1-\xi_1/\xi_2))}\exp \left(- \frac{4 \pi \, \xi_1^{1-a}}{(1-a)\xi_2}\right)  \notag \\
    &\times \int_1^{+\infty} R\left(\frac{1}{z}\right) \exp \left(-\frac{4 \pi z^{-a}}{1-a} \left(\frac{\xi_2}{(\xi_2-\xi_1)}\right)^{-(1-a)}\left(1-\frac{\xi_1}{\xi_2}\right)^a\left(\frac{2|\tx|}{e^{\tau}}\right)^{-a}\right) \frac{dz}{z^2}\notag \\
    & =  \frac{8aM_0}{(1-a) \log \left(\frac{1}{1-a}\right)\xi_2(2|\tx|)^{2+a}e^{-a\tau}(\tau+\log(1-\xi_1/\xi_2))}\exp \left(- \frac{4 \pi \, \xi_1^{1-a}}{(1-a)\xi_2}\right)  \notag \\
    &\times \int_1^{+\infty} R\left(\frac{1}{z}\right) \exp \left(-\frac{4 \pi z^{-a}}{1-a}\left(1-\frac{\xi_1}{\xi_2}\right)\left(\frac{2|\tx|}{e^{\tau}}\right)^{-a}\right) \frac{dz}{z^2}.\label{eq:finalF}
\end{align}

We now check that we recover formula \eqref{eq:GfinalZ} from $\int_{\Rl}F(\xi_1,\xi_2,\xi_3,\tau)d\xi_3$. Using \eqref{eq:finalF} we have 
\begin{align}
    \int_{\Rl}F(\xi_1,\xi_2,\xi_3,\tau) d\xi_3 &=  \frac{8aM_0}{(1-a) \log \left(\frac{1}{1-a}\right)\xi_2 e^{-a\tau}(\tau+\log(1-\xi_1/\xi_2))}\exp \left(- \frac{4 \pi \, \xi_1^{1-a}}{(1-a)\xi_2}\right)  \notag \\
    &\times \int_{\Rl} \frac{d\xi_3}{(2|\tx|)^{2+a}} \int_1^{+\infty} R\left(\frac{1}{z}\right) \exp \left(-\frac{4 \pi z^{-a}}{1-a}\left(1-\frac{\xi_1}{\xi_2}\right)\left(\frac{2|\tx|}{e^{\tau}}\right)^{-a}\right) \frac{dz}{z^2} \notag \\
    & =\frac{8aM_0}{(1-a) \log \left(\frac{1}{1-a}\right)\xi_2e^{-a\tau}(\tau+\log(1-\xi_1/\xi_2))}\exp \left(- \frac{4 \pi \, \xi_1^{1-a}}{(1-a)\xi_2}\right)  \notag \\
    &\times \int_{\Rl} \frac{d\xi_3}{(2|\tx|)^{2+a}} \int_0^{1} R\left(z\right) \exp \left(-\frac{4 \pi }{1-a}\left(1-\frac{\xi_1}{\xi_2}\right)\left(\frac{2|\tx|}{zt}\right)^{-a}\right) dz  \notag \\
    & =\frac{16aM_0}{(1-a) \log \left(\frac{1}{1-a}\right)\xi_2e^{-a\tau}(\tau+\log(1-\xi_1/\xi_2))}\exp \left(- \frac{4 \pi \, \xi_1^{1-a}}{(1-a)\xi_2}\right)  \notag \\
    &\times \Ii \frac{d\xi_3}{(2|\tx|)^{1+a}} \int_{\frac{2|\tx|}{t}}^{+\infty} R\left(\frac{2|\tx|}{yt}\right) \exp \left(-\frac{4 \pi }{1-a}\left(1-\frac{\xi_1}{\xi_2}\right)y^{-a}\right) \frac{dy}{ty^2}
\end{align}
where in the last line we have used the change of variables $y=\frac{2|\tx|}{zt}$. Using Fubini to exchange the order of integration we obtain

\begin{align}
     \int_{\Rl}F(\xi_1,\xi_2,\xi_3,\tau) d\xi_3 & =\frac{16aM_0}{(1-a) \log \left(\frac{1}{1-a}\right)\xi_2e^{-a\tau}(\tau+\log(1-\xi_1/\xi_2))}\exp \left(- \frac{4 \pi \, \xi_1^{1-a}}{(1-a)\xi_2}\right)  \notag \\
    &\times \int_{\frac{2|\xi_2|}{t}}^{+\infty} \exp \left(-\frac{4 \pi }{1-a}\left(1-\frac{\xi_1}{\xi_2}\right)y^{-a}\right) \frac{dy}{ty} \int_{0}^{\sqrt{\left(\frac{yt}{2}\right)^2-\xi_2^2}} R\left(\frac{2|\tx|}{yt}\right)\frac{1}{y} \frac{d\xi_3}{(2|\tx|)^{1+a}} \notag \\
     & =\frac{16aM_0}{(1-a) \log \left(\frac{1}{1-a}\right)\xi_2e^{-a\tau}(\tau+\log(1-\xi_1/\xi_2))}\exp \left(- \frac{4 \pi \, \xi_1^{1-a}}{(1-a)\xi_2}\right)  \notag \\
    &\times \int_{\frac{2|\xi_2|}{t}}^{+\infty} \exp \left(-\frac{4 \pi }{1-a}\left(1-\frac{\xi_1}{\xi_2}\right)y^{-a}\right) H\left(\frac{2\xi_2}{ty}\right) \frac{dy}{y(ty)^{1+a}}
\end{align}
where the last equality follows by the proof of Lemma \ref{lem:bdvalue}. Now changing variable as $y=\left(1-\frac{\xi_1}{\xi_2}\right)^{\frac{1}{a}}x$ inside the integral we get 
\begin{align}
    \int F(\xi_1,\xi_2,\xi_3,\tau)d\xi_3 & =\frac{16aM_0}{(1-a) \log \left(\frac{1}{1-a}\right)\xi_2t(\tau+\log(1-\xi_1/\xi_2))}\exp \left(- \frac{4 \pi \, \xi_1^{1-a}}{(1-a)\xi_2}\right)  \notag \\
    &\times \left(1-\frac{\xi_1}{\xi_2}\right)^{-\left(\frac{1}{a}+1\right)} \int_{\bar x}^{+\infty} \exp \left(-\frac{4 \pi }{1-a}x^{-a}\right) H\left(\frac{2\xi_2}{\left(1-\frac{\xi_1}{\xi_2}\right)^{\frac{1}{a}}xt}\right) \frac{dx}{x^{2+a}}
\end{align}
with $\bar x= \left(1-\frac{\xi_1}{\xi_2}\right)^{-\frac{1}{a}}\frac{2\xi_2}{t}$. 
Moreover, using \eqref{eq:defPhi}-\eqref{eq:Z}, we obtain
\begin{align}
    &\int_{\Rl} F(\xi_1,\xi_2,\xi_3,\tau)d\xi_3  = \frac{16aM_0}{(1-a) \log \left(\frac{1}{1-a}\right)(\tau + \log (1-\xi_1/\xi_2))} \frac{t^{a}}{\xi_2^{2+a}}\exp \left(- \frac{4 \pi \, \xi_1^{1-a}}{(1-a)\xi_2}\right) \notag \\
    & \times \left(\frac{\xi_2-\xi_1}{\xi_2^{1+a}}\right)^{\frac{1}{a}} \bigintsss_{\frac{2}{\left(\frac{\xi_2-\xi_1}{\xi_2^{1+a}}t^a\right)^{\frac{1}{a}}}}^{+\infty} \exp \left(-\frac{4\pi}{1-a}x^{-a}\right)H\left(\frac{2}{\left(\frac{\xi_2-\xi_1}{\xi_2^{1+a}}t^a\right)^{\frac{1}{a}}x}\right) \frac{dx}{x^{2+a}} \notag \\
    & = \frac{M_0}{(1-a) \log \left(\frac{1}{1-a}\right)(\tau + \log (1-\xi_1/\xi_2))}  \frac{e^{a\tau}}{\xi_2^{2+a}} \exp \left(- \frac{4 \pi \, \xi_1^{1-a}}{(1-a)\xi_2}\right)\Phi \left(\frac{1}{\left(\frac{(\xi_2-\xi_1)e^{a\tau}}{\xi_2^{1+a}}\right)^{\frac{1}{a}}}\right)\frac{1}{\left(\frac{(\xi_2-\xi_1)e^{a\tau}}{\xi_2^{1+a}}\right)^{\frac{2}{a}+1}} \notag \\
    & = \frac{M_0}{(1-a) \log \left(\frac{1}{1-a}\right)(\tau + \log (1-\xi_1/\xi_2))}\frac{e^{a\tau}}{\xi_2^{2+a}} Z\left(\frac{(\xi_2-\xi_1)e^{a\tau}}{\xi_2^{1+a}}\right)\exp \left(- \frac{4 \pi \, \xi_1^{1-a}}{(1-a)\xi_2}\right)
\end{align}

Finally we obtain
\begin{align}
   G(\xi_1,\xi_2,\tau) & = \int_{\Rl} F(\xi_1,\xi_2,\xi_3,\tau) d \xi_3 
   \notag \\
   & = \frac{M_0}{(1-a) \log \left(\frac{1}{1-a}\right)(\tau+\log(1-\xi_1/\xi_2))}\frac{e^{a\tau}}{\xi_2^{2+a}} Z\left(\frac{(\xi_2-\xi_1)e^{a\tau}}{\xi_2^{1+a}}\right)\exp \left(- \frac{4 \pi \, \xi_1^{1-a}}{(1-a)\xi_2}\right).
\end{align}

\subsection{Proof of the auxiliary lemmas} \label{sec:auxiliary}

In this section, we prove the technical results used in the argument presented in Section \ref{sec:strategy}. 
We first provide the detailed proof that \eqref{eq:selfSimStrong2} preserves the mass via direct computations. \smallskip

\noindent\textbf{Mass conservation  for \eqref{eq:selfSimStrong2}:} We have
\begin{align*}
    \int d\xi_1 \int d\tilde\xi \ & \left(\mathcal{C}^+F\right)(\xi)= \int d\xi_1 \int d\tilde\xi \ \left[ \frac{8 \delta(\xi_1)}{(2|\tilde \xi|)^{1+a}t^{1-a}}\int_0^{\pi}d\theta \frac{\sin(\theta)}{|\sin(2\theta)|^{1-a}}\int_{\Rp} F\left(\frac{2|\tilde \xi|}{t \sin(2\theta)},\eta, \tau\right)d\eta \right]
   \\&
   = \frac{8 }{t^{1-a}} \int d\xi_1 \delta(\xi_1) \int_0^{\pi}d\theta \frac{\sin(\theta)}{|\sin(2\theta)|^{1-a}}\int d\tilde\xi  \frac 1 {(2|\tilde \xi|)^{1+a}} \int_{\Rp} F\left(\frac{2|\tilde \xi|}{t \sin(2\theta)},\eta, \tau\right)d\eta\\&
   = \frac{8 }{t^{1-a}} \int d\xi_1 \delta(\xi_1) \int_0^{\pi}d\theta \frac{\sin(\theta)}{|\sin(2\theta)|^{1-a}}\int_0^{+\infty} d(|\tilde\xi|)\int_0^{2\pi} d\alpha  \frac {|\tilde\xi| } {(2|\tilde \xi|)^{1+a}} \int_{\Rp} F\left(\frac{2|\tilde \xi|}{t \sin(2\theta)},\eta, \tau\right)d\eta \\&
   = \frac{4 \cdot (2 \pi)}{t^{1-a}} \int d\xi_1 \delta(\xi_1) \int_0^{\pi}d\theta \frac{\sin(\theta)}{|\sin(2\theta)|^{1-a}}\int_0^{+\infty} d(|\tilde\xi|)  \frac {1 } {(2|\tilde \xi|)^{a}} \int_{\Rp} F\left(\frac{2|\tilde \xi|}{t \sin(2\theta)},\eta, \tau\right)d\eta\\&
   = \frac{4 \cdot (2 \pi) t}{t^{1-a}t^a}  \int_0^{\frac \pi 2}d\theta \frac{\sin(\theta)}{|\sin(2\theta)|^{1-a}}\int_0^{+\infty} d\sigma \, \frac{|\sin(2\theta)|}{2} \frac {1 } {( |\sin(2\theta)|)^{a}} \frac {1 } {\sigma^{a}} \int_{\Rp} F\left(\sigma,\eta, \tau\right)d\eta \\& 
   \quad + \frac{4 \cdot (2 \pi) t}{t^{1-a}t^a} \int_{\frac \pi 2}^\pi d\theta \frac{\sin(\theta)}{|\sin(2\theta)|^{1-a}}\int_{-\infty}^0 d\sigma \, \frac{|\sin(2\theta)|}{2} \frac {1 } {( |\sin(2\theta)|)^{a}} \frac {1 } {\sigma^{a}} \int_{\Rp} F\left(\sigma,\eta, \tau\right)d\eta\\& 
   = 4 \pi  \int_0^{\frac \pi 2}d\theta \sin(\theta)  \left[\int_0^{+\infty} d\sigma \, \frac {1 } {\sigma^{a}} \int_{\Rp} F\left(\sigma,\eta, \tau\right)d\eta  + 
   \int_{-\infty}^0 d\sigma \,   \frac {1 } {(-\sigma)^{a}} \int_{\Rp} F\left(\sigma,\eta, \tau\right)d\eta \right]
   \\& 
      = 4 \pi    \int_{\mathbb{R}}  d\sigma  \int_{\Rp} d\eta \, \frac {1 } {|\sigma|^{a}} F\left(\sigma,\eta, \tau\right)=\int d\xi_1 \int d\tilde\xi \  \left(\mathcal{C}^- F\right)(\xi)  
\end{align*}
where the gain and loss operators $\mathcal{C}^+, \, \mathcal{C}^-$ are defined in \eqref{eq:selfSimStrong2_operators}-\eqref{eq:selfSimStrong2_operatorsLoss}.
Notice that, in the computation above,  in the second identity we performed the change of variables to polar coordinates, namely $\tilde \xi=(\xi_2,\xi_3)=|\tilde \xi|(\cos\alpha,\sin\alpha)$  with Jacobian $|\tilde \xi|$. 
Moreover, in the fourth identity we set $\sigma=\frac{2|\tilde \xi|}{t \sin(2\theta)}$ with $d|\tilde \xi|=\frac{1}{2} t |\sin(2\theta)|d\sigma$, so that $|\tilde \xi|=\frac{1}{2} t \sin(2\theta)\sigma$ if $\theta < \frac \pi 2$ and $\sigma>0$ while $|\tilde \xi|=- \frac{1}{2} t \sin(2\theta)\sigma$ if $  \frac \pi 2 < \theta< \pi$ and $\sigma>0$.  
 We also used that $ \int_{-\infty}^0  d\sigma \dots +\int_0^{+\infty} d\sigma \dots=1$ and $\int_0^{\frac \pi 2}d\theta \sin(\theta) =\int_{\frac \pi 2}^\pi d\theta \sin(\theta) =1. $ Therefore, since the integrals of the derivatives with respect to $\xi$ vanish, we obtain the conservation of mass property for  \eqref{eq:selfSimStrong2}.
\medskip

We now prove the Lemmas used in Section \ref{sec:strategy}. We start with the proof of Lemma \ref{lem:bdvalue} and the proof of Lemma \ref{lem:asymHs}. 

\begin{proofof}[Proof of Lemma \ref{lem:bdvalue}] 
We consider the boundary value $G(0,\xi_2,\tau) $ given as in \eqref{eq:Gboundary}, i.e.
\begin{align*} 
    G(0,\xi_2,\tau) & = \frac{16}{\xi_2t^{1-a}}\int_0^{\pi}d\theta \frac{\sin(\theta)}{|\sin(2\theta)|^{1-a}}\Ii \frac{d\xi_3}{(2|\tilde \xi|)^{1+a}}\int_{\Rl} G\left(\frac{2|\tilde \xi|}{t \sin(2\theta)},\eta, \tau\right)d\eta\notag \\& 
    = \frac{32}{\xi_2 t^{1-a}}\int_0^{\frac{\pi}{2}}d\theta \frac{\sin(\theta)}{|\sin(2\theta)|^{1-a}}\Ii \frac{d\xi_3}{(2|\tilde \xi|)^{1+a}}\int_{\Rl} G\left(\frac{2|\tilde \xi|}{t \sin(2\theta)},\eta, \tau\right)d\eta 
\end{align*}
    and rewrite the integral in the $\theta$ variable by setting $\psi=2\theta$. Hence 
   \begin{align} \label{eq:Gbound1}
    G(0,\xi_2,\tau)  = \frac{16}{\xi_2t^{1-a}}\int_0^{\pi}d\psi \frac{\sin(\frac{\psi}{2})}{|\sin(\psi)|^{1-a}}\Ii \frac{d\xi_3}{(2|\tilde \xi|)^{1+a}}\int_{\Rl} G\left(\frac{2|\tilde \xi|}{t \sin(\psi)},\eta, \tau\right)d\eta \, .
\end{align} 
    We divide the integral above separating the cases $\psi \in [0,\pi/2], \psi \in [\pi/2,\pi]$, namely 
$$ 
    \int_0^{\pi}d\psi \frac{\sin(\frac{\psi}{2})}{|\sin(\psi)|^{1-a}}= \frac{1}{2}\int_0^{\pi/2}d \psi_1 \frac{\sin(\psi_1/2)}{|\sin(\psi_1)|^{1-a}} + \frac{1}{2}\int_{\pi/2}^{\pi}d \psi_2 \frac{\sin(\psi_2/2)}{|\sin(\psi_2)|^{1-a}}\, 
$$
where $\psi_1\in [0,\pi/2]$ and $\psi_2\in [\pi/2,\pi]$. 
Performing also the change of variables $y_j = \frac{2|\tx|}{t \sin \psi_j}$ for $j=1,2$, with $dy_j=- \frac{2|\tx|}{t \sin^2 \psi_j}\cos\psi_j \, d\psi_j$ we arrive at  
\begin{align}\label{eq:splitpsi}
    \frac{1}{2} \int_{\frac{2|\tx|}{t}}^{+\infty}  \frac{\sin(\psi_1/2)}{|\sin \psi_1|^{1-a}} \frac{t \sin^2 \psi_1}{2|\tx|\cos \psi_1} dy_1 & + \frac{1}{2} \int_{\frac{2|\tx|}{t}}^{+\infty}  \frac{\sin(\psi_2/2)}{|\sin \psi_2|^{1-a}} \frac{t \sin^2 \psi_2}{2|\tx|\cos \psi_2} dy_2 = \notag \\
    \frac{1}{2} \int_{\frac{2|\tx|}{t}}^{+\infty}  \frac{\sin(\psi_1/2)}{|\sin \psi_1|^{1-a}} \frac{4|\tx|^2}{t^2 y_1^2}\frac{t}{2 |\tx| } \frac{1}{\sqrt{1-\sin^2 \psi_1}} dy_1 & + \frac{1}{2} \int_{\frac{2|\tx|}{t}}^{+\infty}  \frac{\sin(\psi_2/2)}{|\sin \psi_2|^{1-a}} \frac{4|\tx|^2}{t^2 y_2^2}\frac{t}{2 |\tx| } \frac{1}{\sqrt{1-\sin^2 \psi_2}} dy_2 \notag \\
    \int_{\frac{2|\tx|}{t}}^{+\infty}  \frac{\sin(\psi_1/2)}{|\sin \psi_1|^{1-a}} \frac{|\tx|}{ty_1^2} \frac{1}{\sqrt{1-\frac{4|\tx|^2}{t^2y_1^2}}}dy_1 & + \int_{\frac{2|\tx|}{t}}^{+\infty}  \frac{\sin(\psi_2/2)}{|\sin \psi_2|^{1-a}} \frac{|\tx|}{ty_2^2} \frac{1}{\sqrt{1-\frac{4|\tx|^2}{t^2y_2^2}}}dy_2\, .
\end{align}
We remark that $\psi_1,\psi_2$ depend on $y_{1,2},|\tx|,t$ since $\sin\psi_j=\frac{2|\tx|}{t y_j}, \, j=1,2$, $\psi_1\in [0,\pi/2],\,\psi_2\in [\pi/2,\pi]$. 

Therefore, we now integrate with respect to the $\eta$ variable (cf. \eqref{eq:Gbound1}) and, specifically, we consider
$$ \int_0^{\pi}d\psi \frac{\sin(\frac{\psi}{2})}{|\sin(\psi)|^{1-a}}\int_{\Rp} G\left(\frac{2|\tilde \xi|}{t \sin(\psi)},\eta, \tau\right)d\eta. $$
From \eqref{eq:splitpsi} we then obtain
\begin{align}\label{eq:Gbound2}
    & \int_0^{+\infty} d\eta \int_{\frac{2|\tx|}{t}}^{+\infty}  \frac{\sin(\psi_1/2)}{|\sin \psi_1|^{1-a}} \frac{|\tx|}{ty_1^2} \frac{1}{\sqrt{1-\frac{4|\tx|^2}{t^2y_1^2}}} G(y_1,\eta,\tau)dy_1  + \int_0^{+\infty} d\eta\int_{\frac{2|\tx|}{t}}^{+\infty}  \frac{\sin(\psi_2/2)}{|\sin \psi_2|^{1-a}} \frac{|\tx|}{ty_2^2} \frac{1}{\sqrt{1-\frac{4|\tx|^2}{t^2y_2^2}}}G(y_2,\eta,\tau)dy_2  \notag \\
   & = \frac{1}{2}\int_0^{+\infty} d\eta \int_{\frac{2|\tx|}{t}}^{+\infty}  \frac{\sin(\psi_1/2)}{y_1} \left(\frac{2 |\tx|}{ty_1}\right)^{a} \frac{1}{\sqrt{1-\frac{4|\tx|^2}{t^2 y_1^2}}}  G(y_1,\eta,\tau)dy_1 \notag \\
   & +\frac{1}{2}\int_0^{+\infty} d\eta \int_{\frac{2|\tx|}{t}}^{+\infty}  \frac{\sin(\psi_2/2)}{y_2} \left(\frac{2 |\tx|}{ty_2}\right)^{a} \frac{1}{\sqrt{1-\frac{4|\tx|^2}{t^2y_2^2}}}  G(y_2,\eta,\tau)dy_2 \notag \\
  & = \frac{1}{2}\int_0^{+\infty} d\eta \int_{\frac{2|\tx|}{t}}^{+\infty}  \frac{(\sin(\psi_1/2)  + \sin(\psi_2/2))}{y} \left(\frac{2 |\tx|}{ty}\right)^{a} \frac{1}{\sqrt{1-\frac{4|\tx|^2}{t^2y^2}}}  G(y,\eta,\tau)dy.
\end{align}
Thus \eqref{eq:Gbound1} then becomes 
\begin{align*}
    G(0,\xi_2,\tau)= \frac{16}{\xi_2 t^{1-a}} \int_0^{+\infty} d\eta \int_0^{+\infty} \frac{d\xi_3}{(2|\tx|)^{1+a}} \int_{\frac{2|\tx|}{t}}^{+\infty} \frac{(\sin(\psi_1/2)  + \sin(\psi_2/2))}{y} \left(\frac{2 |\tx|}{ty}\right)^{a} \frac{1}{\sqrt{1-\frac{4|\tx|^2}{t^2y^2}}}  G(y,\eta,\tau)dy.
\end{align*}
Applying Fubini on the variables $\xi_3,y$ yields
\begin{align*}
    G(0,\xi_2,\tau)&=\frac{16}{\xi_2 t^{1-a}}  \int_0^{+\infty} d\eta \int_0^{+\infty} \frac{d\xi_3}{(2|\tx|)^{1+a}} \int_{\frac{2|\tx|}{t}}^{+\infty}  \frac{(\sin(\psi_1/2)  + \sin(\psi_2/2))}{y} \left(\frac{2 |\tx|}{ty}\right)^{a} \frac{1}{\sqrt{1-\frac{4|\tx|^2}{t^2y^2}}}  G(y,\eta,\tau)dy \notag \\
    & =  \frac{16}{\xi_2 t^{1-a}} \int_{0}^{+\infty} d \eta \int_{\frac{2|\xi_2|}{t}}^{+\infty} G(y,\eta,\tau) dy \int_0^{\sqrt{\left(\frac{yt}{2}\right)^2-\xi_2^2}} \frac{d\xi_3}{(2|\tx|)^{1+a}}\frac{(\sin(\psi_1/2)  + \sin(\psi_2/2))}{y} \left(\frac{2 |\tx|}{ty}\right)^{a} \frac{1}{\sqrt{1-\frac{4|\tx|^2}{t^2y^2}}}  
\end{align*} 
where we used that $\xi_3= \sqrt{\frac{y^2t^2}{4}-\xi_2^2}$ for $y=\frac{2|\tx|}{t}$.  We now set
\begin{align*}\label{def:Q}
    Q(\xi_2,y,t) &: = \frac{1}{y} \int_0^{\sqrt{\left(\frac{yt}{2}\right)^2-\xi_2^2}} \frac{d\xi_3}{(2|\tx|)^{1+a}}  \left(\frac{2 |\tx|}{ty}\right)^{a} \frac{1}{\sqrt{1-\frac{4|\tx|^2}{t^2 y^2}}}(\sin(\psi_1/2)  + \sin(\psi_2/2)).
\end{align*}
Performing the change of variables $\xi_3=|\xi_2|z$ we can rewrite $Q(\xi_2,y,t)$ above as  
\begin{align}
    Q(\xi_2,y,t) & = \frac{1}{y} \int_0^{\sqrt{\left(\frac{yt}{2\xi_2}\right)^2-1}} \frac{dz}{2^{1+a}(1+z^2)^{\frac{1+a}{2}}} (\sin (\psi_1/2)+\sin(\psi_2/2)) \left(\frac{2(1+z^2)^{\frac{1}{2}}}{ty}\right)^a \frac{1}{\sqrt{1-\frac{4|\xi_2|^2(1+z^2)}{t^2y^2}}}  \notag \\
    & = \frac{1}{2y(ty)^a} \int_0^{\sqrt{\left(\frac{yt}{2\xi_2}\right)^2-1}} \frac{dz}{(1+z^2)^{\frac{1}{2}}} (\sin (\psi_1/2)+\sin(\psi_2/2))  \frac{1}{\sqrt{1-\frac{4|\xi_2|^2(1+z^2)}{t^2y^2}}}.
\end{align}
Using elementary trigonometric formulas, namely $\sin(\psi/2)=\sqrt{\frac{1-\cos\psi}{2}}= \sqrt{\frac{1 \pm \sqrt{1-\sin^2 \psi}}{2}}
$ 
where we consider the plus sign for $\psi_1$ and the negative sign for $\psi_2$, we can rewrite the integrand in the right hand side of the equation above. More precisely,   
\begin{align*}
\sin(\psi_1/2) &+ \sin(\psi_2/2) = \sqrt{\frac{1 + \sqrt{1-\sin^2 \psi}}{2}}+\sqrt{\frac{1 -\sqrt{1-\sin^2 \psi}}{2}} \notag \\
    & = \frac{1}{\sqrt{2}} \left[ \sqrt{1+\sqrt{1-\frac{4|\xi_2|^2(1+z^2)}{t^2y^2}}} + \sqrt{1-\sqrt{1-\frac{4|\xi_2|^2(1+z^2)}{t^2y^2}}}  \right]
\end{align*}
where we used that $|\tx|^2=|\xi_2|^2(1+z^2)$. Hence, we obtain
\begin{align*}\label{eq:Q_H}
    Q(\xi_2,y,t) &=  \frac{1}{2y(ty)^a} \int_0^{\sqrt{\left(\frac{yt}{2\xi_2}\right)^2-1}} \frac{dz}{(1+z^2)^{\frac{1}{2}}}   \frac{1}{\sqrt{1-\frac{4|\xi_2|^2(1+z^2)}{t^2y^2}}} \times\notag \\
    & \frac{1}{\sqrt{2}} \left[ \sqrt{1+\sqrt{1-\frac{4|\xi_2|^2(1+z^2)}{t^2y^2}}} + \sqrt{1-\sqrt{1-\frac{4|\xi_2|^2(1+z^2)}{t^2y^2}}}  \right] \notag \\
    & =: \frac{1}{y(ty)^a} H\left(\frac{2\xi_2}{ty}\right)
\end{align*}
where $H(s):[0,1]\rightarrow \mathbb{R}_+$ is an even function and it has the form
\begin{equation}\label{eq:Hs}
    H(s)  = 2^{-\frac{3}{2}} \int_0^{\sqrt{\frac{1}{s^2}-1}} \frac{dz}{(1+z^2)^{\frac{1}{2}}} \frac{1}{\sqrt{1-s^2(1+z^2)}} 
     \times \left[\sqrt{1+\sqrt{1-s^2(1+z^2)}}+ \sqrt{1-\sqrt{1-s^2(1+z^2)}}\right]. 
\end{equation}
Changing variables again with $\zeta= s \sqrt{1+z^2}$ and  $dz= \frac{\zeta}{s^2\sqrt{\frac{\zeta^2}{s^2}}-1}d \zeta$ yields 
\begin{align*}\label{eq:HsZetaVariable}
    H(s) & = 2^{-\frac{3}{2}} \int_s^1 \frac{\zeta d\zeta}{s^2\sqrt{\frac{\zeta^2}{s^2}-1}} \frac{1}{\sqrt{1-\zeta^2}} \frac{s}{\zeta} \left[\sqrt{1+\sqrt{1-\zeta^2}}+\sqrt{1-\sqrt{1-\zeta^2}}\right] \notag \\
    & = 2^{-\frac{3}{2}} \int_s^1 \frac{d\zeta}{\sqrt{\zeta^2-s^2}} \frac{1}{\sqrt{1-\zeta^2}} \left[\sqrt{1+\sqrt{1-\zeta^2}}+\sqrt{1-\sqrt{1-\zeta^2}}\right].
\end{align*}
Further setting $\sqrt{1-\zeta^2}=x$, we obtain 
\begin{align*}
    H(s)=2^{-\frac{3}{2}} \int_0^{\sqrt{1-s^2}} \frac{dx}{\sqrt{1-x^2}\sqrt{1-s^2-x^2}} [\sqrt{1+x}+\sqrt{1-x}]
\end{align*}
and, finally, denoting $x= \sqrt{1-s^2}\sin \theta$, we arrive at 
\begin{align*}
    H(s)&=2^{-\frac{3}{2}} \int_0^{\pi/2} \frac{d \theta}{\sqrt{1-(1-s^2)\sin^2\theta}}[\sqrt{1+\sqrt{1-s^2}\sin \theta}+ \sqrt{1-\sqrt{1-s^2}\sin \theta}] \notag 
\end{align*}
whence \eqref{eq:defHs} follows, namely 
\begin{align}
    H(s) &=  2^{-\frac{3}{2}} \int_{-\pi/2}^{\pi/2} \frac{d \theta}{\sqrt{1-(1-s^2)\sin^2\theta}}\sqrt{1+\sqrt{1-s^2}\sin \theta} \notag \\
    &=  2^{-\frac{3}{2}} \int_{-\pi/2}^{\pi/2} \frac{d \theta}{\sqrt{1-\sqrt{1-s^2}\sin\theta}}.
\end{align}
\medskip
In the end we obtain
\begin{equation*}
    G(0,\xi_2,\tau)= \frac{16}{\xi_2 t} \int_{0}^{+\infty} d \eta \int_{\frac{2|\xi_2|}{t}}^{+\infty} dy\, G(y,\eta,\tau) \frac{1}{y^{1+a}} \, H\left(\frac{2\xi_2}{ty}\right) \, .
\end{equation*}

\end{proofof}

We now prove Lemma \ref{lem:asymHs} concerning the asymptotics of the function $H(s)$ defined as in \eqref{eq:defHs}. 

\begin{proofof}[Proof of Lemma \ref{lem:asymHs}]
We recall that
\begin{equation}\notag
    H(s) =  2^{-\frac{3}{2}} \int_{-\pi/2}^{\pi/2} \frac{d \theta}{\sqrt{1-\sqrt{1-s^2}\sin\theta}}, \quad s \in [0,1].
\end{equation}
Now, let $ 0<\delta< \pi /2$. In order to study the asymptotic behaviour of $H(s)$ for $s \rightarrow 0^+$ we consider separately the sets where $\theta\in [-\pi/2, \pi/2-\delta]$ and $\theta\in [\pi/2- \delta, \pi/2]$. More precisely, we write
\begin{align*}\label{eq:Hs_asymp}
    H(s) & = 2^{-\frac{3}{2}} \int_{-\pi/2}^{\pi/2-\delta} \frac{d \theta}{\sqrt{1-\sqrt{1-s^2}\sin\theta}}+ 2^{-\frac{3}{2}} \int_{\pi/2-\delta}^{\pi/2} \frac{d \theta}{\sqrt{1-\sqrt{1-s^2}\sin\theta}} = H^{(1)}(s)+H^{(2)}(s) 
\end{align*}
We notice that $H^{(1)}(s)$ is bounded since we have that
\begin{equation}\notag
    \frac{1}{\sqrt{1-\sqrt{1-s^2}\sin\theta}} \leq C_{\delta}
\end{equation}
for $\theta \in [-\pi/2,\pi/2-\delta]$. We now consider $H^{(2)}(s)$. If $\delta$ is chosen to be close enough to $\pi/2$, using the Taylor expansion of the functions $\cos \theta,\, \sin \theta$ at $\theta= \pi/2$ we obtain
    \begin{align*}
    H^{(2)}(s)& \sim 2^{-\frac{3}{2}}\int_{\pi/2-\delta}^{\pi/2} \frac{d \theta}{\sqrt{1-\sqrt{1-s^2}(1-(\pi/2-\theta)^2)}} =  2^{-\frac{3}{2}} \frac{1}{\sqrt[4]{1-s^2}}\int_0^{\sqrt{\frac{1-\sqrt{1-s^2}}{\sqrt{1-s^2}}}\delta} \frac{dX}{\sqrt{1+X^2}} \notag\\
    & = \frac{2^{-\frac{3}{2}}}{\sqrt[4]{1-s^2}} \left[\tanh^{-1}\left(\frac{x}{\sqrt{1+x^2}}\right)\right]_0^{{\sqrt{\frac{1-\sqrt{1-s^2}}{\sqrt{1-s^2}}}\delta}}
    \sim \frac{2^{-\frac{3}{2}}}{\sqrt[4]{1-s^2}} \left(\frac{1}{2}\log \left(\frac{1-\sqrt{1-s^2}}{\sqrt{1-s^2}} \right)+ \frac{1}{2} \log \sqrt{1-s^2}+ \log \delta \right)
    \end{align*}
 where in the first line we have used the change of variables $\sqrt{\frac{1-\sqrt{1-s^2}}{\sqrt{1-s^2}}}(\pi/2-\theta)=X$ and where we have used the following formula for $\tanh^{-1}$
\begin{equation*}
    \tanh^{-1}(z)=\frac{1}{2}(\log(1+z)-\log(1-z)), \quad z \in \mathbb{C}\backslash \{\pm1\}.
\end{equation*}
Hence we obtain the asymptotics 
\begin{equation*}
    H(s) \sim  -\log (s) \quad \text{as $s \rightarrow0^+$.}
\end{equation*}
\end{proofof}\\
We now present the proof of Lemma \ref{lem:Ka}. 
\medskip

\begin{proofof}[Proof of Lemma \ref{lem:Ka}]
We recall the explicit expression of $K(a)$ defined in \eqref{def:Ka} 
\begin{equation*} 
  K(a) =  \Ii d \zeta  \int_{2\zeta}^{+\infty}  \frac{1}{\eta^{a+2}} \,\exp\left(-\frac{4 \pi \,\eta^{-a}}{1-a}\right)   H \left(\frac{2\zeta}{ \eta}\right) d\eta.
\end{equation*}
By performing a  convenient change of variable, namely  $\eta= 2\zeta z$, 
$K(a)$ reads
\begin{align}\label{eq:KaIntermediate}
K(a) & = \Ii d \zeta  \int_{1}^{+\infty}  \frac{1}{(2\zeta)^{a+1}z^{a+2}} \exp\left(-\frac{4 \pi \, (2\zeta)^{-a}z^{-a}}{1-a}\right)   H \left(\frac{1}{z}\right)  dz \notag \\
    &= \Ii d\zeta\, \frac{1-a}{8\pi a} \frac{d}{d\zeta}  \left( \exp\left(-\frac{4 \pi \, (2\zeta)^{-a}z^{-a}}{1-a}\right)\right)  \int_1^{+\infty}  \frac{1}{z^2} H \left(\frac{1}{z}\right) dz \notag \\
    & = \frac{1-a}{8\pi a}\int_1^{+\infty}\frac{1}{z^2}H\left(\frac{1}{z}\right)dz = \frac{1-a}{8 \pi a} \int_0^1 H(s)ds. \notag
\end{align}
where in the second identity we used that $\frac{d}{d\zeta}\left(\exp\left(-\frac{4 \pi (2\zeta)^{-a}z^{-a}}{1-a}\right)\right)=\frac{8\pi a}{(1-a)}z^{-a}(2\zeta)^{-(a+1)}\exp\left(-\frac{4 \pi (2\zeta)^{-a}z^{-a}}{1-a}\right)$ and in the last one we performed the change of variable $s=\frac 1 z$, with $ds=-\frac 1 {z^2}dz$. 

We now compute $\int_0^1H(s)ds$ where $H(s)$ has been defined in \eqref{eq:defHs}. Using Taylor we obtain
\begin{align}
    \int_0^1H(s)ds &= 2^{-\frac{3}{2}} \int_0^1 ds \int_{-\frac{\pi}{2}}^{\frac{\pi}{2}} \frac{d \theta}{\sqrt{1-\sqrt{1-s^2}\sin \theta}} \notag \\
    & = 2^{-\frac{3}{2}} \int_0^1 ds \; \Su \binom{-1/2}{n} (-1)^n (1-s^2)^{\frac{n}{2}} \int_{-\frac{\pi}{2}}^{\frac{\pi}{2}} \sin^n \theta d \theta \notag \\
    & = 2^{-\frac{3}{2}} 2 \int_0^1 ds \; \Su \binom{-1/2}{2n} (-1)^{2n}(1-s^2)^{n} W_{2n} 
\end{align}
where $W_{2n}= \int_0^{\frac{\pi}{2}} \sin^{2n} \theta d \theta$ is the $n-$th Wallis integral, see formula 6.1.49 in \cite{AS}. From the same formula we have that 
\begin{equation}
    W_{2n}= \frac{\sqrt{\pi}\,\Gamma(n+1/2)}{2 \Gamma(n+1)}
\end{equation}
where $\Gamma(z)$ is the \emph{Gamma} function (see 6.1.1 in \cite{AS}). Therefore 
\begin{align}\label{eq:taylorHs}
   2^{-\frac{3}{2}}  2 \int_0^1 ds \; \Su \binom{-1/2}{2n} (-1)^{2n}(1-s^2)^{n} W_{2n} &= 2^{-\frac{3}{2}} 2 \int_0^1 da(s) \; \Su \binom{-1/2}{2n}(1-s^2)^n \frac{\sqrt{\pi}\,\Gamma(n+1/2)}{2 \Gamma(n+1)} \notag \\
    & = 2^{-\frac{3}{2}} 2  \Su \binom{-1/2}{2n} \frac{\sqrt{\pi}\,\Gamma(n+1/2)}{2 \Gamma(n+1)} \int_0^1 (1-s^2)^n ds \notag \\
    & = 2^{-\frac{3}{2}} 2  \Su \binom{-1/2}{2n} \frac{\sqrt{\pi}\,\Gamma(n+1/2)}{2 \Gamma(n+1)} \frac{\sqrt{\pi}}{2} \frac{\Gamma(n+1)}{\Gamma(n+3/2)}  \notag \\
    & = 2^{-\frac{3}{2}} \frac{\pi}{2} \Su \binom{-1/2}{2n} \frac{\Gamma(n+1/2)}{\Gamma(n +3/2)} \notag \\
    & = 2^{-\frac{3}{2}} \frac{\pi}{2}  \Su \binom{-1/2}{2n} \frac{2}{2n+1}
\end{align}
where in the last step we have used the identity
\begin{equation*}
    \Gamma(n+1/2)=2^{1-2n} \sqrt{\pi} \frac{\Gamma(2n)}{\Gamma(n)}
\end{equation*}
and we have used the fact that 
\begin{equation}
    \int_0^1(1-s^2)^nds=\frac{\sqrt{\pi}}{2} \frac{\Gamma(n+1)}{\Gamma(n+3/2)}.
\end{equation}
The former equality follows by noticing that $\int_0^1 (1-s^2)^nds= 2^{2n} B(n,n)$, where $B(z_1,z_2)$ is the \emph{Beta function}, and form \emph{Legendre's duplication formula} (see  6.2.1 and 6.1.18 respectively in \cite{AS}).
Considering the quantity
\begin{equation}\label{eq:auxint}
    \int_0^1 \frac{dx}{\sqrt{1+x}}+\int_0^1 \frac{dx}{\sqrt{1-x}} = 2\sqrt{2}
\end{equation}
and, using the Taylor expansion in the variable $x$ in both integrals in the l.h.s. of \eqref{eq:auxint}, we obtain that the last series in \eqref{eq:taylorHs} is convergent and we have the sum
\begin{equation*}
    \Su \binom{-1/2}{2n} \frac{1}{2n+1}= \Su \binom{4n}{2n} \frac{1}{2^{4n}(2n+1)}= \sqrt{2}.
\end{equation*}
This yields
\begin{align}\label{eq:integralHs}
    \int_0^1 H(s)ds= 2^{-\frac{3}{2}} \pi  \Su \binom{-1/2}{2n} \frac{1}{2n+1}= 2^{-\frac{3}{2}} \pi \sqrt{2} = \frac{\pi}{2}.
\end{align}
Finally we arrive at
\begin{equation}
    K(a)=\frac{(1-a)}{16a}.
\end{equation}

\bigskip 

\end{proofof}
\medskip

We now prove Lemma \ref{lem:R}. 

\begin{proofof}[Proof of Lemma \ref{lem:R}]
The proof of this lemma follows the lines of the proof of Lemma \ref{lem:bdvalue}. In fact, in \eqref{eq:BVPF} set $\psi=2\theta$ so that 
   \begin{align} 
    F(0,\xi_2,\xi_3,\tau)=\frac{16}{\xi_2(2|\tx|)^{1+a}t^{1-a}}\int_0^{\pi}d\psi \frac{\sin(\frac{\psi}{2})}{|\sin(\psi)|^{1-a}}\int_{\Rl} G\left(\frac{2|\tilde \xi|}{t \sin(\psi)},\eta, \tau\right)d\eta \, .\notag
\end{align} 
    Divide then the integral above in the cases $\psi \in [0,\pi/2], \psi \in [\pi/2,\pi]$, namely 
\begin{align}
    \int_0^{\pi}d\psi \frac{\sin(\frac{\psi}{2})}{|\sin(\psi)|^{1-a}}= \frac{1}{2}\int_0^{\pi/2}d \psi_1 \frac{\sin(\psi_1/2)}{|\sin(\psi_1)|^{1-a}} + \frac{1}{2}\int_{\pi/2}^{\pi}d \psi_2 \frac{\sin(\psi_2/2)}{|\sin(\psi_2)|^{1-a}}\, 
\end{align}
with $\psi_1\in [0,\pi/2]$ and $\psi_2\in [\pi/2,\pi]$. Changing also variables as $y = \frac{2|\tx|}{t \sin \psi}$, with $dy=- \frac{2|\tx|}{t \sin^2 \psi}\cos\psi \, d\psi$ and from \eqref{eq:splitpsi}  we obtain
\begin{align}
     &\int_0^{\pi}d\psi \frac{\sin(\frac{\psi}{2})}{|\sin(\psi)|^{1-a}}\int_{\Rl} G\left(\frac{2|\tilde \xi|}{t \sin(\psi)},\eta, \tau\right)d\eta \notag \\
     & = \int_{\frac{2|\tx|}{t}}^{+\infty}  \frac{\sin(\psi_1/2)}{y} \left(\frac{2|\tx|}{ty}\right)^{a} \frac{G(t,\eta,\tau)}{\sqrt{1-\frac{4|\tx|^2}{t^2y^2}}}dy  + \int_{\frac{2|\tx|}{t}}^{+\infty}  \frac{\sin(\psi_2/2)}{y}\left(\frac{2|\tx|}{ty}\right)^{a} \frac{G(y,\eta,\tau)}{\sqrt{1-\frac{4|\tx|^2}{t^2y^2}}}dy. \notag
\end{align}
Arguing as in Lemma \ref{lem:bdvalue} and using elementary trigonometric formulas, namely $\sin(\psi/2)=\sqrt{\frac{1-\cos\psi}{2}}= \sqrt{\frac{1 \pm \sqrt{1-\sin^2 \psi}}{2}}$ we get that   
\begin{align}
    \sin(\psi_1/2) &+ \sin(\psi_2/2) = \sqrt{\frac{1 + \sqrt{1-\sin^2 \psi}}{2}}+\sqrt{\frac{1 -\sqrt{1-\sin^2 \psi}}{2}} \notag \\
    & = \frac{1}{\sqrt{2}}\left(\sqrt{1 + \sqrt{1-\frac{4|\tx|^2}{t^2y^2}}}+\sqrt{1 -\sqrt{1-\frac{4|\tx|^2}{t^2y^2}}}\right) \notag
\end{align}
In the end we have
\begin{align}
   F(0,\xi_2,\xi_3,\tau)= & \frac{8}{\xi_2(2|\tx|)^{1+a}t^{1-a}} \int_0^{\pi}d\theta \frac{\sin(\theta)}{|\sin(2\theta)|^{1-a}}\int_{\Rl} G\left(\frac{2|\tilde \xi|}{t \sin(2\theta)},\eta, \tau\right)d\eta = \notag \\
   &=\frac{8}{\xi_2(2|\tx|)^{1+a}t^{1-a}} \int_0^{+\infty} d\eta \int_{\frac{2|\tx|}{t}}^{+\infty} \frac{dy}{y} G(y,\eta,\tau)  \notag \\
    & \times \frac{1}{\sqrt{2}}\left(\sqrt{1 + \sqrt{1-\frac{4|\tx|^2}{t^2y^2}}}+\sqrt{1 -\sqrt{1-\frac{4|\tx|^2}{t^2y^2}}}\right)\left(\frac{2|\tx|}{ty}\right)^{a} \frac{1}{\sqrt{1-\frac{4|\tx|^2}{t^2y^2}}}\notag \\
    & = \frac{8}{\xi_2(2|\tx|)^{1+a}t^{1-a}} \Ii d \eta \int_{\frac{2|\tx|}{t}}^{+\infty} G(y,\eta,\tau)R\left(\frac{2|\xi|}{ty}\right)\frac{dy}{y}.
\end{align}
\end{proofof}

\subsection{Why we cannot expect the mass of the solution of \eqref{eq:selfSimStrong2} to be contained within a single length scale?} \label{sec:noSelf}

It may be tempting to expect self-similar long-time behaviour for $F$. Indeed, the scaling obtained in Section \ref{ssec:scalinA>1} suggests that all terms in \eqref{eq:selfSimStrong2} balance. To investigate the self-similar behaviour of solutions to \eqref{eq:selfSimStrong2}, in order to study the behaviour of the solution near the origin, i.e. $|\xi|\to 0^+$, one possibility would be to assume that the function
\begin{equation}\label{eq:ansatz}
\lambda(\xi_1,\tau):=\int_{\Rp} F(\xi_1,\tx,\tau)d\tx
\end{equation}
behaves asymptotically like
\begin{equation}\label{eq:ansatz2}
\lambda(\xi_1,\tau) \sim \frac{A}{|\xi_1|^{1-a}} \; \text{as}\quad |\xi_1|\rightarrow 0^+, \text{with $A>0.$} 
\end{equation}
This is a natural assumption to make in order to cancel out the term $\frac{1}{t^{1-a}}$ appearing in the gain term \eqref{eq:selfSimStrong2_operators}. Unfortunately, this ansatz turns out to be incorrect. As we will see later in this section the resulting solution obtained with this ansatz is not consistent since  \eqref{eq:ansatz2} fails in the end due to a mismatch of the constant close to the origin, namely as we will show below, in the end we obtain 
\begin{equation*} \lambda(\xi_1) \sim \frac{(1-a)A}{|\xi_1|^{1-a}} \quad \text{as $|\xi_1| \rightarrow 0^+$.}
\end{equation*}
To prove this, plugging \eqref{eq:ansatz} into \eqref{eq:selfSimStrong2}, we obtain
\begin{gather*}
\partial_{\tau}F-\frac{1}{a}\partial_{\xi_1}(\xi_1F)-\left(\frac{1}{a}-1\right)\partial_{\xi_2}(\xi_2F)-\left(\frac{1}{a}-1\right)\partial_{\xi_3}(\xi_3F)+
    \xi_2 \partial_{\xi_1}F \notag \\
    =\frac{8 \delta(\xi_1)}{(2|\tilde \xi|)^{1+a}t^{1-a}}\int_0^{\pi}d\theta \frac{\sin(\theta)}{|\sin(2\theta)|^{1-a}}\lambda\left(\frac{2|\tilde \xi|}{t \sin(2\theta)},\tau\right)-4\pi |\xi_1|^{-a}F
\end{gather*}
Using the asymptotics \eqref{eq:ansatz2}, we get
$$
\lambda\left(\frac{2|\tilde \xi|}{t \sin(2\theta)},\tau\right)\sim \frac{A (t |\sin(2\theta)|)^{1-a}}{{(2|\tilde\xi|)^{1-a}}} 
$$ 
which, using $\int_0^{\pi}\sin \theta d\theta=2$, yields 
\begin{align*}
  \partial_{\tau}F-\frac{1}{a}\partial_{\xi_1}(\xi_1F)-\left(\frac{1}{a}-1\right)\partial_{\xi_2}(\xi_2F)-\left(\frac{1}{a}-1\right)\partial_{\xi_3}(\xi_3F)+
    \xi_2 \partial_{\xi_1}F =\frac{4 A}{|\tilde \xi|^2}\delta(\xi_1)-4\pi|\xi_1|^{-a}F \notag 
\end{align*}
We now look for stationary solutions of the equation above, namely 
\begin{equation*}
    -\frac{1}{a}\partial_{\xi_1}(\xi_1F)-\left(\frac{1}{a}-1\right)\partial_{\xi_2}(\xi_2F)-\left(\frac{1}{a}-1\right)\partial_{\xi_3}(\xi_3F)+
    \xi_2 \partial_{\xi_1}F 
    =\frac{4 A}{|\tilde \xi|^2}\delta(\xi_1)-4\pi|\xi_1|^{-a}F \ . \end{equation*}
The characteristics are given by \eqref{eq:characteristicsF}, hence we are then led to solve the following boundary value problem
\begin{equation}\label{eq:charsyst2}
    \begin{cases}
        \frac{d F}{ds}=\left(\frac{3}{a}-2\right)F-4\pi |\xi_1|^{-a}F, & \\
        F(0,\tilde \xi)= \frac{4\, A}{|\tilde \xi|^2\xi_2}.    \end{cases}
\end{equation}
We now check whether the asymptotic behaviour given by \eqref{eq:ansatz2} for $\lambda(\xi_1,\tau)$ as $|\xi_1|\to 0^+$. To this end it is sufficient to consider the system of characteristics \eqref{eq:characteristicsF} and $F$ satisfying 
\begin{equation}\label{eq:charsyst2bis}
    \frac{d F}{ds}=-4\pi |\xi_1|^{-a}F, \quad 
    F(0,\tilde \xi)=\frac{{4} \, A}{|\tilde \xi|^2\xi_2 }
\end{equation}
in place of
\begin{equation*}
    \frac{dF}{ds}=\left(\frac{3}{a}-2\right)F-4\pi|\xi_1|^{-a}F.
\end{equation*} 
This approximation is expected to be valid for $|\tilde\xi_0|$ small of order one. The system \eqref{eq:characteristicsF} is solved for
\begin{equation}\label{eq:FboundaryApprox}
    F(s,\xi_{2,0},\xi_{3,0})= \frac{4A}{|\tilde \xi_0|^2\xi_{2,0}} \exp \left(-4 \pi \int_0^s |\xi_1(\sigma)|^{-a}d \sigma\right).
\end{equation}
From \eqref{eq:FboundaryApprox} and \eqref{eq:xi1F}-\eqref{eq:xi3F}, it then follows that
\begin{align}\label{eq:Fxis}
F(\xi,s)&=4A \, \frac{ \exp\left(-\frac{4 \pi \xi_1^{1-a}}{(1-a)\xi_2}\right) }{|\tilde \xi|^2\xi_2 \exp\left(3\left(\frac{1}{a}-1\right)s\right)}
\end{align}
where we have used the fact that from \eqref{eq:lossCharacteristics} we have that $\int_0^s |\xi_1(\sigma)|^{-a} d \sigma= \frac{\xi_{2,0}^{-a}}{1-a}(e^s-1)^{1-a}= \frac{\xi_1^{1-a}}{(1-a)\xi_2}$. 
We observe that the denominator of \eqref{eq:Fxis}, using $e^{-s}-1 = \frac{\xi_1}{\xi_2}  $, can be rewritten as
\begin{equation*}
    |\tilde \xi|^2\xi_2 \exp\left(3\left(\frac{1}{a}-1\right)s\right)={|\tilde \xi|^2\xi_2}{\left(1-\frac{\xi_1}{\xi_2}\right)^{-3\left(\frac{1}{a}-1\right)}}=\frac{|\tilde \xi|^2\xi_2^{1+3 \left(\frac{1}{a}-1\right)}}{(\xi_2-\xi_1)^{3\left(\frac{1}{a}-1\right)}}
\end{equation*}
and then
\begin{equation}\label{eq:denFxis}
    \frac{1}{|\tilde \xi|^2\xi_2 \exp \left(3\left(\frac{1}{a}-1\right)s\right)}=\frac{(\xi_2-\xi_1)^{3\left(\frac{1}{a}-1\right)}}{|\tilde \xi|^2\xi_2^{1+3 \left(\frac{1}{a}-1\right)}}.
\end{equation} 
Combining  \eqref{eq:Fxis} and \eqref{eq:denFxis} we then obtain 
\begin{align}\label{eq:Fxis2}
 F(\xi,s)&=
   4 A \ \frac{(\xi_2-\xi_1)^{3\left(\frac{1}{a}-1\right)} }
{ |\tilde \xi|^2\xi_2^{1+3 \left(\frac{1}{a}-1\right)}}\exp\left(-\frac{4 \pi\xi_1^{1-a}}{(1-a)\xi_2}\right)\ . 
\end{align}
Starting from \eqref{eq:Fxis2}, we are now interested in checking  whether we recover the expected asymptotic behaviour for $\lambda(\xi_1,\tau)$, that is
\begin{equation}\label{eq:lambda}
\lambda(\xi_1,\tau)=\int_{\xi_1}^1 \int_{-1}^1 F(\xi_1,\tilde \xi ,\tau) d\tilde \xi \sim \frac{A}{|\xi_1|^{1-a}} \quad \text{as} \ |\xi_1|\to 0^+. 
\end{equation}
We recall that $\tilde \xi=(\xi_2,\xi_3)$ and that without loss of generality, we are assuming $0\leq \xi_1 \leq \xi_2$. Plugging \eqref{eq:Fxis2} into \eqref{eq:lambda} we obtain
$$\lambda(\xi_1,\tau)
= 4 A \ \int_{\xi_1}^1 \int_{-1}^1 \frac{(\xi_2-\xi_1)^{3\left(\frac{1}{a}-1\right)} }
{ |\tilde \xi|^2\xi_2^{1+3 \left(\frac{1}{a}-1\right)}}\exp\left(-\frac{4 \pi\xi_1^{1-a}}{(1-a)\xi_2}\right)  d\tilde \xi. $$
Furthermore we have that
\begin{align*}   
\lambda(\xi_1,\tau) & =  
\int_{\xi_1}^1\int_{-1}^1\frac{(\xi_2-\xi_1)^{3\left(\frac{1}{a}-1\right)}}{|\xi_2^2+\xi_3^2| \xi_2^{1+3 \left(\frac{1}{a}-1\right)}}\exp\left(-\frac{4 \pi\xi_1^{1-a}}{(1-a)\xi_2}\right)d\xi_2d\xi_3 \notag \\ 
   &= \frac{4A}{\xi_1} \int_1^{+\infty}d\theta_2 \int_{-\infty}^{+\infty}\frac{(\theta_2-1)^{3\left(\frac{1}{a}-1\right)}}{|\theta_2^2+\theta_3^2|\theta_2^{1+3 \left(\frac{1}{a}-1\right)}}\exp\left(-\frac{4 \pi|\xi_1|^{-a}}{(1-a)\theta_2}\right)d\theta_2d\theta_3
\end{align*}
where we have performed the change of variables $\xi_2={\xi_1}\theta_2, \, \xi_3=\xi_1 \theta_3$. Further setting $\theta_3=\theta_2 \eta$ with $d\theta_3=\theta_2 d\eta$ we get
\begin{align*}
   \lambda(\xi_1,\tau) &=  
   \frac{4A}{\xi_1} \int_1^{+\infty}d\theta_2 \int_{-\infty}^{+\infty}\frac{(\theta_2-1)^{3\left(\frac{1}{a}-1\right)}}{\theta_2^2 \ |1+\eta^2|\, \theta_2^{3 \left(\frac{1}{a}-1\right)}}\exp\left(-\frac{4\pi |\xi_1|^{-a}}{(1-a)\theta_2}\right)d\eta \notag \\
   & = 
   \frac{4A}{\xi_1} \int_1^{+\infty}d\theta_2 \frac{(\theta_2-1)^{3\left(\frac{1}{a}-1\right)}}{\theta_2^2  \theta_2^{3 \left(\frac{1}{a}-1\right)}}\exp\left(-\frac{4\pi |\xi_1|^{-a}}{(1-a)\theta_2}\right)d\eta.
\end{align*}
Recalling that $\theta_2= \frac{1}{(1-e^{-s})}$, or equivalently $\theta_2-1=\frac{e^{-s}}{(1-e^{-s})}$, with  $d\theta_2=-\frac{e^{-s}}{(1-e^{-s})^2}ds$, we then arrive at 
\begin{equation*}
\lambda(\xi_1,\tau)=\frac{4\pi A}{\xi_1}
    \int_0^{+\infty}e^{-s}e^{-3\left(\frac{1}{a}-1\right)s}\exp\left(-\frac{4\pi|\xi_1|^{-a}(1-e^{-s})}{1-a}\right)ds.
\end{equation*}
From \eqref{eq:xi1F} we have that the asymptotic for $\lambda(\xi_1,\tau)$ for $|\xi_1|\rightarrow 0^+$ is obtained by taking $s \rightarrow 0^+$ in the previous integral. Therefore we obtain
\begin{align}
\lambda(\xi_1,\tau) &\sim \frac{4\pi A}{\xi_1} \int_0^{\delta}\exp\left(-4\pi|\xi_1|^{-a} \frac{s}{1-a}\right)ds \sim \frac{4\pi A}{\xi_1} \int_0^{+ \infty}\exp\left(-4\pi|\xi_1|^{-a} \frac{s}{1-a}\right)ds \notag \\
& = \frac{4 \pi A}{|\xi_1|} \frac{1-a}{4 \pi|\xi_1|^{-a}}= \frac{(1-a)A}{|\xi_1|^{1-a}}
\end{align}
where $\delta>0$ is small and independent of $\xi_1$. Therefore, we do not recover the expected asymptotic behaviour of $\lambda(\xi_1)$ due to the presence of the  additional factor $(1-a)$ in front of the constant $A$. Thus, the ansatz \eqref{eq:ansatz}-\eqref{eq:ansatz2} is not self-consistent. This is related to the fact that, as discussed in Section \ref{sec:strategy}, we cannot expect a simple power law behaviour.

\section{On the probabilistic interpretation of \eqref{eq:homAdj}}\label{sec:probabilities}

In this section we present the proofs of Lemmas \ref{lem:pFlight}, \ref{lem:pJump} and Lemma \ref{lem:transitionp} stated  in Section \ref{sec:heuristics}. 
\setcounter{subsection}{1}

\subsubsection{Flight Probabilities}

 Here we prove Lemma \ref{lem:pFlight} which provides an explicit formula for the flight probability, namely the probability $p(\xi^m| \, v^m)$ to update the velocity of the tagged particle from $v^m$, the velocity up to  time $t^{m}$, to  $\xi^m$. 
 For simplicity we will drop the indexes in the notation by denoting $v=v^m$, $\xi=\xi^m$, $t^{m}=\tau$. \smallskip
 
\begin{proofof}[Proof of Lemma \ref{lem:pFlight}]
     In order to compute the probability density $p(\xi|\, v)$  we start from \eqref{eq:time}, which we recall reads
\begin{equation*}
        q(\tau| \, v)=|T_{\tau}v|^{-a}\exp\left(-\int_0^{\tau}|T_sv|^{-a}d s \right).
\end{equation*}
    To obtain $p(\xi|v)$, it is sufficient to compute the probability  to arrive at $\xi=(\xi_1,\xi_2,\xi_3)$ at time $\tau$ under the action of the map $T_{\tau}$, that is 
    \begin{equation}
        p(\xi|\, v) := \int_0^{+ \infty} \delta(\xi-T_\tau v)q(\tau|v) d \tau
    \end{equation}
    where
	\begin{equation}
	    \delta(\xi-T_{\tau}v)=\delta(\xi_1-(v_1+\tau v_2))\delta(\xi_2-v_2)\delta(\xi_3-v_3).
	\end{equation}
	If the free flight due to the shear stops at time $\tau$ then $\delta(\xi_1-(v_1+\tau v_2)$ is non-zero if and only if $\xi_1=v_1+\tau v_2$ from which we obtain
	\begin{equation}
	    \tau=\frac{\xi_1-v_1}{|v_2|} \geq 0.
	\end{equation}
	Thus $p(\xi| \, v)$ becomes
\begin{align}\label{flightProb1}
		p(\xi|\, v) & = \int_0^{+\infty}p(\tau|v)\delta(\xi-v)d \tau \notag \\
		&=\left( \frac{1}{|v_2|}\; p\left(\frac{\xi_1-v_1}{v_2}\Big| v \right) \chi_{\{v_2\geq 0\}}+\frac{1}{|v_2|}\; p\left(\frac{v_1-\xi_1}{v_2}\Big| v \right) \chi_{\{v_2< 0\}}\right) \delta(\xi_2-v_2)\delta(\xi_3-v_3).
	\end{align}
	Introducing the function $G(\xi_1|v)$ we can then rewrite  \eqref{flightProb1} as $p(\xi|\,v)=G(\xi_1| \, v)\delta(\xi_2-v_2)\delta(\xi_3-v_3)$, where $G(\xi_1|\,v)$ is defined as in \eqref{eq:defGprob}, namely
    \begin{align*} 
		G(\xi_1|\, v) & := \frac{1}{|v_2|}\left(\xi_1^2+v_2^2+v_3^2 \right)^{-\frac{a}{2}}	\exp\left(- \frac{1}{|v_2|}\int_{v_1}^{\xi_1}\ \left(\eta^2+v_2^2+v_3^2 \right)^{-\frac{a}{2}}d \eta \right)\chi_{\{v_2\geq 0\}}  \notag \\
		& + \frac{1}{|v_2|}\left(\xi_1^2+v_2^2+v_3^2 \right)^{-\frac{a}{2}} \exp\left(-\frac{1}{|v_2|}\int_{\xi_1}^{v_1}(\eta^2+v_2^2+v_3^2 )^{-\frac{a}{2}}d \eta  \right)\chi_{\{v_2<0\}}.
	\end{align*} 
    We notice that if $v_2=0$, we might expect a singularity, but if $v_2=0$ then $T_{\tau}=\mathrm{Id}$ and therefore we simply have $p(\xi|v)=\delta(\xi-v).$ 
    \end{proofof}
    
\subsubsection{Jump Probabilities} 

We compute the probability for the particle to collide with a background particle and to update its velocity from $\xi^m$ to $v^{m+1}$. \smallskip

\begin{proofof}[Proof of Lemma \ref{lem:pJump}] In order to simplify the notation, in what follows, we will set $v^{m+1}=V, \xi^m=\xi, v^m= y$. With this new notation, we have that $V=\mathrm{P}^{\perp}_\omega(\xi) $. Therefore the probability of jumping to $V$ given $\xi$ and given the probability $p(\xi|y)$ for the particle to have flown to $\xi$  is
	\begin{equation}\label{pVk}
		p(V=\bar v| \,\xi)=\frac{1}{4\pi}\int_{S^2}\int_{\R}p(\xi| \, y)\delta(V-\bar v)d \omega d \xi, \quad \overline{v} \in \R.
	\end{equation} 
Here $\bar v$ just plays the role of an auxiliary variable that will be dropped in the end. Notice that we can write any $w \in \R$ using its orthogonal decomposition  with respect to $\omega \in S^2$, i.e.
    \begin{equation}
        v=(w\cdot \omega)\omega+[w-(w\cdot \omega)\omega]=\mathrm{P}_{\omega}w+\mathrm{P}_{\omega}^{\perp}w.
    \end{equation} 
Since $V \perp \omega$ and from the fact that $\delta(a+b)=\delta(a)\delta(b)$ if $a,b \in \R$ are perpendicular to each other,  we have that
        \begin{gather}\label{eq:deltaOrtg}
            \delta(V-\bar v)=\delta(\mathrm{P}_{\omega}^{\perp}(V)-\mathrm{P}_{\omega}^{\perp}(\bar v)-(\bar v\cdot \omega)\omega)=
            \delta(\mathrm{P}_{\omega}^{\perp}(\xi-\bar v))\delta(\bar v \cdot\omega).
        \end{gather}
Plugging into \eqref{pVk} and denoting with $\xi^{\parallel},\xi^{\perp}$ respectively the parallel and orthogonal component of $\xi$ with respect to $\omega$ gives
	\begin{align}
		p(V=\bar v|\xi) & =\frac{1}{4\pi}\int_{S^2}\int_{\R}p(\xi| \, y)\delta(V-\bar v)d \omega d\xi  \notag  \\
		& =\frac{1}{4\pi}\int_{S^2}d\omega\int_{\Pi_{\perp}(\omega)}\int_{\Pi_{\parallel}(\omega)}\; p(\xi^{\perp},\xi^{\parallel} | \, y)\delta(\mathrm{P}_{\omega}^{\perp}(\xi-\bar v))\delta(\bar v \cdot\omega)d \xi^{\perp} d \xi^{\parallel} \label{eq:pParallelOrthogonal}
	\end{align}
where $\Pi_{\perp}(\omega),\Pi_{\parallel}(\omega)$ are, respectively, the plane orthogonal to $\omega$ containing $\omega$ and the line passing through $\omega$ parallel to $\omega$. Integrating the Dirac delta with respect to the orthogonal component yields
    \begin{equation}
		p(V=\bar v | \, \xi)=\frac{1}{4\pi}\int_{S^2}\delta(\bar v \cdot \omega) d \omega \int_{\Pi_{\parallel}(\omega)}p(\overline{v}^{\perp}(\omega),\xi^{\parallel}|\, y)d\xi^{\parallel}
	\end{equation}
where $\bar v^{\perp} = \bar v^{\perp}(\omega)$ since its orthogonal projection depends on $\omega$. Without loss of generality, it is possible to assume $\bar v=|\bar v|(1,0,0)$ from which we have
\begin{equation}\label{eq:deltaJumpOmega1}
    \delta(\bar v\cdot\omega)=\frac{1}{|\bar v|}\delta(\omega_1).
\end{equation}
Furthermore if $\bar v \cdot \omega=0$ then $\bar v^{\perp}(\omega) = \bar v $. Therefore using \eqref{eq:deltaJumpOmega1} and passing into spherical coordinates centred in $\bar v$ gives
	\begin{align}
		p(V=\bar v|\, \xi) &=  \frac{1}{2 \pi|\bar v|}\int_0^{2\pi} d \vp \int_0^{\pi} \delta(\cos \theta) d\theta \int_0^{+\infty}p(\bar v+re_r(\vp,r)|\, y)d r \notag \\
        & = \frac{1}{2 \pi|\bar v|}\int_0^{2\pi} d \vp \int_0^{+\infty}p(\bar v+re_r(\vp,r)|\, y)d r
	\end{align}
where $e_r(\phi,r))$ is a unit vector orthogonal to $\bar v$ and where we have assumed $p$ to be even with respect to $r$. Furthermore
        \begin{align}
           p(V=\bar v|\, \xi) & = \frac{1}{2 \pi |\bar v|}\int_0^{2\pi} d \phi \int_0^{+\infty}p(\bar v+re(\phi,r)| \, y)d r \notag \\
           & = \frac{1}{2 \pi |\bar v|}\int_0^{2\pi}\frac{1}{|\eta|} d \phi \int_0^{+\infty}|\eta| p(\bar v+re(\phi,r)| \, y)d r \
			=\frac{1}{2\pi|\bar v|}\int_{\Pi^{\perp}(\bar v)}\frac{p(\overline{v}+\eta| \, y)}{|\eta|}d S(\eta)
	\end{align}
    where we have set $re(\phi,r)=\eta \perp \bar v$ and we have used $\int_0^{2\pi}d\vp|\eta|dr= dS(\eta)$.
Hence, in the indexed notation, we have obtained
    \begin{equation}\label{eq:jumpInter}
        p(v^{m+1}|\, \xi^m)=\frac{1}{2\pi| v^{m+1}|}\int_{\Pi^{\perp}(v^{m+1})}\frac{p(v^{m+1}+\eta|\, v^m)}{|\eta|}d S(\eta).
    \end{equation}
Notice that we removed the dependence on $\bar v$ since from \eqref{eq:deltaOrtg}-\eqref{eq:pParallelOrthogonal} we have that $\bar v \perp \omega$ and that $\mathrm{P}_{\omega}^{\perp}(V)=\mathrm{P}_{\omega}^{\perp}(\bar v)$ which together imply $\bar v = V$. In order to further simplify \eqref{eq:jumpInter1} it is convenient to determine the support of the operator $\mathrm{P}_{\omega}^{\perp}$. For simplicity we ignore again the indexes in the following computations. We recall that $V=\xi-(\xi \cdot \omega)\omega$ and we notice that $\mathrm{P}_{R\omega}^{\perp}(R\xi) = R \mathrm{P}_{\omega}^{\perp}(\xi)$ for every rotation $R$. Hence without loss of generality it is possible to assume $\xi=r(1,0,0)$, which yields $\mathrm{P}_{\omega}^{\perp}(\xi)$ as
    \begin{equation}\label{support}
	V(\omega)=r\left( \begin{array}{c}
		1-\omega_1^2\\
		-\omega_1\omega_2		\\
			-\omega_1\omega_3	
	\end{array}\right).
	\end{equation}
    Since $\omega \in S^2$ passing into polar coordinates with the North pole aligned to $e_1$ gives
	\begin{equation}\label{supportPolar}
	V(\theta, \varphi)=r\left( \begin{array}{c}
		\sin^2 \theta\\
		-\sin\theta \cos \theta  \sin \varphi		\\
			-\sin\theta\cos \theta \sin \varphi	
	\end{array}\right), \quad \theta \in [0,\pi], \vp \in [0,2\pi].
	\end{equation}
	The surface described by $V(\theta,\varphi)$ will be denoted by $\Gamma(\xi)$. Notice that for $\xi \neq 0$ one has $0,\xi \in \Gamma(\xi)$ and moreover if $A$ is a rotation, then $\Gamma(A\xi)=A(\Gamma(\xi))$. The jump probability is then obtained as follows
	\begin{equation}
		F(v|\, \xi)=\frac{1}{4\pi}\int_{S^2}\delta(V(\xi,\omega)-v)d \omega.
	\end{equation}
	Notice that from the definition the probability measure $F$ is supported on $\Gamma(\xi)$, we can then write $F(v|\xi)=M(v,\xi)\delta_{\Gamma(\xi)}(v)$, with $M(v,\xi)$ a probability density function to be determined. By definition we set
    \begin{align}
        \int_{\R} \delta_{\Gamma(\xi)}(v)\vp(v) dv &:= \int_{\Gamma(\xi)} \vp(v) d S(v) \quad \text{for every $\vp \in C_c^{\infty}(\R)$;} \notag \\
        \int_{\Gamma(\xi)}M(v,\xi)dv & = 1.
    \end{align}
    The identity $F(v|\xi)=M(v,\xi)\delta_{\Gamma(\xi)}(v)$ has to be understood in the sense of operators acting on $\mathcal{P}_+(\R)$ i.e. the space of probability measures, namely
	\begin{equation}\label{eq:jumpInter1}
		\int_{\R}F(v|\, \xi)p(\xi)d \xi=\int_{\R}M(v,\xi)\delta_{\Gamma(\xi)}(v)p(\xi) d \xi, \quad p \in \mathcal{P}_+(\R). 
	\end{equation}
    The two probability densities we have obtained so far, namely \eqref{eq:jumpInter}-\eqref{eq:jumpInter1} must coincide, where we remark that equality has to be understood in the sense of operators acting on $p\in \mathcal{P}_+(\R)$. This leads to the equation
	\begin{equation}\label{jumpingEquality}
		\frac{1}{2\pi|v|}\int_{\Pi^{\perp}(v)}\frac{p(v+\eta)}{|\eta|}d S(\eta)=\int_{\R}F(v|\xi)p(\xi) d \xi=\int_{\R}M(v,\xi)p(\xi)\delta_{\Gamma(\xi)}(v) d \xi.
	\end{equation}
	Choosing $p(\xi)=\delta(\xi-\xi_0)$ in \eqref{jumpingEquality} with $\xi_0=(1,0,0)$ gives
	\begin{equation}
		\frac{1}{2\pi|v|}\int_{\Pi^{\perp}(v)}\frac{\delta(v-\xi_0+\eta)}{|\eta|}d S(\eta)=M(v,\xi_0)\delta_{\Gamma(\xi_0)}(v).
	\end{equation}
	The distribution $\delta(v+\eta-\xi_0)$ is non-zero if and only if $v+\eta=\xi_0=(1,0,0)$ with $\eta \in \Pi_{\perp}(v)$ from which we have
	\begin{equation}\notag
		v=\left(\begin{array}{c}
			v_1\\
			v_2\\
			v_3
		\end{array} \right) \;\;
	\eta =\left(\begin{array}{c}
		1-v_1\\
		-v_2\\
		-v_3
	\end{array} \right)
	\end{equation}
	and, since $\eta \in \Pi_{\perp}(v)$, imposing the condition $v\cdot \eta=0$ yields
	\begin{equation}\label{Gamma}
		v_1^2+v_2^2+v_3^2=v_1.
	\end{equation}
	Since the surface $\Gamma(\xi)$ is invariant under rotations we can choose as $\xi_0$ the vector $\xi \in \R$ originally introduced in \eqref{pVk}. Therefore the surface $\Gamma(\xi)$ can be described as the set of points such that $|v|^2=v\cdot\xi$, or equivalently
	\begin{equation}
		0=|v|^2-v\cdot\xi=|v|^2-2v \cdot\frac{\xi}{2}+\frac{|\xi|^2}{4}-\frac{|\xi|^2}{4}=\Big|v-\frac{\xi}{2}\Big|^2-\frac{|\xi|^2}{4}.
	\end{equation}
	This shows that$\Gamma(\xi)=S_{\frac{|\xi|}{2}}\left(\frac{\xi}{2}\right)$ where $S_r(x)$ is the sphere centred at $x$ of radius $r$. Furthermore, since $v \cdot \eta=0$, we have $dS(\eta)=\frac{1}{|v|}\delta(v \cdot \eta)$ and assuming, without loss of generality, $v=|v|(1,0,0)$ in \eqref{jumpingEquality} we obtain
	\begin{gather}
		\frac{1}{2 \pi |v|}\int_{\Pi^{\perp}(v)}\frac{\delta(v+\eta-\xi)}{|\eta|}d S(\eta)
		=\frac{1}{2 \pi}\int_{\R}\frac{\delta(v+\eta-\xi)}{|\eta|}\delta(v\cdot \eta)d \eta. \notag
	\end{gather}
    The above equation has to be understood in the sense of probability measures $\mathcal{P}_+(\R)$ acting over test functions $\vp \in C^{\infty}_c(\R)$, therefore taking a test function $\vp$ we get
	\begin{align}\notag \label{intermediate1}
		\int_{\R}M(v,\xi)\varphi(v)\delta_{\Gamma(\xi)}(v)d v &= \int_{\Gamma(\xi)}\varphi(v)M(v,\xi)d v \\
		& = \frac{1}{2 \pi}\int_{\R}dv\int_{\R}\delta(v+\eta-\xi)\delta(\eta \cdot v)\frac{\varphi(v)}{|\eta|}d\eta.
	\end{align}
	Consider now the product of Dirac deltas $\delta(v+\eta-\xi)\delta(v \cdot \eta).$ From the equations $v+\eta-\xi=0, \; v\cdot \eta=0$
	it follows that $|v|^2-v\cdot\xi=0$ as in \eqref{Gamma}. Thus
	\begin{equation}
	    \delta(\eta \cdot v)=\delta((v-\xi+\eta)\cdot v-(v-\xi)\cdot v)=\delta(|v|^2-v\cdot \xi)
	\end{equation}
	which yields
	\begin{equation}\label{deltaProd}
		\delta(v+\eta-\xi)\delta(v \cdot \eta)=\delta(|v|^2-v\cdot\xi).
	\end{equation}
	Plugging \eqref{deltaProd} into \eqref{intermediate1} gives
	\begin{align}
		\int_{\R}M(v,\xi)\varphi(v)\delta_{\Gamma(\xi)}(v)d v  & = \frac{1}{2 \pi}\int_{\R}\delta(|v|^2-v\cdot \xi)\varphi(v)d v\int_{\R}\frac{\delta(v+\eta-\xi)}{|\eta|}d \eta  \notag \\
		& = \frac{1}{2 \pi}\int_{\R}\frac{1}{|v-\xi|}\varphi(v)\delta\left( \Big|v-\frac{\xi}{2}\Big|^2 -\frac{|\xi|^2}{4}\right) \frac{1}{2 \pi|\xi|}\int_{\R}\frac{1}{|v-\xi|}\varphi(v) \notag \\
        &= \frac{1}{2 \pi|\xi|}\int_{\Gamma(\xi)}\frac{1}{|v-\xi|}\varphi(v)
	\end{align}
	where we used the formula $\delta(x^2-a^2)=\frac{1}{2a}(\delta(x-a) +\delta(x+a))$. Finally \eqref{intermediate1} leads to
	\begin{equation}\label{M}
		M(v,\xi):=\frac{1}{2\pi |v-\xi|}\frac{1}{|\xi|}\delta\left( \Big|v-\frac{\xi}{2}\Big|^2 -\frac{|\xi|^2}{4}\right).
	\end{equation}
     In the end we have obtained the following formula for \eqref{eq:defpJump}, namely 
        \begin{align}
            p(v| \xi) &= \frac{1}{2\pi |\xi|} \frac{1}{|v-\xi|} \delta\left( \Big|v-\frac{\xi}{2} \Big|^2-\frac{|\xi|^2}{4} \right)=H(v, \xi)\delta\left(\Big |v-\frac{\xi}{2}\Big |-\frac{|\xi|}{2} \right);\label{eq:jumpProbFinal} \notag \\
            H(v,\xi) &:= \frac{1}{2\pi|\xi|} \frac{1}{|v-\xi|}.
        \end{align}
\end{proofof}

\subsubsection{Transition Probability}
    \begin{proofof}[Proof of Lemma \ref{lem:transitionp}] We now compute the transition probabilities for the stochastic process $(\xi^m)_m$. The key observation is that the process only consisting of $\{\xi^m\}_m$ is a Markov process and it is possible to compute the transition probability via the classical Chapman-Kolmogorov equation
\begin{equation}\label{transition1}
		p(\xi_{m+1}| \, \xi_m)=\int_{\R}p(\xi_{m+1}|\, v_m)p(v_m|\, \xi_m)d v_m.
	\end{equation}
    To simplify the formulas set $\xi_{m+1}=\xi,\xi_m=\eta, v_m=v$, then the Chapman-Kolmogorov equation \eqref{transition1} reads
	\begin{align}
		p(\xi|\, \eta) & = \int_{\R}p(\xi|\, v)p(v|\, \eta)dv \notag \\
		& = \int_{\R}G(\xi_1,v)H(v,\eta)\delta(\xi_2-v_2)\delta(\xi_3-v_3)\delta\left(\Big |v-\frac{\eta}{2}\Big |-\frac{|\eta|}{2} \right) dv \notag \\
		& = \int_{\R}\frac{G(\xi_1,v)}{2 \pi |v-\eta||\eta|}\delta(\xi_2-v_2)\delta(\xi_3-v_3)\delta\left(\Big |v-\frac{\eta}{2}\Big |-\frac{|\eta|}{2} \right)dv
	\end{align}
    where $G $ is as in \eqref{eq:defGprob}. 
	Integrating with respect to $(v_2,v_3)$ gives
	\begin{gather}\label{transitionInterme}
		\int_{\mathbb{R}}d v_1\frac{G(\xi_1,v_1,\xi_2,\xi_3)}{2 \pi |\eta|\sqrt{(v_1-\eta_1)^2+(\xi_2-\eta_2)^2+(\xi_3-\eta_3)^2}}\delta\left(\Big |\tilde v-\frac{\eta}{2}\Big |-\frac{|\eta|}{2} \right)
	\end{gather}	
    where we have set $\tilde v= (v_1,\xi_2,\xi_3)$. Moreover we have
	\begin{equation}\nonumber
		\delta\left(\Big |v-\frac{\eta}{2}\Big |-\frac{|\eta|}{2} \right)=\delta\left( \sqrt{\left( v_1-\frac{\eta_1}{2}\right) ^2+\left( \xi_2-\frac{\eta_2}{2}\right) ^2+\left( \xi_3-\frac{\eta_3}{2}\right) ^2}-\frac{|\eta|}{2}\right).
	\end{equation}
	We recall that for a smooth function $f$ with zeroes $\{x_j\}_{j \in J}$ we have that
        $$ \delta(f(x))= \sum_{j\in J}\frac{1}{|f'(x_j)|}\delta(x-x_j) . $$
    In particular, we have 
	\begin{align}
		&\delta\left( \sqrt{\left( v_1-\frac{\eta_1}{2}\right) ^2+\left( \xi_2-\frac{\eta_2}{2}\right) ^2+\left( \xi_3-\frac{\eta_3}{2}\right) ^2}-\frac{|\eta|}{2}\right) \notag \\
        &=\mathlarger{\mathlarger{\sum}}_{\ell=1,2}\frac{\sqrt{\left( v_1^{\ell}-\frac{\eta_1}{2}\right) ^2+\left( \xi_2-\frac{\eta_2}{2}\right) ^2+\left( \xi_3-\frac{\eta_3}{2}\right) ^2}}{\Big|v_1^{\ell}-\frac{\eta^{(1)}}{2} \Big|}\delta(v_1-v_1^{\ell}) = \mathlarger{\mathlarger{\sum}}_{\ell=1,2 }\frac{|\eta|/2}{\Big|v_1^{\ell}-\frac{\eta_1}{2} \Big|}\delta(v_1-v_1^{\ell})  \label{transitionDelta}
	\end{align}
	where we have used the fact that 
    \begin{equation}
        \sqrt{\left( v_1-\frac{\eta_1}{2}\right) ^2+\left( \xi_2-\frac{\eta_2}{2}\right) ^2+\left( \xi_3-\frac{\eta_3}{2}\right) ^2}-\frac{|\eta|}{2} =0
    \end{equation}
    if and only if $v_1$ satisfies
	\begin{equation}\label{zeroes}
		v_1^{\ell}=\frac{\eta_1}{2}\pm\sqrt{\frac{|\eta|^2}{4}-\left(\xi_2-\frac{\eta_2}{2} \right)^2-\left(\xi_3- \frac{\eta_3}{2}\right)^2}, \quad \ell \in \{1,2\}.
	\end{equation}
		Plugging \eqref{transitionDelta} into \eqref{transitionInterme} then leads to \eqref{transitionFinal}, i.e. 
	\begin{align}
		p(\xi|\eta) 
		= \frac{1}{\sqrt{\frac{|\eta|^2}{4}-\left(\xi_2-\frac{\eta_2}{2} \right)^2-\left(\xi_3- \frac{\eta_3}{2}\right)^2}} 
		\mathlarger{\sum}_{\ell=1,2}\frac{G(\xi_1,v_1^{\ell},\xi_2,\xi_3)}{4 \pi \sqrt{(v_1^{\ell}-\eta_1)^2+(\xi_2-\eta_2)^2+(\xi_3-\eta_3)^2}} 
	\end{align}
    and this concludes the proof.
\end{proofof}

 \bigskip


\section{Expected scenario in the case of frozen collisions, for homogeneity $\gamma <-1$}\label{sec:frozenCollision}

In this paper we focused on the analysis of solutions to \eqref{eq:Bshear} in the regime $\gamma\in (-1,0)$.
However, as noted in the introduction, the expected long-time behaviour of solutions to \eqref{eq:Bshear} for $\gamma<-1$ differs greatly from that described in the previous sections.  Here we  provide a heuristic justification for the emergence of these different long-time asymptotics when $\gamma < -1$. More precisely, we consider \eqref{eq:Bshear}, namely
\begin{align*}
	{\partial_{t}f + Kw_2  \partial_{w_1} f  }&  =\int_{\mathbb{R}^{3}}dw_{\ast}\int_{S^{2}%
	}d\omega |w-w_*|^\gamma
	\left(  M_*'g_{\ast}^{\prime}-M_*g\right), \, \quad f=f(w,t)\, .
\end{align*} 
Without loss of generality, we assume $K = 1$, as in the previous sections.
We then perform the following change of variables
\begin{equation}\label{eq3:CV2}
		f(w,t)=\frac{1}{t}F(\xi,t), \;\; \xi_1=\frac{w_1}{t}, \;\; w_2=\xi_2, \;\; w_3=\xi_3.
	\end{equation}
Plugging \eqref{eq3:CV2} into \eqref{eq:Bshear} and changing the time scale as $\tau=\log(t)$ we arrive at the following evolution equation
\begin{equation}\label{eq3:evF}
   \partial_{\tau}F {+} \partial_{\xi_1}((\xi_2 {-}\xi_1)F)=e^{\tau}\int_{\R}\int_{S^2}B(n \cdot \omega, |T(t)\xi-v_*|)(M_*'F'-M_*F)d \omega d v_* 
\end{equation}
with 
\begin{equation}
    T(t)=\left(\begin{array}{ccc}
        t & 0 & 0  \\
        0 & 1 & 0 \\
        0 & 0 & 1
    \end{array}\right), \quad S(t)=(T(t))^{-1}
\end{equation}
and where
\begin{equation}
    \begin{cases}
        \xi'= \xi-((T(t)\xi-v_*)\cdot \omega) S(t)\omega, & \\
        v_*' = v_*+ ((T(t)\xi-v_*) \cdot \omega) S(t) \omega.
    \end{cases}
\end{equation}
In the same spirit as the $\gamma\in (-1,0)$ case, we consider \eqref{eq:Bshear} in the large-velocity approximation. As discussed in Section \ref{ssec:AsymModel}, we use the adjoint equation \eqref{eq:homAdj} to derive the following approximate model  
\begin{equation}\label{eq:appr_FC}
    \partial_{t}\vp -w_2  \partial_{w_1} \vp  =|w|^{\gamma} \int_{S^2} (\varphi(w-(w \cdot \omega)\omega)-\varphi(w)) d \omega \, .
    \end{equation} 
Then, using the change of variables \eqref{eq3:CV2} in \eqref{eq:appr_FC}, we obtain 
\begin{equation}\label{eq3:evvp}
    \partial_{\tau}\psi - (\xi_2-{\xi_1}) \partial_{\xi_1}\psi=e^{\tau}|T(t)\psi|^{\gamma} \int_{S^2} (\psi(\xi-(T(t)\xi\cdot\omega)S(t)\omega)-\psi(\xi))d \omega.
\end{equation}
If $\gamma<-1$, the right-hand side of \eqref{eq3:evF} decays exponentially fast as $\tau \rightarrow + \infty$. Hence, for $\gamma<-1$, the collision term is expected to be negligible for long times; this regime is usually called the \emph{frozen collision} regime (see \cite{JNV1,JNV2,JNV3,MNV, NV}).  A rigorous proof of the smallness of this term requires a careful analysis of the interplay between collisions and shear and will be addressed in future work. Below we state the conjecture on the expected long-time behaviour of the solution.  

\begin{mainresult}[Asymptotics]\label{thm:ConvergenceFrozen}
    Let $f \in C([0,+\infty),\M)$ be a weak solution of \eqref{eq:Bshear}. 
    Then 
    \begin{equation}\label{eq:thm:ConvergenceFrozen}
       tf(t\,\xi_1,\xi_2,\xi_3,t) \overset{\ast}{\rightharpoonup} h(\xi_2,\xi_3)\delta(\xi_1-\xi_2) \quad \text{in $\M$ as $t \rightarrow + \infty.$} 
    \end{equation}
    where $h\in \mathscr{M}_{+}(\mathbb{R}^2)$ is a measure depending on the initial datum $f_0$ and $(\xi_1,\xi_2,\xi_3)= \left(\frac{v_1}{t},v_2,v_3\right)$. 
\end{mainresult}

We further remark that also in this case it is possible to obtain a probabilistic interpretation of \eqref{eq:appr_FC}, computing explicitly
the transition probabilities of the Markov process underlying \eqref{eq:appr_FC}. The results are analogous to the ones presented in Section \ref{sec:probabilities}. Specifically, Lemmas \ref{lem:pJump} and \ref{lem:transitionp} hold true also when $\gamma<-1$.  
The only difference occurs at the level of the flight probabilities of the process (cf. Lemma \, \ref{lem:pFlight}). In fact, when $\gamma<-1$, there is a non-zero probability that the velocity $v_m$ of the tagged particle could escape at infinity, i.e. $v_m=\infty$.  This turns out to be possible only when $\gamma<-1$.  More precisely, going back to the proof of Lemma \ref{lem:pFlight}, we observe that when $-1 < \gamma<0$ we have $\int_{\R}p(\xi|v)d \xi = 1$ where $p$ is the flight probability, i.e. \eqref{flightProb}. Conversely, when $\gamma<-1$ and $v_2 \neq 0$, then $\int_{\R}p(\xi|v)d \xi <1$ since there is a non-zero probability for a particle to escape at infinity, i.e. $\xi= \pm \infty$. 
The probabilities to reach $ \xi= \pm \infty$ are then given by
	\begin{align}
		p(+\infty|v) &= 1-\int_{v_1}^{+\infty}\; p(\xi_1|v)d \xi_1= \exp\left(- \frac{1}{|v_2|}\int_{v_1}^{+\infty}(\eta^2+v_2^2+v_3^2)^{-\frac{a}{2}}d \eta \right); \\
		p(-\infty|v) &= 1-\int_{-\infty}^{v_1}\; p(\xi_1|v)d \xi_1= \exp\left(- \frac{1}{|v_2|}\int_{-\infty}^{v_1}(\eta^2+v_2^2+v_3^2)^{-\frac{a}{2}}d \eta \right)\, ,
	\end{align}
    respectively. 
	The states $\xi=\pm \infty$ are absorbing states for the process. %
    Therefore, when $\gamma<-1$ we can define the flight probability as follows
\begin{equation}\label{eq:flightCompactified}
            \bar p(\xi|v)=
            \begin{cases}
                1 & \text{if $\xi=v=\infty$}, \\
                0 & \text{if $\xi \in \R, v=\infty$}, \\
                p(\xi|v) & \text{if $\xi,v \in \R$}, \\
                1-\int_{|v|}^{+\infty}\int_{|\rho|=1}p(R\rho|v)d\rho dR & \text{if $\xi=\infty, v \in \R$}. \\
            \end{cases}
        \end{equation}
        Notice that in the above definition we consider the $1-$point compactification of $\R$, namely $\Rc_c= \R \cup \{\infty\}$ and treat $\{\pm \infty\}$ as only one point $\{\infty\}$. With this definition, $\bar p(\xi|v)$ given as in \eqref{eq:flightCompactified} turns out to be a probability density.

\bigskip

\noindent\textbf{Acknowledgements.} 
A.~Nota  gratefully
acknowledges the support
by the project PRIN 2022 (Research Projects of National Relevance)
- Project code 202277WX43.   J.~J.~L.~Vel\'azquez gratefully acknowledges the support by the Deutsche Forschungsgemeinschaft (DFG)
through the collaborative research centre The mathematics of emerging effects (CRC 1060, Project-ID
211504053) and the DFG under Germany’s Excellence Strategy -EXC2047/1-390685813.
\bigskip




\bigskip

\def\adresse{
\begin{description}

\item[N. Miele:]{ Gran Sasso Science Institute,\\ Viale Francesco Crispi 7, 67100 L’Aquila, Italy  \\
E-mail: \texttt{nicola.miele@gssi.it}}

\item[A. Nota:] {Gran Sasso Science Institute,\\ Viale Francesco Crispi 7, 67100 L’Aquila, Italy  \\
E-mail: \texttt{alessia.nota@gssi.it}}

\item[J.~J.~L. Vel\'azquez:] { Institute for Applied Mathematics, University of Bonn, \\ Endenicher Allee 60, D-53115 Bonn, Germany\\
E-mail: \texttt{velazquez@iam.uni-bonn.de}}

\end{description}
}

\adresse


\begin{thebibliography}{99}        
    \bibitem{AS} M.~Abramowitz, I.~A.~Stegun, Handbook of mathematical functions, Dover editions, (1972)

     \bibitem {BLLS} H.~van Beijeren, O.~E.~Lanford, J.~L.~Lebowitz, H. Spohn, Equilibrium time correlation functions in the low-density limit, \textit{J. Stat. Phys.} \textbf{22}(2), 237--257 (1980)
    
     \bibitem{BNV} A.~Bobylev, A.~Nota, J.~J.~L.~Vel\'azquez, Self-similar Asymptotics for a Modified Maxwell–Boltzmann Equation in Systems Subject to Deformations, \textit{Comm. Math. Phys.} \textbf{380}, 409--448, (2020)


	\bibitem{CCDL}
    J.~A.~Carrillo, K.~ F.~Chan, R.~Duan, Z.~Li, 
    Inelastic Boltzmann equation under shear heating, \emph{Kinet. Relat. Model.}
    doi: 10.3934/krm.2026002	(2025)


	\bibitem {CercArchive}C.~Cercignani, Existence of homoenergetic affine flows for the Boltzmann equation, \textit{Arch. Rat. Mech. Anal.}
	\textbf{105}(4), 377--387, (1989)
	
	\bibitem {Cerc2000}C.~Cercignani, Shear Flow of a Granular Material,
	\textit{J. Stat. Phys.} \textbf{102}(5), 1407--1415, (2001)

    \bibitem{DJ}  K.~Dayal,  R.~D.~James, Nonequilibrium molecular dynamics for bulk materials and nanostructures, \textit{Journal of the Mechanics and Physics of Solids} \textbf{58}, 145–163, (2010) 

    \bibitem{DJ1} K.~Dayal, R.~D.~James, Design of a viscometers corresponding to a universal molecular simulation method, \textit{J. Fluid Mechanics}, \textbf{691}, 461--486,(2012)
    
    \bibitem{DL1} R.~Duan and S.~Liu, The Boltzmann equation for uniform shear flow, \emph{Arch. Rat. Mech. Anal.}, \textbf{242}(3), 1947--2002, (2021)

    \bibitem{DL2} R.~Duan and S.~Liu, Uniform shear flow via the Boltzmann equation with hard potentials, \emph{Discrete Contin. Dyn. Syst.} \textbf{45}(10),  4071--4118,  (2025) 


 
	\bibitem {Galkin1}V.~S.~Galkin, On a class of solutions of Grad's moment equation, \textit{PMM}, \textbf{22}(3), 386--389, (1958). 
	
	
	
	\bibitem {garzo} V.~Garz\'{o} and A.~Santos, Kinetic Theory of
	Gases in Shear Flows: Nonlinear Transport, {Kluwer Academic Publishers}, (2003)

    \bibitem{Kall} O.~Kallenberg, Foundations of Modern Probability, \emph{Probability Theory and Stochastic Modelling}, Springer Nature Switzerland, (2021)

    \bibitem {LNW} B.~Lods, A.~Nota, R.~Winter, Kinetic Description of a Rayleigh Gas with Annihilation, \textit{J. Stat. Phys.} \textbf{176} (6), 1434--1462, (2019)
	
	\bibitem {JNV1} R.~D.~James, A.~Nota, J.~J.~L.~Vel\'azquez, Self-similar profiles for homoenergetic solutions of the Boltzmann equation: particle velocity distribution and entropy, \emph{Arch. Rat. Mech. Anal.}, \textbf{231}, 787--843 (2019)
 
    \bibitem {JNV2}  R.~D.~James, A.~Nota, J.~J.~L.~Vel\'azquez, Long time asymptotics for homoenergetic solutions of the Boltzmann equation: Collision-dominated case, \emph{J. Nonlinear Sci.}, 29, 1943--1973, (2019)
    
    \bibitem {JNV3} R.~D.~James, A.~Nota, J.~J.~L.~Vel\'azquez, Long time asymptotics for homoenergetic solutions of the Boltzmann equation. Hyperbolic-dominated case, \emph{Nonlinearity}, 33, 3781--3815 (2020)  

    \bibitem{JQW} 
    R.~D.~James, K.~Qi, L.~Wang, On the kinetic description of the objective molecular dynamics, \textit{Multiscale Modeling \& Simulation}, \textbf{22}(4), 1646--1682, (2024)


    \bibitem {K1} B.~Kepka, Self-similar Profiles for Homoenergetic Solutions of the Boltzmann Equation for Non-cutoff Maxwell Molecules, \textit{J. Stat. Phys.} \textbf{190} (27) (2023)

     \bibitem {K2} B.~Kepka, Longtime behaviour of homoenergetic solutions in the collision dominated regime for hard potentials, \textit{Pure and Applied Analysis} \textbf{6}(2), 415--454 (2024)

    \bibitem{MT} K.~Matthies, F.~Theil, Rescaled Objective Solutions of Fokker-Planck and Boltzmann
equations, \emph{SIAM J. Math. Anal.}, \textbf{51}(2), 1321–1348. (2019)
    
    \bibitem{MST} K.~Matthies, G.~Stone, F.~Theil, The derivation of the linear Boltzmann equation from a Rayleigh gas particle model, \textit{Kinetic and Related Models}, \textbf{11}, 137--177 (2018)

    \bibitem{Miele} N.~Miele, PhD thesis, Gran Sasso Science Institute L'Aquila

    \bibitem{MNV} N.~Miele, A.~Nota,  J.~J.~L.~Vel\'{a}zquez, Homoenergetic solutions for the Rayleigh-Boltzmann equation: existence of a stationary non-equilibrium solution, \emph{J. Stat. Phys.} \textbf{192}, 97 (2025)
    
    \bibitem {NVW} A.~Nota, J.~J.~L.~ Vel\'azquez,  R.~Winter, Interacting particle systems with long-range interactions: scaling limits and kinetic equations, \textit{Atti Accad. Naz. Lincei Cl. Sci. Fis. Mat. Natur.} \textbf{32} (2), 335--377 (2021)

 \bibitem{NV} A.~Nota, J.~J.~L.~Vel\'azquez, Homoenergetic solutions of the Boltzmann equation: the case of simple-shear deformations, \textit{Mathematics in Engineering}, \textbf{5}(1), 1--25, (2023) 
    
    


    \bibitem{Spohn80} H.~Spohn, Kinetic equations from Hamiltonian dynamics: Markovian limits, \textit{Review of Modern Physics} \textbf{53}, 569--615 (1980)
	
	\bibitem {T}C.~Truesdell, On the pressures and flux of energy in a gas according to Maxwell's kinetic theory, II, \textit{J. Rat. Mech. Anal.} \textbf{5}, 55--128 (1956)
	
	\bibitem {TM}C.~Truesdell and 
    R.~G.~Muncaster, Fundamentals of Maxwell's
	Kinetic Theory of a Simple Monatomic Gas, {Academic Press}, (1980)
	
	\bibitem {V02}C.~Villani, A review of mathematical topics in collisional kinetic theory, {Hand-book of mathematical fluid dynamics}, \textbf{1}, 71--305, North-Holland, Amsterdam, (2002)

    

    
    
\end{thebibliography}
\end{document}